\documentclass{amsart}

\usepackage{amsfonts}
\usepackage{amsmath}
\usepackage{amssymb}
\usepackage{amsthm}
\usepackage{array}
\usepackage{enumerate}
\usepackage{fancyhdr}
\usepackage{fullpage}
\usepackage{graphicx}
\usepackage{hyperref}
\usepackage{ifthen}
\usepackage{longtable}
\usepackage{mathrsfs}
\usepackage{pgfplots}
\usepackage{tikz}

\usetikzlibrary{arrows.meta,backgrounds,calc,cd,decorations.pathmorphing,fit,patterns,positioning}

\pgfplotsset{compat=1.15}

\hypersetup{
    colorlinks=true,
    linkcolor=black,
    citecolor=black,
    filecolor=black,
    urlcolor=blue
}

\usepackage[
  backend=biber,
  hyperref=false,
  url=false,
  backref=false,
  isbn=false,
  doi=false,
  eprint=true,
  useprefix=true,
  maxcitenames=99,
  maxbibnames=99,  
  maxalphanames=99, 
  minalphanames=99,
  safeinputenc,
  style=alphabetic,
  citestyle=alphabetic,
  block=space,
  datamodel=preamble/ext-eprint
]{biblatex}
\title{Degree One Contributions and Open Gromov-Witten Invariants}
\date{}
\author{Sarah McConnell}

\counterwithin{figure}{section}

\newtheorem{theorem}{Theorem}[section]
\newtheorem{corollary}[theorem]{Corollary}
\newtheorem{proposition}[theorem]{Proposition}
\newtheorem{lemma}[theorem]{Lemma}
\newtheorem{hypothesis}[theorem]{Hypothesis}

\theoremstyle{definition}
\newtheorem{definition}[theorem]{Definition}
\newtheorem{example}[theorem]{Example}

\newtheorem{remark}[theorem]{Remark}

\newcommand{\nc}{\newcommand}

\nc{\C}{\mathbb{C}}
\nc{\E}{\mathbb{E}}
\renewcommand{\L}{\mathcal{L}}
\nc{\M}{\mathcal{M}}
\nc{\N}{\mathcal{N}}
\nc{\Q}{\mathbb{Q}}
\nc{\R}{\mathbb{R}}
\renewcommand{\S}{\mathcal{S}}
\nc{\X}{\mathcal{X}}
\nc{\Z}{\mathbb{Z}}

\nc{\Mbar}{\ov{\M}}
\nc{\Nbar}{\ov{\N}}

\renewcommand{\epsilon}{\varepsilon}
\nc{\ep}{\epsilon}

\nc{\arr}{\rightarrow}

\renewcommand{\hat}{\widehat}

\nc{\oo}{\infty}

\nc{\f}[2]{\frac{#1}{#2}}
\nc{\df}[2]{\dfrac{#1}{#2}}

\nc{\rt}[1]{\sqrt{#1}}

\nc{\ov}[1]{\overline{#1}}

\newcommand{\norm}[1]{\left\lVert#1\right\rVert}

\nc{\dlim}[1]{\lim\limits_{#1}}
\nc{\limto}[2]{\dlim{{#1} \arr {#2}}}

\nc{\sumto}[2]{\sum\limits_{#1}^{#2}}
\nc{\prodto}[2]{\prod\limits_{#1}^{#2}}

\nc{\dint}{\displaystyle\int}

\nc{\delbar}{\ov{\partial}}

\DeclareMathOperator{\Aut}{Aut}
\DeclareMathOperator{\codim}{codim}
\DeclareMathOperator{\coker}{coker}
\DeclareMathOperator{\Diff}{Diff}
\DeclareMathOperator{\im}{Im}
\DeclareMathOperator{\ind}{index}

\DeclareMathOperator{\Map}{Map}

\DeclareMathOperator{\proj}{Proj}
\DeclareMathOperator{\rk}{rk}

\nc\sphere{}
\def\sphere[#1](#2)(#3)(#4){
    \begin{scope}[shift={(#2)}]
        \def\xRadius{#3}
        \def\yRadius{#4}
	    \draw[#1] (-\xRadius,0) to [bend right=15] (\xRadius,0);
	    \draw[#1] [dashed] (-\xRadius,0) to [bend left=10] (\xRadius,0);
	    \draw[#1] (0,0) ellipse ({\xRadius} and {\yRadius});
    \end{scope} 
}

\nc\disk{}
\def\disk[#1](#2)(#3)(#4){
    \begin{scope}[shift={(#2)}]
        \def\xRadius{#3}
        \def\yRadius{#4}
	    \draw[#1] (-\xRadius,0) to [bend right=15] (\xRadius,0);
	    \draw[#1] [dashed] (-\xRadius,0) to [bend left=10] (\xRadius,0);
	    \draw[#1] (\xRadius,0) arc (0:180:{\xRadius} and {\yRadius});
    \end{scope} 
}

\nc\addGenus{}
\def\addGenus[#1](#2)(#3)(#4){
    \coordinate (genus) at (#2);
    \begin{scope}[scale=#4]
        \draw[#1] [rotate=#3] ($(genus)+(-1,0)$) to [bend left] ($(genus)+(1,0)$);
        \draw[#1] [rotate=#3] ($(genus)+(-1.2,0.1)$) to [bend right] ($(genus)+(1.2,0.1)$);
    \end{scope}
}

\nc\torus{}
\def\torus[#1](#2)(#3)(#4){
    \node at (#2) {$\bullet$};
    \def\torusXrad{0.6*#4}
    \begin{scope}[shift={(#2)}, rotate around={#3:(#2)}]
        \draw[#1] (0,#4) ellipse ({\torusXrad} and {#4});
        \addGenus[#1](0,#4)(90)(0.375*#4);
    \end{scope}
}

\nc\closedGhost{}
\def\closedGhost[#1](#2)(#3)(#4)(#5)(#6){
    \node at (#2) {$\bullet$};
    \begin{scope}[shift={(#2)}, rotate around={#3:(#2)}]
        \draw[#1] (0,#5) ellipse ({#4} and {#5});
        
	    \edef\genusDist{0}
	    \pgfmathparse{2*#5/(#6+1)} 
        \edef\genusDist{\pgfmathresult}

        \def\genusScale{0.5*#4}
        
        \foreach \idx in {1,...,#6}
            \addGenus[#1](0,\idx*\genusDist)(0)(\genusScale);
    \end{scope}
}

\nc\openGhostR{}
\def\openGhostR[#1](#2)(#3)(#4)(#5){
    \node at (#2) {$\bullet$};
    
    \begin{scope}[shift={(#2)}]
        \draw[#1] (2*#3,0) arc (0:180:{#3} and {#4});
        
        \edef\bdryRad{0}
    	\pgfmathparse{#3/(2*#5-1)} 
        \edef\bdryRad{\pgfmathresult}
        
        \foreach \idx in {1,...,#5}
            \draw [bend right=15] (4*\idx*\bdryRad-4*\bdryRad,0) to (4*\idx*\bdryRad-2*\bdryRad,0);
        \foreach \idx in {1,...,#5}
            \draw [bend left=10, dashed] (4*\idx*\bdryRad-4*\bdryRad,0) to (4*\idx*\bdryRad-2*\bdryRad,0);
        \ifthenelse{#5=1}{}{
            \foreach \idx in {2,...,#5}
                \draw (4*\idx*\bdryRad-4*\bdryRad,0) to [out=90, in=0] (4*\idx*\bdryRad-5*\bdryRad,0.3*#4) to [out=180, in=90] (4*\idx*\bdryRad-6*\bdryRad,0);
            }
    \end{scope}
}

\nc\openGhost{}
\def\openGhost[#1](#2)(#3)(#4)(#5)(#6){
    \ifthenelse{\equal{#6}{left}}{
        \begin{scope}[xscale=-1]
            \openGhostR[#1](#2)(#3)(#4)(#5)
        \end{scope}
    }{\openGhostR[#1](#2)(#3)(#4)(#5)}
}

\begin{document}

\begin{abstract}
We show that it is possible to define the contribution of degree one covers of a disk to open Gromov-Witten invariants. We build explicit sections of obstruction bundles in order to extend the algebro-geometric techniques of Pandharipande to the case of domains with boundary.
\end{abstract}

\maketitle

\thispagestyle{fancy} 

\tableofcontents

\section{Introduction} \label{intro}

In essence, Gromov-Witten invariants are counts of holomorphic curves. In the closed case (i.e., where curves do not have boundary), these invariants can be defined as integrals over moduli spaces of maps. These curve counts generally take values in $\Q$ (rather than $\Z$) for two reasons: first, that underlying domains may have automorphisms, and second, that some holomorphic maps factor through non-trivial branched covers of Riemann surfaces.

Spaces of domains are well understood: domains without boundary are described in \cite{harris} and \cite{mumford}, and domains with boundary are described in \cite{katzLiu} and \cite{liu}. This knowledge is enough to count simple maps in many cases by showing that they comprise an oriented manifold of dimension zero (see \cite{msBig}). However, the contributions of multiply covered curves are degenerate in the sense that moduli spaces of such maps are not of the expected dimension. One way to compute these contributions is via obstruction bundles. The fiber of the obstruction bundle over a moduli space of curves is the cokernel of the linearization of the $\delbar$ operator. If the rank of this bundle (or some relative of it) is equal to the dimension of the base, one may determine the contribution of the moduli space to Gromov-Witten invariants by computing the Euler number of the bundle (again taking values in $\Q$ because the spaces we consider may be orbifolds).

The link between Gromov-Witten invariants and obstruction bundles is based on Ruan-Tian perturbations (described in \cite{rt}). Rather than studying the holomorphic curve equation $\delbar_J(u)=0$, we study the perturbed equation $\delbar_J(u)=\nu$. The contribution of degree one covers of a curve $C$ (maps which are obtained by adding constant components to $C$) is precisely the count of those curves which perturb to a $t\nu$-holomorphic curve for all small $t$.

Open Gromov-Witten invariants are similar to their closed counterparts, but we allow domains to have boundary. For a thorough treatment of moduli spaces of open curves, see \cite{liu}. Two new issues arise when counting curves with boundary. The first problem is that of orientation. We do not address this rather thorny obstacle to defining invariants. Under suitable hypotheses (as in \cite{ortn}), the moduli space of open curves is orientable.

The second problem in the open case is that moduli spaces of domains may have codimension one boundary strata. In particular, the techniques used in \cite{pand} may fail because of this boundary. Indeed, suppose that $E$ is a vector bundle over an oriented manifold $X$. If $X$ is closed, then the Euler class of $E$ is the Poincar\'{e} dual of the zero locus of a generic section of $E$, and in the case $\rk(E)=\dim(X)$ this zero locus is a finite number points whose signed count is independent of the choice of section. This argument fails in the case where $X$ has boundary, as the number of zeros in a generic section may vary. The goal of this thesis is to solve this problem by adapting the techniques of \cite{pand}.

\subsection{Statement of Results} \label{resultsS}

We compute the contribution $C(g,h)$ of degree one covers of (regular, pseudoholomorphic, embedded) disks to Gromov-Witten invariants of type $(g,h)$ in a Calabi-Yau $3$-fold (cf. \cite{pand}).
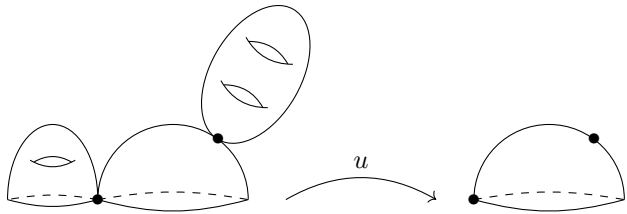
\begin{figure}[ht]
\centering
\begin{tikzpicture}

\def\r{1}
\def\h{1}
\def\w{0.6}
\def\d{5}

\coordinate (c2) at (0,0);
\coordinate (im) at (\d,0);

\disk[](c2)(\r)(\r)

\coordinate (LNode2) at ($(c2)+(-\r,0)$);
\openGhost[](LNode2)(\w)(\h)(1)(left);
\begin{scope}[shift={(LNode2)}]
	\addGenus[](-1*\w,0.5*\h)(0)(0.25);
\end{scope}

\coordinate (RNode2) at ($(c2)+(0.6*\r,0.8*\r)$);
\closedGhost[](RNode2)(-30)(\w)(\h)(2);


\coordinate (f) at ($0.5*(c2)+0.5*(im)$);
\draw [->] ($(f)+(-1,0)$) to [bend left] ($(f)+(1,0)$);
\node at ($(f)+(0,0.5)$) {$u$};


\disk[](im)(\r)(\r)

\coordinate (LNodeI) at ($(im)+(-\r,0)$);
\node at (LNodeI) {$\bullet$};

\coordinate (RNodeI) at ($(im)+(0.6*\r,0.8*\r)$);
\node at (RNodeI) {$\bullet$};

\end{tikzpicture}
\caption{A degree one cover of a disk.}
\end{figure}
We show that it is possible to define a contribution despite the codimension one boundary strata in the moduli space because it is always possible to construct a non-vanishing section of the obstruction bundle near these problematic strata. By relating the algebro-geometric techniques of \cite{pand} and \cite{niuZinger} to explicit sections of appropriate bundles, we avoid entirely the issue of defining characteristic classes of bundles over spaces with boundary.

We first prove a special case in Section~\ref{11s} in order to illustrate the main principles.

{
\renewcommand{\thetheorem}{\ref{calc11}}
\begin{theorem}
The contribution of degree one covers of a disk $\Sigma_0$ to type $(1,1)$ Gromov-Witten invariants is
\[
C(1,1) = \df{1}{2} \mu(T\Sigma_0,T\partial\Sigma_0) \cdot \chi(\E_1),
\]
where $\E_1$ is the Hodge bundle for genus $1$ curves.
\end{theorem}
\addtocounter{theorem}{-1}
}

In the remaining sections we consider the general case and prove the following theorems (stated slightly differently here for the sake of clarity; cf. Proposition~\ref{calc} and Corollaries~\ref{perm} and \ref{zero}).

\begin{theorem}
The contribution of degree one covers of a disk to Gromov-Witten invariants of genus $g$ with $h$ boundary components is zero whenever $h>1$.
\end{theorem}

\begin{theorem}
The contribution $C(g,1)$ of degree one covers of a disk to Gromov-Witten invariants of genus $g$ with $1$ boundary component is given by the generating function
\[
\sumto{g=0}{\oo} C(g,1)t^{2g-1} = \left( \df{\sin(t/2)}{t/2} \right)^{-1}.
\]
\end{theorem}

\begin{remark}
We expect that these techniques will extend to the case where the main component is not a disk. We predict a similar vanishing result: the contribution of degree one covers should be zero except when the domains of the cover and the main component have the same number of boundary components. The contribution in the case with the same number of boundary components will again reduce to the same flavor of generating function as the closed case.
\end{remark}

\subsection{Outline}

In Sections~\ref{prelim} and \ref{mainS} we include definitions, hypotheses, and standard results.

Section~\ref{11s} covers the simplest non-trivial case. The main principles of the argument all appear in this section, with minimal technical detail.

The remainder of this paper covers the general case. We examine moduli spaces $\Nbar$ of holomorphic maps in Section~\ref{baseS} and build an obstruction bundle $Ob$ in Section~\ref{obS}. In Sections~\ref{glueS} and \ref{lotS} we determine the relationship of this bundle to the contribution of these maps to Gromov-Witten invariants; in particular we describe a space $\L$ of gluing parameters and relate those maps which can be perturbed to a particular section $\alpha$ of $\pi_{\L}^*Ob$. Finally, we compute the contribution in Section~\ref{calcS}.
\[
\begin{tikzcd}
\Omega^{0,1}(\Nbar) \arrow[r, "\pi_{Ob}"] & Ob \arrow[d] & \pi_{\L}^*Ob \arrow[l, "\pi_{\L}"] \arrow[d]
\\
& \Nbar \arrow[ul, bend left, dashed, "\nu"] \arrow[u, bend left, dashed, "\ov{\nu}"] & \L \arrow[l, "\pi_{\L}"] \arrow[u, bend left, dashed, "\alpha"]
\end{tikzcd}
\]

The reader may wish to refer periodically to the list of symbols on page~\pageref{dict}.

The code used to generate figures is available at \url{https://www.overleaf.com/read/bdfwpzxsfqzc}.

\section{Preliminaries} \label{prelim}

In this section we give basic definitions and assumptions. Most are fairly standard; we include them for the sake of completeness. Definitions and hypotheses should be consistent with \cite{katzLiu}, \cite{liu}, and \cite{ortn} in order to ensure that moduli spaces of curves are sufficiently well-behaved.

We let $(M,\omega)$ be a symplectic manifold and $J$ an almost complex structure on $M$ which is tamed by $\omega$.

\begin{hypothesis} \label{hypL}
We assume that  $\dim_\R(M)=6$ and $c_1(M)=0$.

We assume that $J$ is generic and that all simple $J$-holomorphic maps (below some fixed energy bound) are regular embeddings with disjoint images.

We assume that $L$ is a spin Lagrangian submanifold of $M$ with Maslov class zero.
\end{hypothesis}

We endow $M$ with a metric so that a neighborhood of $L$ can be identified with $T^*L$ and $L$ is totally geodesic.

\begin{hypothesis} \label{hypMain}
We assume that $u_0:(\Sigma_0,\partial\Sigma_0)\arr(M,L)$ is an embedded disk. We assume that $u_0$ is $J$-holomorphic, that $u(\Sigma_0 \setminus \partial\Sigma_0) \cap L = \emptyset$, and that $u_0$ is regular in the sense that the linearization $D_0:\Gamma(\Sigma_0,\partial\Sigma_0;u^*TM,u^*TL)\arr\Omega^{0,1}(\Sigma,u^*TM)$ is surjective.
\end{hypothesis}

Throughout this paper we discuss degree one covers of maps satisfying Hypothesis~\ref{hypMain}.

\begin{definition}
For $U \subset \{z \in \C: \text{Im}(z) \geq 0\}$, a function $U \arr \C$ is holomorphic if it extends to a holomorphic function on an open neighborhood of $U$ in $\C$.
\end{definition}

\begin{proposition}[Schwarz Reflection Principle]
If $f:U \arr \C$ is holomorphic in the usual sense away from $U \cap \R$ and $f(U \cap \R) \subset \R$, then $f$ is holomorphic.
\end{proposition}

\begin{lemma} \label{holFnRS}
If $S$ is a compact Riemann surface (possibly with boundary) and $f:S \arr \C$ is holomorphic on $S$ with $f|_{\partial S} \subset \R$, then $f$ is constant.
\end{lemma}

\begin{definition} \label{modDomain}
A compact Riemann surface $\Sigma$ is \emph{closed} if $\partial\Sigma=\emptyset$ and \emph{open} otherwise.

For $g,n \in \mathbb{N}$, we denote by $\Mbar_{g,n}$ the moduli space of closed genus $g$ surfaces with $n$ marked points (see \cite{harris} and \cite{mumford}).

For $g,n \in \mathbb{N}$, $h \in \Z_+$, and $\vec{m} \in \mathbb{N}^h$, we denote by $\Mbar_{(g,h),n,\vec{m}}$ the moduli space of open surfaces of topological type $(g,h)$ with $(n,\vec{m})$ marked points. More specifically, 
\begin{itemize}
\item $g$ is the genus,
\item $h$ is the number of boundary components,
\item $n$ is the number of interior marked points, and
\item $\vec{m}=(m_1,\ldots,m_h)$, where $m_j$ is the number of marked points on the $j^{\text{th}}$ boundary component.
\end{itemize}
(For a complete definition of nodal curves with boundary, see Section~3 of \cite{liu}.)
\end{definition}

Observe that
\[
\dim_\R\Mbar_{(g,h),n,\vec{m}} = 3(2g+h-1)-3+2n+m_1+\ldots+m_h = \dim_\C\Mbar_{\tilde{g},2n+m_1+\ldots+m_h},
\]
where $\tilde{g}=2g+h-1$ is the genus of a doubled $(g,h)$-type curve.

The following definition is from Section 3.3.3 of \cite{katzLiu}.

\begin{definition} \label{cplxDouble}
For $\Sigma$ a bordered Riemann surface, the \emph{complex double} of $\Sigma$ is a closed Riemann surface $\Sigma^{(\C)}$ equipped with
\begin{enumerate}[(i)]
\item an antiholomorphic involution $c:\Sigma^{(\C)}\arr\Sigma^{(\C)}$,
\item a covering map $\pi:\Sigma^{(\C)} \arr \Sigma$ of degree two satisfying $\pi \circ c = \pi$, and
\item an embedding $\phi:\Sigma\arr\Sigma^{(\C)}$ such that $\pi\circ\phi=\text{I}_\Sigma$.
\end{enumerate}
The triple $(\Sigma^{(\C)},c,\pi)$ is unique up to isomorphism.
\end{definition}

Over a space $\mathcal{B}$ of domains with smooth maps into $M$, there is a vector bundle $\mathcal{E}$ whose fiber over $u:(\Sigma,\partial\Sigma)\arr(M,L)$ is
\[
\mathcal{E}_u = \Omega^{0,1}(\Sigma,u^*TM)
\]
(see Definition~\ref{defForm} and Section~3.1 of \cite{msBig}). 
We define a section $\delbar_J$ of this bundle by
\[
\delbar_J(\Sigma,u) = \df{1}{2}\left( du+J\circ du\circ j\right),
\]
where $j$ is the complex structure on $\Sigma$. The moduli space of $J$-holomorphic maps in $\mathcal{B}$ is the zero set of $\delbar_J$, up to automorphism (as in Section~2.1.3 of \cite{wendl}).

In order to understand this moduli space, we linearize $\delbar$. We first need to explain the tangent spaces to $\mathcal{B}$ and $\mathcal{E}$. If $\M$ is space of domains, we can view $T_\Sigma\M$ as variations in the complex structure on $\Sigma$. Aside from variations in the domain, we must consider variations in $u$, which are just vector fields along the image of $u$.

\begin{definition} \label{defVF}
For a smooth Riemann surface $\Sigma$ and a smooth map $u:(\Sigma,\partial\Sigma)\arr(M,L)$, a variation in $u$ is a section in
\[
\Gamma(\Sigma,\partial\Sigma;u^*TM,u^*TL) = \{\xi \in \Gamma(\Sigma,u^*TM): \xi(\partial\Sigma) \subset u^*TL\}.
\]
(If $\partial\Sigma=\emptyset$, then this space is just $\Gamma(\Sigma,u^*TM)$.)

For $\Sigma$ a nodal Riemann surface, label the smooth components $\Sigma_0,\ldots,\Sigma_r$. A variation in a smooth map $u:(\Sigma,\partial\Sigma)\arr(M,L)$ is a section $\xi=\xi_0\cup\ldots\cup\xi_r$, where $\xi_j$ is a variation in $u|_{\Sigma_j}$, such that these vector fields agree at the nodes (i.e., if two components are attached at a node $z \sim y$, then $\xi(z)=\xi(y) \in T_{u(z)}M$).
\end{definition}

\begin{definition} \label{defForm}
For $\Sigma$ nodal with smooth components $\Sigma_0,\ldots,\Sigma_r$ and a smooth map $u:(\Sigma,\partial\Sigma)\arr(M,L)$, we define $\Omega^{0,1}(\Sigma;u^*TM)$ to be the space of forms $\nu_0 \cup \ldots \cup \nu_r$, where $\nu_j$ is a $(0,1)$-form with values in $(u|_{\Sigma_j})^*TM$. We do not impose boundary conditions or require that the forms on the components match at the nodes.
\end{definition}

For a holomorphic map $u:(\Sigma,\partial\Sigma)\arr(M,L)$, where $\Sigma$ lies in some moduli space of domains $\M$, the linearization of $\delbar$ at $u$ is the map
\[
D:T_\Sigma\M \oplus \Gamma(\Sigma,\partial\Sigma;u^*TM,u^*TL) \arr\Omega^{0,1}(\Sigma,u^*TM)
\]
given by
\[
D(k,\xi) = \df{1}{2}\left( J \circ df \circ k + \nabla\xi + J(f) \circ \nabla\xi \circ j + (\nabla_\xi J) \circ df \circ j \right).
\]

\begin{theorem}[Riemann-Roch]
Assume that $S$ is a compact Riemann surface, $E \arr S$ a complex vector bundle with $F$ a totally real sub-bundle along $\partial S$, and $D$ a real linear Cauchy-Riemann operator on $E$. After completing in appropriate Sobolev norms (as Appendices~B and C of \cite{msBig}), $D$ is Fredholm with index
\[
\ind(D) = \rk_\C(E)\chi(S)+\mu(E,F).
\]
\end{theorem}

\begin{remark}
For thorough discussions of Sobolev spaces and their relevance to these arguments, see \cite{audinDamian}, \cite{dw}, \cite{liu}, and Appendices~B and C of \cite{msBig}.
\end{remark}

\begin{definition} \label{extTensor}
Given bundles $E \arr X$ and $F \arr Y$, the \emph{exterior tensor product} of $E$ and $F$ is
\[
E \boxtimes F = \pi_X^*(E) \otimes \pi_Y^*(F),
\]
where $\pi_X:X \times Y \arr X$ and $\pi_Y:X \times Y \arr Y$ are the natural projection maps.
\end{definition}

\begin{definition} \label{relTan}
Let $\phi_{g,n}:\mathcal{C}_{g,n} \arr \Mbar_{g,n}$ be the universal algebraic curve (see \cite{harris} and Appendix~D.6 of \cite{msBig}). 
The \emph{relative tangent bundle} $\mathcal{T}_{g,1}$ is $\ker(d\phi_{g,1})$; its fiber over $(\Sigma,z)$ is $T_z\Sigma$.

Fix $g \geq 0$, $h \geq 1$. Set $\vec{m}=(1,0,\ldots,0) \in \mathbb{N}^h$. The \emph{universal curve} $\phi_{(g,h),(0,\vec{m})}:\mathcal{C}_{(g,h),(0,\vec{m})} \arr \Mbar_{(g,h),(0,\vec{m})}$ is the real locus of the universal curve over $\Mbar_{\tilde{g},1}$ for $\tilde{g}=2g+h-1$. The \emph{relative tangent bundle} $\mathcal{T}_{(g,h),0,\vec{m}}$ is $\ker(d\phi_{(g,h),0,\vec{m}})$; its fiber over $(\Sigma,z)$ is $T_z\partial\Sigma$.
\end{definition}

\begin{definition} \label{hodge}
Fix $\Sigma \in \Mbar_{g,0}$ for $g \geq 0$. The \emph{$(0,0)$-Dolbeault operator} is
\[
\delbar_\Sigma:\Gamma(\Sigma,\C) \arr \Omega^{0,1}(\Sigma,\C).
\]
The \emph{Hodge bundle} $\E_g \arr \Mbar_{g,0}$ is the complex rank $g$ bundle of holomorphic differentials. Its dual $\E_g^*$ is the bundle whose fiber over $\Sigma$ is $\coker(\delbar_\Sigma)$.

Fix $\Sigma \in \Mbar_{(g,h),0,\vec{0}}$ for $g \geq 0$, $h \geq 1$. The \emph{$(0,0)$-Dolbeault operator} is
\[
\delbar_\Sigma:\Gamma(\Sigma,\partial\Sigma;\C,\R) \arr \Omega^{0,1}(\Sigma,\C).
\]
The \emph{Hodge bundle} $\E_{g,h}^* \arr \Mbar_{(g,h),0,\vec{0}}$ is the real rank $2g+h-1$ bundle which is the real locus of $\E_{\tilde{g}}$ for $\tilde{g}=2g+h-1$. Its dual $\E_{g,h}^*$ is the bundle whose fiber over $\Sigma$ is $\coker(\delbar_\Sigma)$.
\end{definition}

\begin{remark} \label{orbifold}
In general, $\Mbar_{g,n}$ is an orbifold and $\E_g^*$ an orbibundle. Similarly, $\Mbar_{(g,h),(n,\vec{m})}$ is an orbifold with boundary and $\E_{g,h}^*$ an orbibundle over it. However, these spaces have finite covers which are smooth, allowing us to ignore the orbifold structure.
\end{remark}

For discussion of the following proposition, see \cite{harris} or \cite{mumford}.

\begin{proposition} \label{eSquared}
The Euler classes of the duals of the Hodge bundles satisfy $e(\E_g^*)^2=0$ and $e(\E_{g,h}^*)^2=0$.
\end{proposition}

\section{Main Component and Linearization} \label{mainS}

In this section we establish some properties of linearizations. All of the holomorphic curves we consider will have a disk satisfying Hypothesis~\ref{hypMain} as a main component. The goal of this section is to compute the linearization for this disk and determine how the addition of constant components affects this linearization. These results are all standard; some readers may wish to skip this section entirely.

Let $C=u_0(\Sigma_0)$ and $T_0=TC_0$. We split 
\[
u_0^*TM = T_0 \oplus N_0
\]
and
\[
u_0^*TL = T_0^{(\R)} \oplus N_0^{(\R)},
\]
where $T_0^{(\R)}=T\partial C_0$ and $N_0^{(\R)}=N_0 \cap L$. Note that if we add ghosts to $u_0$, this splitting extends naturally to the trivial pullback bundle over any constant component. 

If $\mathcal{M}$ is a space of domains, we can then split the linearization into two pieces:
\begin{align*}
& D^T: \Gamma(\Sigma,TC) \oplus T_\Sigma\mathcal{M} \arr \Omega^{0,1}(\Sigma,TC)
\\
& D^N: \Gamma(\Sigma,N) \arr \Omega^{0,1}(\Sigma,N).
\end{align*}
(For both of these operators, we restrict in the domain and project in the codomain.)

\begin{lemma} \label{ker0}
Assume that $u_0:(\Sigma_0,\partial\Sigma_0)\arr(M,L)$ satisfies Hypothesis~\ref{hypMain} and let
\[
D_0~:~\Gamma(\Sigma_0,\partial\Sigma_0;u_0^*TM,u_0^*TL)\arr\Omega^{0,1}(\Sigma,u_0^*TM)
\]
be its linearization, with $D_0^N$ and $D_0^T$ its normal and tangent components. Then 
\[
\begin{array}{rclcrcl}
\ker(D_0^N) & = & 0 & \qquad\qquad & \coker(D_0^N) & = & 0
\\
\ker(D_0^T) & = & T_{\text{I}}\Aut(\Sigma_0) && \coker(D_0^T) & = & 0
\\
\ker(D_0) & = & T_{\text{I}}\Aut(\Sigma_0) && \coker(D_0) & = & 0.
\end{array}
\]
In particular, $D$ descends to an isomorphism $T_{u_0}(\mathcal{B}_0/Aut(\Sigma_0)) \arr \Omega^{0,1}(\Sigma_0,u_0^*TM)$ along the moduli space $\mathcal{B}_0$ of holomorphic disks, and the map $u_0$ is isolated in $\mathcal{B}_0$.
\begin{proof}
Let $C_0=u_0(\Sigma_0) \subset M$, $T_0=u_0^*TC_0$, and $N_0=u_0^*TM/u_0^*TC_0$. We split into tangent and normal directions:
\begin{align*}
\Gamma_0^T & = \Gamma(\Sigma_0,\partial\Sigma_0;T_0,T_0^{(\R)})
\\
\Gamma_0^N & = \Gamma(\Sigma_0,\partial\Sigma_0;N_0,N_0^{(\R)}).
\end{align*}

Let $\mathcal{B}_0=C^\oo(\Sigma_0,\partial\Sigma_0;M,L)$ (there are no variations in the complex structure on $\Sigma_0$). Around $u_0$, pick a slice $S \subset \mathcal{B}_0$ for the action of $\Aut(\Sigma_0)$:
\[
\begin{tikzcd}
& \mathcal{E} \arrow[d]
\\
S \arrow[r, hook] & \mathcal{B}_0 \arrow[d] \arrow[u, bend left, dashed, "\delbar"]
\\
& \mathcal{B}_0/\Aut(\Sigma_0) \arrow[ul, bend left, dashed, "\cong"]
\end{tikzcd}
\]
Near $u_0$, $S$ is isomorphic to the moduli space. 
If $O(u_0)$ is the orbit of $u_0$, then
\[
T_{u_0}S \cong T_{[u_0]}\left( \mathcal{B}_0/\Aut(\Sigma_0) \right) \cong T_{u_0}\mathcal{B}_0/T_{u_0}O(u_0).
\]
Note that $T_{u_0}O(u_0)$ lies entirely in the tangent direction. Therefore
\[
T_{u_0}S \cong (\Gamma_0^T/T_{u_0}O(u_0)) \oplus \Gamma_0^N.
\]
We define $D_0:T_{u_0}\mathcal{B}_0 \arr \Omega^{0,1}(u_0^*TM)$ in the usual way. 
Because $u_0$ is an embedding we have
\[
\mu(T_0,T_0^{(\R)})=\mu(T\Sigma_0,T\partial\Sigma_0)=2.
\]
Hypothesis~\ref{hypL} implies
\[
\mu(T_0,T_0^{(\R)})+\mu(N_0,N_0^{(\R)})=\mu(u_0^*TM,u_0^*TL)=0.
\]
Riemann-Roch gives
\begin{align*}
\ind(D_0^T)&=\rk_\C(T_0)\cdot\chi(\Sigma_0)+\mu(T_0,T_0^{(\R)})=3
\\
\ind(D_0^N)&=\rk_\C(N_0)\cdot\chi(\Sigma_0)+\mu(N_0,N_0^{(\R)})=0
\\
\ind(D_0)&=\rk_\C(u_0^*TM)\cdot\chi(\Sigma_0)+\mu(u_0^*TM,u_0^*TM^{(\R)})=3.
\end{align*}

First we address $D_0^T$. Because $u_0$ is an embedding, $D_0^T$ can be identified with the linearization of
\[
\delbar_{\Sigma_0}:\Diff(\Sigma_0) \arr \Omega^{0,1}(\Sigma_0,T\Sigma_0).
\]
Since the zero set is precisely $\Aut(\Sigma_0)$, we have $\ker(D_0^T)=T_{u_0}O(u_0)$ and $\coker(D_0^T)=0$.

By hypothesis $\coker(D_0^N)=0$, so $\ind(D_0^N)=0$ implies $\ker(D_0^N)=0$.

Finally we examine $D_0$, which is surjective by hypothesis. Therefore $\coker(D_0)=0$ and $\dim_\R(\ker(D_0))=3$. But since $T_{u_0}O(u_0) \subset \ker(D_0)$ also has dimension $3$ this inclusion must in fact be an equality.

It is clear from the computation that $D$ descends to an isomorphism after dividing out by equivalence, so the tangent space to the space of holomorphic disks at $u_0$ must be zero, meaning that $u_0$ is isolated.
\end{proof}
\end{lemma}

Linearizations of constant maps are straightforward; all that remains is to determine what happens when we add these constant maps to the main component. The proof of the following proposition is, unfortunately, somewhat technical in nature. Its purpose is merely to demonstrate that the kernel of the linearization is always the tangent space to the moduli space of curves and that the cokernel (i.e., the fiber of the obstruction bundle) depends only on the cokernels of the linearizations at the constant components. The moduli space of curves $\Nbar$ is defined in Section~\ref{baseS}, and the cokernels over the constant components appear in Section~\ref{obS}.

\begin{proposition} \label{nodalLin}
Assume that $u_0:(\Sigma_0,\partial\Sigma_0)\arr(M,L)$ satisfies Hypothesis~\ref{hypMain} and let $D_0$ be its linearization (as in Lemma~\ref{ker0}).

For $i=1,\ldots,r+q$, assume that $\Sigma_i \in \Mbar_{\sigma_i}$ is a domain with one marked point $y_i$ attached to $\Sigma_0$ at $z_i$. For $i \leq r$, we assume that $\Sigma_i$ is closed (i.e., $\sigma_i=(g_i,1)$) and attached at an interior point of $\Sigma_0$. For $i>r$, we assume that $\Sigma_i$ is open with marked point along the boundary (i.e., $\sigma_i=((g_i,h_i),0,(1,0,\ldots,0))$) and attached at a boundary point of $\Sigma_0$. 
Extend $u_0$ to a map $u$ defined on $\Sigma=\Sigma_0 \cup \Sigma_1 \ldots \cup \Sigma_{r+q}$ by requiring that $u$ be constant on all components but $\Sigma_0$. 
Let
\[
D_i:T_{\Sigma_i}\Mbar_{\sigma_i} \oplus \Gamma(\Sigma_i,\partial\Sigma_i;u_i^*TM,u_i^*TL)\arr\Omega^{0,1}(\Sigma_i,u_i^*TM)
\]
be the linearization of $\delbar$ at $u_i=u|_{\Sigma_i}$, with $D_i^N$ and $D_i^T$ its tangent and normal components. 

Let $(g,h)$ be the topological type of $\Sigma$ and $\S$ the (real codimension $2r+q$) nodal stratum of $\Mbar_{(g,h),0,\vec{0}}$ which contains $\Sigma$. Set
\[
\Nbar = \prod\limits_{i=1}^{r} (\Sigma_0 \times \Mbar_{\sigma_i}) \times \prod\limits_{i=r+1}^{r+q} (\partial\Sigma_0\times\Mbar_{\sigma_i}).
\]
If
\[
D:T_\Sigma\S \oplus \Gamma(\Sigma,\partial\Sigma;u^*TM,u^*TL)\arr\Omega^{0,1}(\Sigma,u^*TM)
\]
is the linearization of $\delbar$ at $u$, then
\[
\begin{array}{rclcrcl}
\ker(D^N) & = & 0 & \qquad\qquad & \coker(D^N) & = & \bigoplus\limits_{i=1}^{r+q} \coker(D_i^N)
\\
\ker(D^T) & = & T_\Sigma\Nbar && \coker(D^T) & = & \bigoplus\limits_{i=1}^{r+q} \coker(D_i^T)
\\
\ker(D) & = & T_\Sigma\Nbar && \coker(D) & = & \bigoplus\limits_{i=1}^{r+q} \coker(D_i).
\end{array}
\]
\begin{proof}
The analyses of the component linearizations appear in Lemmas~\ref{ker0}, \ref{kerClosed}, and \ref{kerOpen}. The primary goal of this proof is to determine how a collection of linearizations over components fit together into a linearization for a nodal Riemann surface. The key observation is that for $E$ a bundle over a nodal domain $\Sigma$, we cannot break $\Gamma(\Sigma,E)$ into a direct sum of factors of the form $\Gamma(\Sigma_i,E|_{\Sigma_i})$, but we can decompose in this manner if we restrict to sections which vanish at the nodal points. Using long exact sequences to relate these special sections to general sections, we can determine the relationship between the total linearization and its components (as in Appendix~A of \cite{splitting}).

Let $S$ be a domain with marked points ${\bf{x}}=(x_1,\ldots,x_k,x_{k+1},\ldots,x_{k+l})$ such that $x_i \in \partial S$ if and only if $i>k$. Let $E \arr S$ be a complex vector bundle with totally real sub-bundle $E^{(\R)}$ along $\partial S$. Set
\[
\Gamma_{\bf{x}}(S,\partial S;E,E^{(\R)}) = \{\xi \in \Gamma(S,\partial S;E,E^{(\R)}): \xi(x_i)=0 \text{ for all } i\}.
\]
Let
\[
E_{\bf{x}} = \left( \bigoplus\limits_{i=1}^{k} E_{x_i}  \right) \oplus \left( \bigoplus\limits_{i=k+1}^{k+l} E_{x_i}^{(\R)} \right)
\]
be the direct sum of the fibers over the marked points and $\text{ev}_{\bf{x}}$ the evaluation map at the marked points. For an operator $A:\Gamma(S,\partial S;E,E^{(\R)}) \arr \Omega^{0,1}(S,E)$, let $A_{\bf{x}}$ be its restriction to those sections which vanish at the marked points. Then there is a commutative diagram
\begin{equation}
\begin{tikzcd}
0 \arrow[r] & \Gamma_{\bf{x}}(S,\partial S;E,E^{(\R)}) \arrow[r] \arrow[d, "A_{\bf{x}}"] & \Gamma(S,\partial S;E,E^{(\R)}) \arrow[r, "\text{ev}_{\bf{x}}"] \arrow[d, "A"] & E_{\bf{x}} \arrow[r] & 0
\\
0 \arrow[r] & \Omega^{0,1}(S,E) \arrow[r, "\text{I}"] & \Omega^{0,1}(S,E) \arrow[r] & 0 \arrow[r] & 0
\end{tikzcd} \label{diagDecomp}
\end{equation}
with exact rows, allowing us to apply the snake lemma. 

We proceed with this construction and application of the snake lemma for each component of $\Sigma$ in both the normal and tangent directions. We begin with the computation in the normal direction, as it is fairly straightforward.

First, consider the main component $\Sigma_0$ with its marked points ${\bf z}=(z_1,\ldots,z_{r+q})$. We examine diagram~(\ref{diagDecomp}) with $S=\Sigma_0$, $E=N_0$, and $A=D_0^N$. By Lemma~\ref{ker0}, we have $\ker(D_0^N)=0$ and $\coker(D_0^N)=0$. Therefore the snake lemma yields the following exact sequence:
\begin{equation}
\begin{tikzcd}
0 \arrow[r] & \ker(D_{0,\bf{z}}^N) \arrow[r] & 0 \arrow[r] & (u_0^*NC)_{\bf{z}} \arrow[r, "\delta_0^N"] & \coker(D_{0,\bf{z}}^N) \arrow[r] & 0
\end{tikzcd}
\label{snake0N}
\end{equation}
It follows immediately that $\ker(D_{0,\bf{z}}^N)=0$ and $\coker(D_{0,\bf{z}}^N)=(u_0^*NC)_{\bf{z}}$.

Next we apply the same process to $D_i^N$ for $1 \leq i \leq r+q$. By Lemmas~\ref{kerClosed} and \ref{kerOpen}, $\ker(D_i^N)$ consists of constant sections (with values in $(u_i^*NC)_{y_i}$, where $(u_i^*NC)_{y_i}=N_{p_i}M$ for $i \leq r$ and $(u_i^*NC)_{y_i}=N_{p_i}L$ for $i>r$). Therefore when we restrict to those sections which vanish at a point, the kernel disappears: $\ker(D_{i,y_i}^N)=0$ for all $i$. The snake lemma gives the following exact sequence:
\begin{equation}
\begin{tikzcd}
0 \arrow[r] & (u_i^*NC)_{y_i} \arrow[r] &(u_i^*NC)_{y_i} \arrow[r, "\delta_i^N"] & \coker(D_{i,y_i}^N) \arrow[r] & \coker(D_i^N) \arrow[r] & 0
\end{tikzcd}
\label{snakeiN}
\end{equation}
It follows that $\delta_i^N=0$ and $\coker(D_{i,y_i}^N)=\coker(D_i^N)$.

As observed at the beginning of this proof, while the linearizations $D_i^N$ do not fit together in a straightforward fashion, their restrictions to sections which vanish at marked points do. That is, for any bundle $E$ we have
\[
\Gamma_{\bf{z},\bf{y}}(\Sigma,\partial\Sigma; E,E^{(\R)}) \cong \Gamma_{\bf{z}}(\Sigma_0,\partial\Sigma_0; E_0,E_0^{(\R)}) \oplus \bigoplus\limits_{i=1}^{r+q} \Gamma_{y_i}(\Sigma_i,\partial\Sigma_i; E_i,E_i^{(\R)}).
\]
Since $(0,1)$-forms are not required to match at nodes, $\Omega^{0,1}(\Sigma,E)$ also decomposes as a direct sum over components. Thus the operator $D_{\bf{z},\bf{y}}^N$ is precisely the direct sum of the operators $D_{0,\bf{z}}^N$ and $D_{i,y_i}^N$. It follows that
\[
\ker(D_{\bf{z},\bf{y}}^N) = \ker(D_{0,\bf{z}}^N) \oplus \ker(D_{1,y_1}^N) \oplus \ldots \oplus \ker(D_{r+q,y_{r+q}}^N) = 0
\]
and
\[
\coker(D_{\bf{z},\bf{y}}^N) = \coker(D_{0,\bf{z}}^N) \oplus \coker(D_{1,y_1}^N) \oplus \ldots \oplus \coker(D_{r+q,y_{r+q}}^N).
\]
As for $D^N$, the kernel is straightforward to analyze because if $D^N(\xi_0\cup\ldots\cup\xi_{r+q})=0$ then $D_i^N(\xi_i)=0$ for $0 \leq i \leq r+q$. Thus we see $\ker(D^N)=0$. Therefore applying the snake lemma to diagram~(\ref{diagDecomp}) with $S=\Sigma$, $E=u^*NC$, and $A=D^N$ yields the following exact sequence:
\begin{equation}
\begin{tikzcd}
0 \arrow[r] & (u^*NC)_{\bf{z}} \arrow[r, "\delta^N"] & \coker(D_{0,\bf{z}}^N) \oplus \bigoplus\limits_{i=1}^{r+q} \coker(D_{i,y_i}^N) \arrow[r] & \coker(D^N) \arrow[r] & 0
\end{tikzcd}
\label{snakeN}
\end{equation}
All that remains to be seen is that the image of $\delta^N$ is precisely $\coker(D_{0,\bf{z}}^N)$. Pick $v \in (u^*NC)_{\bf{z}}$ and choose some vector field $\xi=\xi_0\cup\ldots\cup\xi_{r+q}$ on $\Sigma$ so that $\text{ev}_{\bf{z}}(\xi)=v$. We can assume without loss of generality that $\xi$ is constant on $\Sigma_i$ for all $i>0$. Then
\[
\delta^N(v) \equiv D_0^N(\xi_0)\cup D_1^N(\xi_1)\cup\ldots\cup D_{r+q}^N(\xi_{r+q}) \equiv D_0^N(\xi_0)\cup 0 \cup\ldots\cup 0\pmod{Im(D_{\bf{z},\bf{y}}^N)},
\]
so $\delta^N$ certainly lands in $\coker(D_{0,\bf{z}}^N)$. But since $\coker(D_{0,\bf{z}}^N)=(u_0^*NC)_{\bf{z}}=(u^*NC)_{\bf{z}}$, we see that $\delta^N$ maps onto $\coker(D_{0,\bf{z}}^N)$. This proves
\[
\coker(D^N) = \bigoplus\limits_{i=1}^{r+q} \coker(D_{i,y_i}^N) = \bigoplus\limits_{i=1}^{r+q} \coker(D_i^N).
\]

The analysis of the operators in the tangent direction is similar in essence, but the argument is complicated by the addition of variations in the domain. 

Lemma~\ref{ker0} gives $\ker(D_0^T)=T_{\text{I}}\Aut(\Sigma_0)$ and $\coker(D_0^T)=0$. Restricting the kernel to those elements which vanish at marked points, we see $\ker(D_{0,\bf{z}}^T)=T_{\text{I}}\Aut(\Sigma_0,\bf{z})$. Therefore the snake lemma gives
\begin{equation}
\begin{tikzcd}
0 \arrow[r] & T_{\text{I}}\Aut(\Sigma_0,\bf{z}) \arrow[r] & T_{\text{I}}\Aut(\Sigma_0) \arrow[r] & (u_0^*TC)_{\bf{z}} \arrow[r, "\delta_0^T"] & \coker(D_{0,\bf{z}}^T) \arrow[r] & 0
\end{tikzcd}
\label{snake0T}
\end{equation}

Next we apply the same process to $D_i^T$ for $1 \leq i \leq r+q$. Since constant maps are holomorphic regardless of the domain, we see immediately that $T_{\Sigma_i}\Mbar_{\sigma_i}$ is contained in $\ker(D_i^T)$; we can therefore ignore variations in the domain. Let $\hat{D}_i^T$ be the linearization with fixed domain. By Lemmas~\ref{kerClosed} and \ref{kerOpen}, $\ker(\hat{D}_i^T)$ consists of constant sections (with values in $(u_i^*TC)_{y_i}$, where $(u_i^*TC)_{y_i}=T_{p_i}C$ for $i \leq r$ and $(u_i^*TC)_{y_i}=T_{p_i}\partial C$ for $i>r$). Therefore when we restrict to those sections which vanish at a point, the kernel disappears: $\ker(\hat{D}_{i,y_i}^T)=0$ and $\ker(D_{i,y_i}^T)=T_{\Sigma_i}\Mbar_{\sigma_i}$ for all $i$. The snake lemma gives the following exact sequence:
\begin{equation}
\begin{tikzcd}
0 \arrow[r] & (u_i^*TC)_{y_i} \arrow[r] & (u_i^*TC)_{y_i} \arrow[r, "\delta_i^T"] & \coker(D_{i,y_i}^T) \arrow[r] & \coker(D_i^T) \arrow[r] & 0
\end{tikzcd}
\label{snakeiT}
\end{equation}
It follows that $\delta_i^T=0$ and $\coker(D_{i,y_i}^T)=\coker(D_i^T)$.

Before analyzing $D^T$, we first examine the linearization $\hat{D}^T$ with fixed domain. As with $D^N$, we can split sections which vanish at nodal points and $(0,1)$-forms to see
\[
\ker(\hat{D}_{\bf{z},\bf{y}}^T) = \ker(D_{0,\bf{z}}^T) \oplus \ker(\hat{D}_{1,y_1}^T) \oplus \ldots \oplus \ker(\hat{D}_{r+q,y_{r+q}}^T) = T_{\text{I}}\Aut(\Sigma_0,\bf{z})
\]
and
\[
\coker(\hat{D}_{\bf{z},\bf{y}}^T) = \coker(D_{0,\bf{z}}^T) \oplus \coker(\hat{D}_{1,y_1}^T) \oplus \ldots \oplus \coker(\hat{D}_{r+q,y_{r+q}}^T).
\]
If $\hat{D}^T(\xi_0 \cup\ldots\cup \xi_{r+q})=0$ then $\hat{D}_i^T(\xi_i)=0$ for each $i$, so $\ker(\hat{D}^T)=\ker(D_0^T)=T_{\text{I}}\Aut(\Sigma_0)$. Therefore applying the snake lemma to diagram~(\ref{diagDecomp}) with $S=\Sigma$, $E=u^*TC$, and $A=\hat{D}^T$ yields the following exact sequence:
\begin{equation}
\begin{tikzcd}[column sep=small]
0 \arrow[r] & T_{\text{I}}\Aut(\Sigma_0,\bf{z}) \arrow[r] & T_{\text{I}}\Aut(\Sigma_0) \arrow[r] & (u^*TC)_{\bf{z}}
\\
& \arrow[r, "\delta^T"] & \coker(D_{0,\bf{z}}^T) \oplus \bigoplus\limits_{i=1}^{r+q} \coker(\hat{D}_{i,y_i}^T) \arrow[r] & \coker(\hat{D}^T) \arrow[r] & 0
\end{tikzcd}
\label{snakeT}
\end{equation}
But we can see (using an argument similar to that used for $D^N$) that the image of $\delta^T$ is precisely $\coker(D_{0,\bf{z}}^T)$, so the first half of the diagram is the same as sequence~(\ref{snake0T}). Therefore 
\[
\coker(\hat{D}^T) = \bigoplus\limits_{i=1}^{r+q} \coker(\hat{D}_{i,y_i}^T) = \bigoplus\limits_{i=1}^{r+q} \coker(D_i^T).
\]
Finally, we must determine how varying the domain changes the linearization. Note that we restrict to the nodal stratum $\S$ because we wish to examine only those maps with main component $(u_0,\Sigma_0)$. Moreover, varying the domain within $\S$ only changes the constant domains and the choices of marked points along the main component, so $T_\Sigma\S$ is contained in $\ker(D^T)$. Therefore $\ker(D^T)=\ker(\hat{D}^T) \oplus T_\Sigma\S$ and  $\coker(D^T)=\coker(\hat{D}^T)$. All that remains is to observe that the identification $\ker(D^T)=\ker(\hat{D}^T) \oplus T_\Sigma\S$ is really just a splitting of the following short exact sequence:
\[
\begin{tikzcd}
0 \arrow[r] & T_{\text{I}}\Aut(\Sigma_0) \arrow[r] & T_\Sigma\Nbar \arrow[r] & T_\Sigma\S \arrow[r] & 0
\end{tikzcd}
\]
Therefore $\ker(D^T)=T_\Sigma\Nbar$.

Finally, we observe that because $\ker(D^T)=0$ and $\ker(D^N)=0$, we must have $\ker(D)=0$ and hence $\coker(D)=\coker(D^T) \oplus \coker(D^N)$.
\end{proof}
\end{proposition}

\section{A Special Case: (1,1) Domains} \label{11s}

In this section we work through one example in detail. The results described here are all special cases of results in later sections. The purpose of this section is to illustrate the main principles in a situation which is less complicated than the general case.

Throughout this section, we assume that the domain $\Sigma$ has topological type $(g,h)=(1,1)$. As always, the main component $(\Sigma_0,u_0)$ satisfies Hypothesis~\ref{hypMain}. In Subsection~\ref{base11ss} we compute the moduli space $\Nbar$ of maps with such domains. 
We wish to compute the contribution $C(1,1)$ of maps in $\Nbar$ to type $(1,1)$ Gromov-Witten invariants. To do so, we fix a generic perturbation $\nu$, a section of the bundle of $(0,1)$ forms over $\Nbar$. We wish to count those maps $P \in \Nbar$ for which there exists a $t\nu$-holomorphic perturbation for every small $t$ (without loss of generality we can project $\nu$ to $Ob$, the cokernel bundle).

To count such maps, we will build obstruction and gluing bundles $Ob$ and $\L$ over $\Nbar$ in Subsections~\ref{ob11ss} and \ref{glue11ss}, respectively:
\[
\begin{tikzcd}
\Omega^{0,1}(\Nbar) \arrow[r, "\pi_{Ob}"] & Ob \arrow[d] & \pi_{\L}^*Ob \arrow[l, "\pi_{\L}"] \arrow[d]
\\
& \Nbar \arrow[ul, bend left, dashed, "\nu"] \arrow[u, bend left, dashed, "\ov{\nu}"] & \L \arrow[l, "\pi_{\L}"] \arrow[u, bend left, dashed, "\alpha"]
\end{tikzcd}
\]
We let $\ov{\nu}=\pi_{Ob}\circ\nu$. In Subsection~\ref{lot11ss} we build a linear section $\alpha$ of $\pi_\L^*Ob$ and let $Ob^F$ be the complement of its image $Im(\pi_\L \circ \alpha)$ in $Ob$. We will then show that a map $P\in\Nbar$ perturbs precisely when $\ov{\nu}(P) \in Im(\pi_\L \circ \alpha)$\footnote{There is also a technical requirement which we will show to be irrelevant: that $\ov{\nu}(P)$ land in the positive part of the fiber when the ghost is attached along $\partial\Sigma_0$.}. This will allow us to relate the count of perturbable maps to the zeros of the projection of $\ov{\nu}$ to $Ob^F$. Finally, in Subsection~\ref{calc11ss} we will demonstrate that it is possible to determine the zero count of a (particular type of) generic section of $Ob^F$ and compute this value.

\subsection{Moduli Space of (1,1) Domains} \label{base11ss}

In this section we compute the moduli space of holomorphic curves of type $(1,1)$ with main component $(u_0,\Sigma_0)$. This a special case of the analysis presented in Section~\ref{baseS}. We will show that the moduli space of such curves has two cells of real dimensions four and five, respectively, and that their intersection has real dimension three. These cells will exhibit rather different behaviors with respect to obstruction and gluing. For this reason we take care to explain precisely how these cells intersect; it will later allow us to reconcile these different types of behavior.

\begin{remark}
If $\Sigma$ is of type $(1,1)$, then its double $\Sigma^{(\C)}$ has ghosts with total genus $2$. These ghosts may be two conjugate tori or a single genus $2$ surface attached along the real locus.

However, not all symmetric closed curves with total genus $2$ are doubles of the open curves we consider. 
When we restrict ourselves to domains of type $(1,1)$, we exclude the case of a genus $2$ surface whose upper half has genus $0$ and three boundary components.
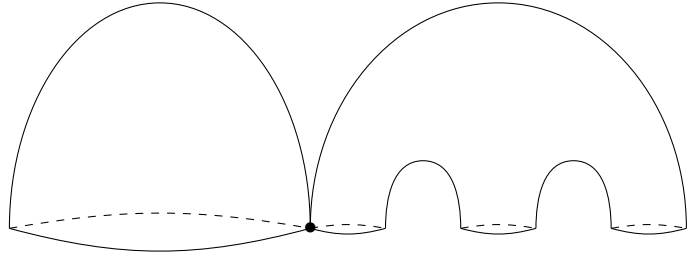
\begin{figure}[ht]
\centering
\begin{tikzpicture}

\disk[](-4,0)(2)(3);
\coordinate (node) at (-2,0);
\openGhost[](node)(2.5)(3)(3)();

\end{tikzpicture}
\caption{A curve of type $(0,3)$.}
\end{figure}
Such a domain is of type $(0,3)$; although its double is also of genus $2$, it cannot be obtained from those we consider in this section. The double of a $(1,1)$ curve can be deformed to the double of a $(0,3)$ curve, but such a deformation necessarily passes through curves which are not doubles of bordered Riemann surfaces. See \cite{katzLiu} for a more detailed discussion of open curves and their doubles.
\end{remark}

\begin{lemma} \label{11sigma}
There is a canonical inclusion $\rho:\Mbar_{1,1} \hookrightarrow \Mbar_\sigma$. Under this inclusion, $\E_{1,1}^*|_{\Mbar_{1,1}} = \E_1^*$, where $\E_{1,1}$ and $\E_1$ are the Hodge bundles for curves of type $(1,1)$ and genus $1$, respectively.
\begin{proof}
Each curve in $\Mbar_\sigma$ has a double in the real locus of $\Mbar_{2,1}$ (see Definition~\ref{cplxDouble}). There is a subset of the real locus of $\Mbar_{2,1}$ obtained by contracting two conjugate loops $\gamma,\ov{\gamma}$ as shown in Figure~\ref{m21}.
\begin{figure}[ht]
\centering
\begin{tikzpicture}

\coordinate (g2) at (0,0);
\node at ($(g2)+(-1.5,0)$) {$\bullet$};
\sphere[](g2)(1.5)(3);
\addGenus[]($(g2)+(0,2)$)(90)(0.5);
\addGenus[yscale=-1]($(g2)+(0,-2)$)(90)(0.5);

\draw [bend right=15] (-1.45,0.75) to (1.45,0.75);
\draw [bend left=10, dashed] (-1.45,0.75) to (1.45,0.75);
\draw [bend right=15] (-1.45,-0.75) to (1.45,-0.75);
\draw [bend left=10, dashed] (-1.45,-0.75) to (1.45,-0.75);
\node [right] at (1.45,0.75) {$\gamma$};
\node [right] at (1.45,-0.75) {$\ov{\gamma}$};

\coordinate (arrow) at ($(g2)+(3,0)$);
\draw [->, bend left] ($(arrow)+(-0.75,0)$) to ($(arrow)+(0.75,0)$);

\coordinate (nodal) at ($(g2)+(5.5,0)$);
\sphere[](nodal)(1)(1);
\node at ($(nodal)+(-1,0)$) {$\bullet$};

\coordinate (upnode) at ($(nodal)+(0,1)$);
\torus[](upnode)(0)(1.2);

\coordinate (downnode) at ($(nodal)+(0,-1)$);
\begin{scope}[yscale=-1]
    \torus[](downnode)(0)(1.2);
\end{scope}

\end{tikzpicture}
\caption{Embedding $\Mbar_{1,1}$ in $\Mbar_{(1,1),0,1}$.}
\label{m21}
\end{figure}
Such nodal symmetric genus two curves are in bijection with curves in $\Mbar_{1,1}$ because the central sphere lies in $\Mbar_{0,3}=\{\text{pt}\}$.

Take $\Sigma \in \Mbar_{1,1}$. For any compact Riemann surface $S$ we have
\[
\coker(\delbar_S)=H^{0,1}(S)=(H^{1,0}(S))^*.
\]
It follows that there is an injective map $\coker(\delbar_{\Sigma}) \arr \coker(\delbar_{\rho(\Sigma)})$. Since $\E_{1,1}$ has real rank $2(1)+(1)-1=2$ and $\E_1$ has complex rank $1$, this map must be an isomorphism.
\end{proof}
\end{lemma}

\begin{proposition} \label{base11}
Assume that $u_0:(\Sigma_0,\partial\Sigma_0) \arr (M,L)$ satisfies Hypothesis~\ref{hypMain}. Let $\Nbar$ be the moduli space of holomorphic curves $u:(\Sigma,\partial\Sigma)\arr(M,L)$ such that
\begin{enumerate}[(i)]
\item $\Sigma_0$ is a component of $\Sigma$ with $u|_{\Sigma_0}=u_0$,
\item $u$ is constant on every component other than $\Sigma_0$, and 
\item $\Sigma$ is of type $(1,1)$ (i.e., a smoothing has genus one with one boundary component).
\end{enumerate}
Then
\[
\Nbar = \Nbar_{1,1} \bigcup\limits_{\Nbar_{1,1} \cap \Nbar_\sigma} \Nbar_\sigma
\]
where
\begin{align*}
\Nbar_{1,1} & \cong \Sigma_0 \times \Mbar_{1,1}
\\
\Nbar_\sigma & \cong \partial\Sigma_0 \times \Mbar_\sigma
\\
\Nbar_{1,1} \cap \Nbar_\sigma & \cong \partial\Sigma_0 \times \Mbar_{1,1}
\end{align*}
for $\sigma=((1,1),0,1)$.
\begin{proof}
The domain $\Sigma$ of any curve in $\N$ has a complex double $\Sigma^{(\C)}$, with a sphere as a main component and ghosts with total genus $2$. This symmetric double has either two conjugate ghost tori or one symmetric genus $2$ ghost attached along the real locus. In the latter case, the ghost component of $\Sigma$ must be of type $(1,1)$ or $(0,3)$. Observe that a $(0,3)$ ghost cannot be perturbed to a $(1,1)$ ghost nor to a pair of genus $1$ ghosts, so we can exclude such domains from consideration. However, as we will see, pairs of ghost tori must be considered along with $(1,1)$ ghosts.

In the stratum where $\Sigma^{(\C)}$ has two ghost tori, $\Sigma_1 \in \Mbar_{1,1}$ is a closed torus, attached along the interior of $\Sigma_0$.
\begin{figure}[ht]
\centering
\begin{tikzpicture}

\disk[](0,0)(2)(3);
\coordinate (node) at (1.49,2);
\torus[](node)(-60)(2);

\end{tikzpicture}
\caption{Typical curve in $\Nbar_{1,1}$.}
\label{pic11}
\end{figure}
We cannot perturb $(u_0,\Sigma_0)$, so in the open stratum we can only vary the curve by varying $\Sigma_1 \in \Mbar_{1,1}$ or the point at which it is attached. Thus we obtain a stratum of the moduli space of the form $(\Sigma_0\setminus\partial\Sigma_0) \times \Mbar_{1,1}$ (which is of real dimension $4$).

Next we consider what happens when the marked point $z_1 \in \Sigma_0$ approaches the boundary. In $\Sigma^{(\C)}$, the ghost tori collide, producing a sphere bubble. Thus as $z_1 \arr \partial\Sigma_0 \cong S^1$, we see curves as in Figure~\ref{pic11sigma}.
\begin{figure}[ht]
\centering
\begin{tikzpicture}

\disk[](-4,0)(2)(3);
\node at (-2,0) {$\bullet$};
\disk[](0,0)(2)(3);

\coordinate (node) at (1.49,2);
\torus[](node)(-60)(2)

\draw [decorate,decoration={brace,amplitude=10pt,mirror},xshift=0.4pt,yshift=-0.4pt] (-2,-0.5) -- (5,-0.5) node[black,midway,yshift=-0.6cm] {constant};

\end{tikzpicture}
\caption{Typical curve in $\Nbar_{1,1} \cap \Nbar_\sigma$.}
\label{pic11sigma}
\end{figure}
The sphere bubble in $\Sigma^{(\C)}$ has three marked points, so there is only one such curve for each choice of node in $\partial\Sigma_0$ and each choice of ghost torus. Therefore $\Nbar$ has a cell\footnote{See Section~\ref{baseS} for terminology related to moduli spaces.} of the form $\Sigma_0 \times \Mbar_{1,1}$.

But we can also smooth the node between the two ghost components, yielding one ghost of type $(1,1)$, as in Figure~\ref{m11}.
\begin{figure}[ht]
\centering
\begin{tikzpicture}

\disk[](-4,0)(2)(3);

\openGhost[](-2,0)(2)(3)(1)();
\addGenus[](0,1.5)(30)(0.75);

\end{tikzpicture}
\caption{Typical curve in $\Nbar_\sigma$.}
\label{m11}
\end{figure}
This gives an additional cell $\partial\Sigma_0 \times \Mbar_\sigma$ of real dimension $5$, glued along the $3$-dimensional stratum $\partial\Sigma_0 \times \Mbar_{1,1}$. In particular, we glue via the inclusion of $\Mbar_{1,1}$ into $\Mbar_\sigma$ described in Lemma~\ref{11sigma}.
\end{proof}
\end{proposition}

\subsection{Obstruction Bundle for (1,1) Domains} \label{ob11ss}

We can build an obstruction bundle over the moduli space $\Nbar$ of curves, the fiber of which is the cokernel of the linearization (cf. Section~\ref{obS}). We begin by computing the fibers over the open strata of $\Nbar_{1,1}$ and $\Nbar_\sigma$ in the standard fashion. In the proof of Proposition~\ref{ob11} we examine the entire bundle over $\Nbar$.

\begin{lemma} \label{ob11torus}
Fix a map $u:(\Sigma,\partial\Sigma) \arr (M,L)$ in $\N_{1,1}$, so $\Sigma = \Sigma_0 \cup_{z \sim y} \Sigma_1$ for $\Sigma_0$ a disk, $(\Sigma_1,y) \in \M_{1,1}$, and $u_1=u|_{\Sigma_1}$ constant with value $p$. If
\[
\hat{D}:\Gamma(\Sigma,\partial\Sigma;u^*TM,u^*TL) \arr \Omega^{0,1}(\Sigma,u^*TM)
\]
is the linearization of the $\delbar$ operator at $u$, then $\ker(\hat{D})=0$ and
\[
\coker(\hat{D}) \cong T_pM \otimes_\C H^{0,1}(\Sigma_1).
\]
\begin{proof}
Let $\hat{D}_i$ be the linearization at $u$ over $\Sigma_i$, with no variations in the domain. We split into tangent and normal directions, as in Section~\ref{mainS}. For the analysis of $\hat{D}_0=D_0$, see Lemma~\ref{ker0}.

The linearization $D_1^V$ is straightforward for $V=u_1^*TC$ or $V=u_1^*NC$ (with $V_y$ the fiber of $V$) because $u_1$ is constant. Indeed, we can ignore the choice of domain because clearly the linearization is identically zero along $T_{\Sigma_1}\Mbar_{1,1}$. Because the bundle $u_1^*TM$ is trivial, the domain and codomain of $\hat{D}_1^V$ are
\begin{align*}
\Gamma(\Sigma_1,V) & \cong C^\oo(\Sigma_1, V_y)
\\
\Omega^{0,1}(\Sigma_1,V) & \cong \Omega^{0,1}(\Sigma_1) \otimes_\C V_y.
\end{align*}
Under this identification the linearization is just
\[
\hat{D}_1^V\xi_1 = \df{1}{2}\left( \nabla\xi_1+J(u_1) \circ \nabla\xi_1 \circ j \right) = \delbar\xi_1.
\]
It follows that $\ker(\hat{D}_1^V)$ is the set of holomorphic functions $\Sigma_1 \arr V_y$, all of which are constant by Lemma~\ref{holFnRS}. This proves $\ker(\hat{D}_1^V) \cong V_y$. We can also see that the cokernel is precisely $V_y \otimes_\C \coker_{\Sigma_1}(\delbar)$.

To complete our analysis of $\hat{D}$, all that remains is to determine how the components $\hat{D}_0$ and $\hat{D}_1$ fit together. This follows from an argument using long exact sequences; see Proposition~\ref{nodalLin}.
\end{proof}
\end{lemma}

\begin{lemma} \label{ob11open}
Fix a map $u:(\Sigma,\partial\Sigma) \arr (M,L)$ in $\N_\sigma$, so $\Sigma = \Sigma_0 \cup_{z \sim y} \Sigma_1$ for $\Sigma_0$ a disk, $z \in \partial\Sigma_0$, $(\Sigma_1,y) \in \M_\sigma$, and $u_1=u|_{\Sigma_1}$ constant with value $p \in L$. If
\[
\hat{D}:\Gamma(\Sigma,\partial\Sigma;u^*TM,u^*TL) \arr \Omega^{0,1}(\Sigma,u^*TM)
\]
is the linearization of the $\delbar$ operator at $u$\, then $\ker(\hat{D})=0$ and
\[
\coker(\hat{D}) \cong T_pL \otimes_\R H^{0,1}(\Sigma_1).
\]
\begin{proof}
Let $\hat{D}_i$ be the linearization at $u$ over $\Sigma_i$, with no variations in the domain. We split into tangent and normal directions, as in Section~\ref{mainS}. For the analysis of $\hat{D}_0=D_0$, see Lemma~\ref{ker0}.

Because $u_1$ is constant, the bundle pairs $(u_1^*TC,u_1^*T\partial C)$ and $(u_1^*NC,u_1^*N\partial C)$ are both trivial. If $(V,V^{(\R)})$ is either of these bundle pairs, then the linearization $\hat{D}_1^V$ is straightforward. Indeed, we can ignore variations in the domain because the linearization is identically zero along $T_{\Sigma_1}\Mbar_\sigma$. Because the bundle pair $(u_1^*TM,u_1^*TL)$ is trivial, the domain and codomain of $\hat{D}_1^V$ are
\begin{align*}
\Gamma(\Sigma_1, \partial\Sigma_1; V, V^{(\R)}) & \cong C^\oo(\Sigma_1,\partial\Sigma_1; V_y,V_y^{(\R)})
\\
\Omega^{0,1}(\Sigma_1,V) & \cong \Omega^{0,1}(\Sigma_1) \otimes_\R V_y^{(\R)}.
\end{align*}

Under this identification the linearization is just
\[
\hat{D}_1^V\xi_1 = \df{1}{2}\left( \nabla\xi_1+J(u_1) \circ \nabla\xi_1 \circ j \right) = \delbar\xi_1.
\]
It follows that $\ker(\hat{D}_1^V)$ is the set of holomorphic functions $(\Sigma_1,\partial\Sigma_1) \arr (V_y,V_y^{(\R)})$, all of which are constant by Lemma~\ref{holFnRS}. This proves $\ker(\hat{D}_1^V) \cong V_y^{(\R)}$. We can also see that the cokernel is precisely $V_y^{(R)} \otimes_\R \coker_{\Sigma_1}(\delbar)$.

To complete our analysis of $\hat{D}$, all that remains is to determine how the components $\hat{D}_0$ and $\hat{D}_1$ fit together. This follows from an argument using long exact sequences; see Proposition~\ref{nodalLin}.
\end{proof}
\end{lemma}

\begin{proposition} \label{ob11}
With $\Nbar$ as in Proposition~\ref{base11}, the obstruction bundle is of the form
\[
\begin{tikzcd}
u_0^*TM \boxtimes_\C \E_1^* \arrow[d] & \bigcup & u_0^*TL \boxtimes_\R \E_{1,1}^* \arrow[d]
\\
(\Sigma_0 \times \Mbar_{1,1}) & \bigcup\limits_{S^1 \times \Mbar_{1,1}} & (\partial\Sigma_0 \times \Mbar_\sigma)
\end{tikzcd}
\]
where $\E_1 \arr \Mbar_{1,1}$ and $\E_{1,1} \arr \Mbar_\sigma$ are the Hodge bundles. Over the intersection $\mathcal{N}_{1,1} \cap \mathcal{N}_\sigma \cong \partial\Sigma_0 \times \Mbar_{1,1}$, we identify the fibers as follows:
\begin{align*}
(u_0^*TM)|_{\partial\Sigma_0} \boxtimes_\C \E_1^* & \cong \left( u_0^*TL \otimes_\R \C \right) \boxtimes_\C \E_1^*
\cong u_0^*TL \boxtimes_\R \E_{1,1}^*.
\end{align*}
(These are isomorphisms of real vector spaces.)
\begin{proof}
By Lemmas~\ref{ob11torus} and \ref{ob11open}, the kernel of the linearization at any map $u \in \Nbar$ is zero\footnote{More precisely, the kernel of the linearization with fixed domain is zero. If we allow variations in the domain, the kernel is exactly the set of these variations.}. It follows that there exist bundles over $\Nbar_{1,1}$ and $\Nbar_\sigma$ whose fiber over $u$ is $\coker(D_u)$; these fibers were also computed in Lemmas~\ref{ob11torus} and \ref{ob11open}. It is clear how fibers fit together along $\N_{1,1}$ or $\N_\sigma$, so all that remains is to examine the intersection $\Nbar_{1,1} \cap \Nbar_\sigma$. This intersection sits inside $\Nbar_{1,1} \cong \Sigma_0 \times \Mbar_{1,1}$ as $\partial\Sigma_0 \times \Mbar_{1,1}$, and it sits inside $\Nbar_\sigma \cong \partial\Sigma_0 \times \Mbar_\sigma$ via the inclusion $\Mbar_{1,1} \hookrightarrow \Mbar_\sigma$.

Fix maps $u \in \Nbar_{1,1}$ and $v \in \Nbar_\sigma$ which represent the same curve in $\Nbar_{1,1} \cap \Nbar_\sigma$. Then $u$ and $v$ have domains $\Sigma_0 \cup_{z \sim y_{1,1}} \Sigma_{1,1}$ and $\Sigma_0 \cup_{z \sim y_\sigma} \Sigma_\sigma$, respectively, where $\Sigma_{1,1} \in \Mbar_{1,1}$ is identified with $\Sigma_\sigma$ under the inclusion $\Mbar_{1,1} \hookrightarrow \Mbar_\sigma$. Moreover, the restrictions $u_0$ and $v_0$ of $u$ and $v$ to the main component $\Sigma_0$ are equal, and both maps send their constant components to $p=u_0(z)$. 

Regardless of which space of curves we consider, Lemmas~\ref{ob11torus} and \ref{ob11open} show that the fiber of the obstruction bundle is the cokernel of the linearization over the ghost component. Let
\begin{align*}
D_{1,1}: & \Gamma(\Sigma_{1,1},\Sigma_{1,1} \times T_pM) \arr \Omega^{0,1}(\Sigma_{1,1},\Sigma_{1,1} \times T_pM)
\\
D_\sigma: & \Gamma(\Sigma_\sigma, \partial\Sigma_\sigma; \Sigma_\sigma \times T_pM, \Sigma_\sigma \times (T_pM)^{(\R)}) \arr \Omega^{0,1}(\Sigma_\sigma,\Sigma_\sigma \times T_pM)
\end{align*}
be these linearizations. But since these bundles are trivial, we can identify
\begin{align*}
\Gamma(\Sigma_{1,1},\Sigma_{1,1} \times T_pM) & \cong C^\oo(\Sigma_{1,1}, T_pM)
\\
\Omega^{0,1}(\Sigma_{1,1},\Sigma_{1,1} \times T_pM) & \cong \Omega^{0,1}(\Sigma_{1,1}) \otimes_\C T_pM
\\
\Gamma(\Sigma_\sigma, \partial\Sigma_\sigma; \Sigma_\sigma \times T_pM, \partial\Sigma_\sigma \times (T_pM)^{(\R)}) & \cong C^\oo(\Sigma_\sigma,\partial\Sigma_\sigma; T_pM, (T_pM)^{(\R)})
\\
\Omega^{0,1}(\Sigma_\sigma,\Sigma_\sigma \times T_pM) & \cong \Omega^{0,1}(\Sigma_\sigma) \otimes_\C T_pM.
\end{align*}
Under these identifications, $D_{1,1}$ and $D_\sigma$ just become $(0,0)$-Dolbeault operators for $\Sigma_{1,1}$ and $\Sigma_\sigma$, respectively. Since
\[
T_pM \cong \C \otimes_\R (T_pM)^{(\R)},
\]
all that remains is to apply Lemma~\ref{11sigma}:
\begin{align*}
T_pM \otimes_\C \coker(\delbar_{\Sigma_{1,1}}) & \cong (T_pM)^{(\R)} \otimes_\R \C \otimes_\C \coker(\delbar_{\Sigma_{1,1}})
\\
& \cong (T_pM)^{(\R)} \otimes_\R \coker(\delbar_{\Sigma_\sigma}).
\end{align*}
\end{proof}
\end{proposition}

\subsection{Gluing Parameters for (1,1) Domains} \label{glue11ss}

Our goal now is to determine the relationship between invariants of the bundle $Ob \arr \Nbar$ we built in Subsection~\ref{ob11ss} and the contribution of curves in $\Nbar$ to Gromov-Witten invariants. As described at the beginning of Section~\ref{11s}, we perturb the equation $\delbar(u) = 0$ via some $\nu$ and count those $P \in \Nbar$ which perturb to a $t\nu$-holomorphic map for all small $t$.

We would like to view the solution space as the zero locus of a generic section of a bundle. Ideally, we would be able to use $Ob \arr \Nbar$. Unfortunately, the rank of the bundle is too large: 
\[
\begin{array}{rcr}
\rk_\C(Ob_{1,1})	 = 3 & \qquad & \rk_\R(Ob_\sigma) = 6
\\
\dim_\C(\Nbar_{1,1}) = 2 && \dim_\R(\Nbar_\sigma) = 5.
\end{array}
\]
For this reason, we must introduce an extra line bundle (complex over $\Nbar_{1,1}$ and real over $\Nbar_\sigma$) in order to resolve this difference.

This bundle will consist of gluing parameters, which give ways to smooth out nodes to yield new (non-holomorphic) curves. In this subsection we construct this line bundle (cf. Section~\ref{glueS}); in Subsection~\ref{lot11ss} we will examine its relationship to the contribution we wish to compute.

\begin{definition}
The bundle $\L_{1,1}$ of gluing parameters over $\Nbar_{1,1}$ is $T\Sigma_0 \boxtimes_\C \mathcal{T}_{1,1}$, where $\mathcal{T}_{1,1}$ is the relative tangent bundle over $\Mbar_{1,1}$.

The bundle $\L_\sigma$ of gluing parameters over $\Nbar_\sigma$ is $T\partial\Sigma_0 \boxtimes_\R \mathcal{T}_\sigma$, where $\mathcal{T}_\sigma$ is the relative tangent bundle over $\Mbar_\sigma$. 
A real gluing parameter $\tau_z \otimes_\R \tau_y \in T_z\partial\Sigma_0 \otimes_\R T_y\partial\Sigma_1$ is \emph{positive} if $-j(\tau_z) \in T_z\Sigma_0$ and $j(\tau_y) \in T_y\Sigma_1$ are both inward pointing or both outward pointing.

The bundle of gluing parameters $\pi_\L:\L \arr \Nbar$ is obtained by attaching $\L_\sigma$ to $\L_{1,1}$ along $\Nbar_{1,1} \cap \Nbar_\sigma$. 
Over this intersection, the fibers of $\L_{1,1}$ are glued along $\Nbar_\sigma$ in a direction normal to $\Nbar_{1,1} \cap \Nbar_\sigma$ and the fibers of $\L_\sigma$ are glued along $\Nbar_{1,1}$ in a direction normal to $\Nbar_{1,1} \cap \Nbar_\sigma$ (see Remark~\ref{howToGlue11}).
\end{definition}

\begin{remark} \label{howToGlue11}
We see that $\L_{1,1}$ is a complex line bundle whose fiber over $(\Sigma_0 \bigcup_{z \sim y} \Sigma_1, u)$ is $T_z\Sigma_0 \otimes_\C T_y\Sigma_1$, and $\L_\sigma$ is a real line bundle whose fiber over $(\Sigma_0 \bigcup_{z \sim y} \Sigma_1, u)$ is $T_z\partial\Sigma_0 \otimes_\R T_y\partial\Sigma_1$.

Over the intersection $\Nbar_{1,1} \cap \Nbar_\sigma$, the direct sum $\L_{1,1} \oplus \L_\sigma$ has real rank three. 
Fix a map $(u,\Sigma) \in \Nbar_{1,1} \cap \Nbar_\sigma$. Its ghost component can be viewed as $[\Sigma_{1,1},y_{1,1}] \in \Mbar_{1,1}$ or $[\Sigma_\sigma,y_\sigma] \in \Mbar_\sigma$. 
Then the fiber of the direct sum over $u$ is $\C_u \oplus \R_u$, where
\begin{align*}
\C_u & = (T_z\Sigma_0 \otimes_\C T_{y_{1,1}}\Sigma_{1,1})
\\
\R_u & = (T_z\partial\Sigma_0 \otimes_\R T_{y_\sigma}\partial\Sigma_\sigma).
\end{align*}
We can split
\[
T\Nbar|_{\Nbar_{1,1} \cap \Nbar_\sigma} \cong T(\Nbar_{1,1} \cap \Nbar_\sigma) \oplus V_{1,1} \oplus V_\sigma,
\]
where $V_{1,1}$ is the normal bundle to $\Nbar_{1,1} \cap \Nbar_\sigma$ in $\Nbar_{1,1}$ and $V_\sigma$ is the normal bundle to $\Nbar_{1,1} \cap \Nbar_\sigma$ in $\Nbar_\sigma$. Observe that $\dim_\R(V_{1,1})=1$ and $\dim_\C(V_\sigma)=1$; we wish to identify these bundles with $\R_u$ and $\C_u$, respectively.

Lemma~\ref{glue11} gives instructions for identifying smoothing parameters with maps. A complex gluing parameter $\tau \in \C_u$ can be used to smooth the interior node of $\Sigma$, and a real gluing parameter $\tau \in \R_u$ can be used to smooth the boundary node of $\Sigma$ (see Figure~\ref{pic11sigma}). Thus we can identify 
\begin{align*}
T_u\Nbar_{1,1} & \cong T_u(\Nbar_{1,1} \cap \Nbar_\sigma) \oplus \R_u
\\
T_u\Nbar_\sigma & \cong T_u(\Nbar_{1,1} \cap \Nbar_\sigma) \oplus \C_u.
\end{align*}
Therefore it makes sense to identify the fibers of gluing parameters from the two pieces of the moduli space with directions normal to the intersection. 
This process allows us to build the bundle $\L$ of gluing parameters over all of $\Nbar$. Although $\Nbar$ has two pieces of different dimensions, the total space of $\L$ has constant real dimension $6$ (which is also the real rank of the obstruction bundle).
\end{remark}

While smoothing a given curve, we choose some small constant $R_0>0$ which satisfies all the hypotheses for gluing in \cite{dw}. In the case of an interior node $z \in \Sigma_0 \setminus \partial\Sigma_0$, we also add the hypothesis that $R_0$ is small enough to guarantee that the ball of radius $4R_0$ around $z$ does not intersect $\partial\Sigma_0$.

\begin{lemma} \label{glue11}
For $(u,\Sigma) \in \Nbar$, fix an element $\tau$ of the fiber $\L_{(u,\Sigma)}$ and assume
\begin{enumerate}[(i)]
\item $|\tau|<R_0$, and
\item $\tau$ is positive if $(u,\Sigma) \in \Nbar_\sigma$.
\end{enumerate}
Then $\tau$ yields a Riemann surface $(\Sigma_\tau,j_\tau)$ and a smooth map $\tilde{u}_\tau:(\Sigma_\tau,\partial\Sigma_\tau) \arr (M,L)$ such that $\norm{\delbar(\tilde{u}_\tau)}_{L^p}$ is small in the sense of Proposition~5.8 of \cite{dw}.
\begin{proof}
First consider $(u,\Sigma) \in \Nbar_{1,1} \setminus \Nbar_\sigma$. We can use $\tau$ to smooth the node $z \sim y$ and build a smooth map of this new Riemann surface into $M$ as in \cite{dw}. Because the domain and map are only altered near the node, the analysis in Sections~4.2 and 5.2 of \cite{dw} still applies.

However, we must take care when the node sits in $\partial\Sigma_0$. We may apply the results of \cite{dw} only to the double of the curve in the case of a boundary node.

Fix $(u,\Sigma) \in \Nbar_\sigma$ and let $(\Sigma^{(\C)},c,\pi)$ be the complex double of $\Sigma$. Choose a metric on $\Sigma^{(\C)}$ so that the fixed locus of the involution $c$ is totally geodesic.
\begin{figure}[ht]
\centering
\begin{tikzpicture}

\def\r{2}
\def\w{0.6*\r}

\coordinate (sigma0) at (0,0);
\sphere[](sigma0)(\r)(\r);
\node [left] at ($(sigma0)+(-\r,0)$) {Fix$(\sigma)$};

\node at ($(sigma0)+(\r,0)$) {$\bullet$};

\coordinate (sigma1) at ($(sigma0)+(\r+\w,0)$);
\sphere[](sigma1)(\w)(\r)

\begin{scope}[shift=(sigma1)]
	\addGenus[](0,0.375*\r)(20)(0.5);
	\addGenus[yscale=-1](0,-0.375*\r)(20)(0.5);
\end{scope}

\begin{scope}[shift=(sigma1)]
	\draw [decorate,decoration={brace,amplitude=10pt,mirror},xshift=0.4pt,yshift=-0.4pt] (\w+0.5,0) -- (\w+0.5,\r) node[black,midway,xshift=0.6cm] {$\Sigma$};
	\draw [decorate,decoration={brace,amplitude=10pt,mirror},xshift=0.4pt,yshift=-0.4pt] (\w+1.5,-\r) -- (\w+1.5,\r) node[black,midway,xshift=0.8cm] {$\Sigma^{(\C)}$};
\end{scope}

\end{tikzpicture}
\caption{The complex double of a curve in $\Nbar_\sigma$.}
\end{figure}
We must analyze smoothings of $\Sigma^{(\C)} = \Sigma_0^{(\C)} \bigcup_{z \sim y} \Sigma_1^{(\C)}$ to understand smoothings of $\Sigma$. Gluing parameters for $\Sigma^{(\C)}$ are (small) elements of $T_z\Sigma_0^{(\C)} \otimes_\C T_y\Sigma_1^{(\C)}$. We smooth $\Sigma^{(\C)}$ by removing small neighborhoods of $z$ and $y$ from $\Sigma_0^{\C}$ and $\Sigma_1^{\C}$, respectively, and then identifying small collars $A_z$ and $A_y$ around these removed neighborhoods via a map $\iota_\tau$. 
\begin{figure}[ht]
\centering
\begin{tikzpicture}

\def\r{2}
\def\w{0.6*\r}

\coordinate (sigma0) at (-1,0);

\begin{scope}[shift=(sigma0)]
	\draw (0,0) circle (\r);
	\draw [thick, bend right = 15] (-\r,0) to (\r,0);
\end{scope}

\begin{scope}[shift=(sigma0)]
	\draw [white, fill=white] (1.8,-0.872) to [bend left=20] (1.8,0.872) -- (2.1,1) -- (2.1,-1) -- cycle;
	
	\draw (1.8,-0.872) to [bend right=5] (1.8,0.872);
	\draw (1.8,-0.872) to [bend left=20] (1.8,0.872);
	\draw (1.6,-1.2) to [bend left=20] (1.6,1.2);
\end{scope}
\begin{scope}[shift=(sigma0)]
	\clip (sigma0) circle (\r);
	\draw [pattern = north east lines] (1.8,-0.872) to [bend left=20] (1.8,0.872) to [bend right=50] (1.6,1.2) to [bend right=20] (1.6,-1.2) to [bend right=50] cycle;
\end{scope}
\node [below right] at ($(sigma0)+(1.7,-0.872)$) {$A_z$};

\coordinate (sigma1) at (4,0);
\begin{scope}[shift=(sigma1)]
	\draw (0,0) ellipse ({\w} and {\r});
	\draw [thick, bend right = 15] (-\w,0) to (\w,0);
	\addGenus[](0,0.375*\r)(20)(0.5);
	\addGenus[yscale=-1](0,-0.375*\r)(20)(0.5);
\end{scope}

\begin{scope}[shift={(-\w+0.2,0)}]
    \draw [white, fill=white] (3.92,-0.872) to [bend right=20] (3.92,0.872) -- (3.7,1) -- (3.7,-1) -- cycle;
    
    \draw (3.92,-0.872) to [bend left=5] (3.92,0.872);
    \draw (3.92,-0.872) to [bend right=20] (3.92,0.872);
    \draw (4.04,-1.2) to [bend right=20] (4.04,1.2);
\end{scope}
\begin{scope}[shift={(-\w+0.2,0)}]
    \node [below left] at (3.92,-0.872) {$A_y$};
    \clip (sigma1) ellipse ({\w} and {\r});
    \draw [pattern = north east lines] (3.92,-0.872) to [bend right=20] (3.92,0.872) to [bend left=50] (4.04,1.2) to [bend left=20] (4.04,-1.2) to [bend left=50] cycle;
    \node [below left] at (3.92,-0.872) {$A_{y_i}$};
\end{scope}

\end{tikzpicture}
\caption{Collars near nodes.}
\label{collars11}
\end{figure}
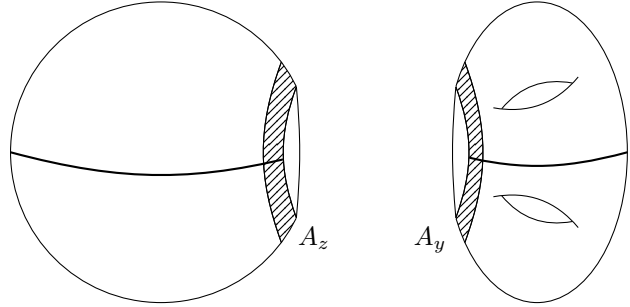

If $\tau=\tau_0 \otimes_\C \tau_1$ for $\tau_i$ tangent to $\Sigma_i$, then
\[
v \otimes_\C (\exp_y^{-1} \circ \iota_\tau \circ \exp_z(v)) = \tau_0 \otimes_\C \tau_1
\]
for all $v \in T_z\Sigma_0$. In particular, we have
\begin{align*}
\iota_\tau(\exp_z(t\tau_0))&=\exp_y(\tfrac{1}{t}\tau_1)
\\
\iota_\tau(\exp_z(-j(t\tau_0)))&=\exp_y(j(\tfrac{1}{t}\tau_1)).
\end{align*}
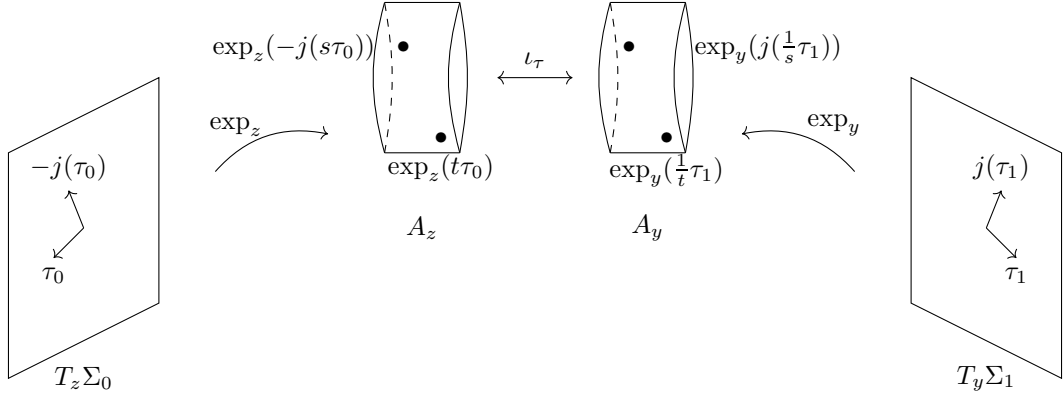
\begin{figure}[ht]
\centering
\begin{tikzpicture}

\def\d{1}
\def\w{1}
\def\r{3}
\def\h{2}

\draw (-\d,-1) to [bend right=10] (-\d,1);
\draw (-\d,-1) to [bend left=15] (-\d,1);
\draw [dashed] (-\d-\w,-1) to [bend right=10] (-\d-\w,1);
\draw (-\d-\w,-1) to [bend left=15] (-\d-\w,1);
\draw (-\d-\w,-1) -- (-\d,-1);
\draw (-\d-\w,1) -- (-\d,1);

\node at (-\d-0.75*\w,0.4) {$\bullet$};
\node [left] at (-\d-\w,0.4) {$\exp_z(-j(s\tau_0))$};
\node at (-\d-0.25*\w,-0.8) {$\bullet$};
\node [below] at (-\d-0.25*\w,-0.9) {$\exp_z(t\tau_0)$};
\node at (-\d-0.5*\w,-2) {$A_z$};

\coordinate (center0) at (-\r-3,-\h);
\coordinate (jt0) at ($(center0)+(-0.2,0.5)$);
\coordinate (t0) at ($(center0)+(-0.4,-0.4)$);
\draw [->] (center0) to (t0);
\draw [->] (center0) to (jt0);
\node [below] at (t0) {$\tau_0$};
\node [above] at (jt0) {$-j(\tau_0)$};
\draw ($(center0)+(-1,-2)$) -- ($(center0)+(-1,1)$) -- ($(center0)+(1,2)$) -- ($(center0)+(1,-1)$) -- cycle;
\node at ($(center0)+(0,-2)$) {$T_z\Sigma_0$};

\coordinate (arrow0) at ($0.5*(center0)+0.5*(-\d,0)$);
\draw [->] ($(arrow0)+(-0.75,-0.25)$) to [bend left] ($(arrow0)+(0.75,0.25)$);
\node [above left] at ($(arrow0)+(0,0.1)$) {$\exp_z$};

\draw [dashed] (\d,-1) to [bend right=10] (\d,1);
\draw (\d,-1) to [bend left=15] (\d,1);
\draw (\d+\w,-1) to [bend right=10] (\d+\w,1);
\draw (\d+\w,-1) to [bend left=15] (\d+\w,1);
\draw (\d+\w,-1) -- (\d,-1);
\draw (\d+\w,1) -- (\d,1);

\node at (\d+0.25*\w,0.4) {$\bullet$};
\node [right] at (\d+\w,0.4) {$\exp_y(j(\tfrac{1}{s}\tau_1))$};
\node at (\d+0.75*\w,-0.8) {$\bullet$};
\node [below] at (\d+0.75*\w,-0.9) {$\exp_y(\tfrac{1}{t}\tau_1)$};
\node at (\d+0.5*\w,-2) {$A_y$};

\coordinate (center1) at (\r+3,-\h);
\coordinate (jt1) at ($(center1)+(0.2,0.5)$);
\coordinate (t1) at ($(center1)+(0.4,-0.4)$);
\draw [->] (center1) to (t1);
\draw [->] (center1) to (jt1);
\node [below] at (t1) {$\tau_1$};
\node [above] at (jt1) {$j(\tau_1)$};
\draw ($(center1)+(1,-2)$) -- ($(center1)+(1,1)$) -- ($(center1)+(-1,2)$) -- ($(center1)+(-1,-1)$) -- cycle;
\node at ($(center1)+(0,-2)$) {$T_y\Sigma_1$};

\coordinate (arrow1) at ($0.5*(center1)+0.5*(\d,0)$);
\draw [->] ($(arrow1)+(0.75,-0.25)$) to [bend right] ($(arrow1)+(-0.75,0.25)$);
\node [above right] at ($(arrow1)+(0,0.1)$) {$\exp_y$};

\draw [<->] (-0.5*\d,0) -- (0.5*\d,0);
\node [above] at (0,0) {$\iota_\tau$};

\end{tikzpicture}
\caption{Identifying collars via $\tau$.}
\label{iota11}
\end{figure}
Now we must determine whether this smoothing of $\Sigma^{(\C)}$ yields a smoothing $\Sigma_\tau$ of $\Sigma$. The doubled curve $\Sigma^{(\C)}$ is equipped with an anti-holomorphic involution $c$ and a double cover $\pi:\Sigma^{(\C)}\arr\Sigma$. In order for the smoothing of $\Sigma^{(\C)}$ to yield a smoothing of $\Sigma$, the smoothing must respect these structures. That is, when we identify the collars in $\Sigma_0^{(\C)}$ and $\Sigma_1^{(\C)}$, the halves of the collars which lie in $\Sigma_0$ and $\Sigma_1$ must be identified. This occurs precisely when the gluing parameter $\tau$ lies in the positive half of the real locus of $T_z\Sigma_0^{(\C)} \otimes_\C T_y\Sigma_1^{(\C)}$.

Indeed, the smoothing $\Sigma_\tau^{\C}$ yields a smoothing of $\Sigma$ precisely when $\iota_\tau(A_z \cap \Sigma_0)=A_y \cap \Sigma_1$. If $\tau=\tau_0 \otimes_\C \tau_1$, we can assume without loss of generality that $\tau_1$ is tangent to the fixed locus $\text{Fix}(c)$ and that $j(\tau_1)$ points inward along $\Sigma_1$. It follows that $\exp_y(\tfrac{1}{t}\tau_1)$ must lie in the fixed locus and that $\exp_y(j(\tfrac{1}{t}\tau_1))$ must lie in $\Sigma_1$ (for appropriate values $t \in \R^+$). Since the points $\exp_y(\tfrac{1}{t}\tau_1)$ and $\exp_y(j(\tfrac{1}{t}\tau_1))$ are identified under $\iota_\tau$ with $\exp_z(t\tau_0)$ and $\exp_z(-j(t\tau_0))$, respectively, we can smooth $\Sigma$ via $\tau$ precisely when $\tau_0$ is also tangent to $\text{Fix}(c)$ and $-j(\tau_0)$ is inward pointing along $\Sigma_0$. When we embed $\Sigma \arr \Sigma^{(\C)}$, we see that $\text{Fix}(c)$ is precisely $\partial\Sigma$, so a smoothing of $\Sigma_\tau^{(\C)}$ yields a smoothing of $\Sigma$ if and only if the gluing parameter lies in the positive part of $T_z\partial\Sigma_0 \otimes_\R T_y\partial\Sigma_1$.

When $\tau$ is real and positive, we define a smoothed map $u_\tau:(\Sigma_\tau,\partial\Sigma_\tau)\arr(M,L)$ precisely as in \cite{dw}. This map still sends $\partial\Sigma$ to $L$ because $L$ is totally geodesic. It is evident from the construction that the estimates computed in \cite{dw} still apply.
\end{proof}
\end{lemma}

\subsection{Leading Order Term for (1,1) Domains} \label{lot11ss}

In order to relate gluing parameters to the perturbable maps we wish to count, we pull the obstruction bundle back over the map $\pi_{\L}:\L\arr \Nbar$ and build a section of this new bundle.

\begin{lemma} \label{lot11}
There is a section $\alpha:\L \arr \pi_{\L}^*Ob$ whose restriction to $\L_{1,1}$ is the leading order term of the obstruction map constructed in Section~5.7 of \cite{dw}. A curve $P$ perturbs to a $t\nu$-holomorphic map if and only if there exists a gluing parameter $\tau$ (which must be positive if the curve lies in $\mathcal{N}_\sigma$) such that $\alpha(P;\tau)=t\ov{\nu}_P$.
\begin{proof}
Fix $(\Sigma,u) \in \Nbar_{1,1}$. For $\tau=\tau_0 \otimes_\C \tau_1$, we define $\alpha_{1,1}(\Sigma,u;\tau)$ by
\[
\langle \alpha_{1,1}(\Sigma,u;\tau),v \otimes_\C \zeta \rangle_{L^2} = \langle (\ov{\zeta} \otimes_\C d_yu)(\tau_1),v \rangle.
\]
Similarly, if $(\Sigma,u) \in \Nbar_\sigma$ and $\tau=\tau_0 \otimes_\R \tau_1$, we define $\alpha_\sigma(\Sigma,u;\tau)$ by
\[
\langle \alpha_\sigma(\Sigma,u;\tau),v \otimes_\R \zeta \rangle_{L^2} = \langle (\ov{\zeta} \otimes_\R d_y(u|_{\partial\Sigma_0}))(\tau_1),v \rangle.
\]
See the proof of Lemma~\ref{lot} for further details.
\end{proof}
\end{lemma}

\subsection{Contribution for (1,1) Domains} \label{calc11ss}

\begin{proposition} \label{calc11}
The contribution of $\Nbar$ is
\begin{equation}
C(1,1) = \df{1}{2}\mu(N_0,N_0^{(\R)}) \cdot \chi(\E_1^*) = \df{1}{2}\mu(T\Sigma_0,T\partial\Sigma_0) \cdot \chi(\E_1). \label{cont11}
\end{equation}
\begin{proof}
What is written below applies if we first pass to a smooth cover of each orbifold. Since each space of domains has a finite smooth cover, we can ignore the orbifold structure entirely.

We will show that the contribution of $\mathcal{N}_\sigma$ is zero (cf. Proposition~\ref{calc}), which will leave only the contribution of $\mathcal{N}_{1,1}$. This second contribution can be computed in a fairly straightforward manner because the ghost is attached along the interior of the embedded component (and in particular, there is no issue of whether gluing parameters are positive).

We computed the obstruction bundle $Ob$ in Proposition~\ref{ob11} and bundle $\L$ of gluing parameters in Lemma~\ref{glue11}. Let $\alpha$ be the leading order term from Lemma~\ref{lot11}. Its image in $Ob_{1,1}$ is the complex line bundle $T\Sigma_0 \boxtimes_\C \E_1^*$, and its image in $Ob_\sigma$ is the real line bundle $T\partial\Sigma_0 \boxtimes_\R F_\sigma^\perp$, where $F_\sigma^\perp \subset \E_{1,1}^*$ is the (real rank $1$) complement of the bundle generated by $\zeta_y=0$. Let $Ob^F$ be the complement of the image of $\alpha$ in $Ob$:
\[
\begin{tikzcd}
N_0 \boxtimes_\C \E_1^* \arrow[d] & \bigcup & \left( (T_0^{(\R)} \boxtimes_\R F_\sigma) \oplus (N_0^{(\R)} \boxtimes_\R \E_{1,1}^*) \right) \arrow[d]
\\
(\Sigma_0 \times \Mbar_{1,1}) & \bigcup\limits_{S^1 \times \Mbar_{1,1}} & (\partial\Sigma_0 \times \Mbar_\sigma)
\end{tikzcd}
\]
Over each stratum, the rank of the bundle is equal to the dimension of the base, so a generic section has a finite number of zeros. Observe that in general a moduli space of open curves, such as $\Mbar_\sigma$, may have codimension one boundary (meaning that the zero count may vary from one section to the next). However, we will construct a non-vanishing section in order to demonstrate that the contribution of this cell of the moduli space is zero.

By Lemma~\ref{lot11}, we only need to construct a generic section $\rho$ of $Ob$ such that
\begin{enumerate}[(i)]
\item $\proj_{Ob_\sigma^F}(\rho_\sigma)$ is non-vanishing, and
\item the number of zeros of $\proj_{Ob_{1,1}^F}(\rho_{1,1})$ is (\ref{cont11}).
\end{enumerate}

Since the bundle over each stratum is a tensor product of bundles, we consider the factors separately. First we decompose the tangent bundles $TM|_{\Sigma_0}$ and $TL|_{\partial\Sigma_0}$. We can split into directions tangent and normal to the curve. Because the normal bundle is a complex rank two bundle over a surface, we can split off a trivial line bundle. Therefore, we can write
\begin{align*}
u_0^*TM & = V_1 \oplus V_2 \oplus V_3
\\
u_0^*TL & = V_1^{(\R)} \oplus V_2^{(\R)} \oplus V_3^{(\R)},
\end{align*}
where $V_1=T\Sigma_0$, $V_3$ is trivial, and $V_j^{(\R)}=V_j \cap L$. 

Pick generic sections $v_j$ of $V_j$ so that
\begin{enumerate}[(i)]
\item $v_3$ is non-vanishing,
\item $v_j|_{\partial\Sigma_0}$ lands in $V_j^{(\R)}$, and
\item $v_j|_{\partial\Sigma_0}$ is non-vanishing as a section of $V_j^{(\R)}$.
\end{enumerate}
It is possible to insist that the projections onto the real sub-bundles be non-vanishing because every (orientable) bundle over $\partial\Sigma_0 \cong S^1$ is trivial.

Next, choose sections $\eta_1,\eta_2,\eta_3$ of $\E_{1,1}^*$ which are transverse to the zero section so that
\begin{enumerate}[(i)]
\item $Z(\eta_2) \cap Z(\eta_3) = \emptyset$ and
\item $\proj_{F_\sigma}(\eta_1)$ is transverse to the zero section of $F_\sigma$.
\end{enumerate}
Note that by restricting to $\Mbar_{1,1}$ each $\eta_j$ yields a section of $\E_1^*$ satisfying the same properties (see Lemma~\ref{11sigma}).

Finally, we set
\[
\rho = (v_1 \boxtimes_\R \eta_1) \oplus (v_2 \boxtimes_\R \eta_2) \oplus (v_3 \boxtimes_\R \eta_3).
\]
The contribution from $\Nbar_\sigma$ is the signed count of positive zeros of
\[
\proj_{Ob_\sigma^F}(\rho_\sigma) = (v_1|_{\partial\Sigma_0} \boxtimes_\R \proj_{F_\sigma}(\eta_1)) \oplus (v_2|_{\partial\Sigma_0} \boxtimes_\R \eta_2) \oplus (v_3|_{\partial\Sigma_0} \boxtimes_\R \eta_3).
\]
It is this positivity criterion which makes counting difficult. 
But $v_2$ and $v_3$ are non-vanishing along $\partial\Sigma_0$ and $Z(\eta_2) \cap Z(\eta_3) = \emptyset$, implying that $\proj_{Ob_\sigma^F}(\rho_\sigma)$ has no zeros. In particular, the non-vanishing of $\proj_{Ob_\sigma^F}(\rho_\sigma)$ eliminates the issue of positivity of gluing parameters.

The total contribution from $\Nbar$ is the contribution from $\Nbar_{1,1}$, which is the signed count of zeros of
\[
\proj_{Ob_{1,1}^F}(\rho_{1,1}) = (v_2 \boxtimes_\C \eta_2|_{\Mbar_{1,1}}) \oplus (v_3 \boxtimes_\C \eta_3|_{\Mbar_{1,1}}).
\]
Since $v_3$ is non-vanishing and $\eta_2$ and $\eta_3$ have disjoint zero loci, the set of zeros is
\[
Z(v_2) \times Z(\eta_3).
\]
Because $\eta_3$ is a generic section of the complex line bundle $\E_1^* \arr \Mbar_{1,1}$, its zero locus represents the Euler class of this bundle. On the other hand, $v_2$ is a section over a disk, so we cannot use Chern classes to represent its zero locus. However, using the doubling constructions described in Section~3.3.3 of \cite{katzLiu}, we can see that $\#Z(v_2)$ is precisely half the Maslov index:
\[
\df{1}{2}\mu(V_2,V_2^{(\R)}) = \df{1}{2}\mu(N_0,N_0^{(\R)}) = -\df{1}{2}\mu(T\Sigma_0,T\partial\Sigma_0).
\]
This completes the proof.
\end{proof}
\end{proposition}

\begin{remark}
Let $V_2 \arr \Sigma_0$ and $v_2$ be as in the proof of Proposition~\ref{calc11}. Assume the target manifold $M$ has an anti-symplectic involution whose fixed locus is $L$. We can double bundles and sections, as in Section~3.3.3 of \cite{katzLiu}. Therefore the contribution is
\[
\df{1}{2}\mu(V_2,V_2^{(\R)})\chi(\E_1^*) = \#Z(v_2)\chi(\E_1^*) = \df{1}{2}c_1(V_2^{(\C)})\chi(\E_1^*) = \df{1}{2}c_1(N_0^{(\C)})\chi(\E_1^*)
\]
(cf. \cite{niuZinger}).
\end{remark}

\section{Moduli Spaces of Curves} \label{baseS}

Fix $g \in \N$ and $h \in \Z_+$. In the remaining sections we examine the moduli space of holomorphic curves with main component $(u_0,\Sigma_0)$ and topological type $(g,h)$. 
In each case the domain is modeled on a tree with root $\Sigma_0$ and (possibly nodal) constant branches. These curves are indexed by the distribution of genus and boundary components across ghost branches.

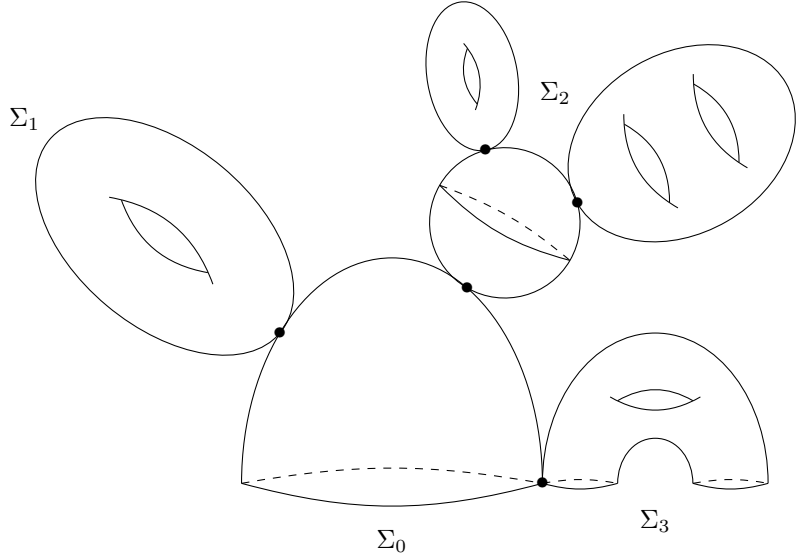
\begin{figure}[ht]
\centering
\begin{tikzpicture}

\coordinate (sigma0) at (0,0);
\def\rzero{2}
\disk[](sigma0)(\rzero)(1.5*\rzero);
\begin{scope}[shift=(sigma0)]
	\node [below] at (0,-0.5) {$\Sigma_0$};
\end{scope}

\coordinate (node1) at ($(sigma0)+(-1.49,2)$);
\def\rone{2}
\torus[](node1)(50)(\rone)
\begin{scope}[shift=(node1), rotate around={50:(node1)}]
    \node [above left] at ($(node1)+(0,2*\rone)$) {$\Sigma_1$};
\end{scope}

\coordinate (node2) at ($(sigma0)+(1,2.6)$);
\def\rtwo{1}

\begin{scope}[shift=(node2), rotate around={-30:(node2)}]
    \coordinate (node2a) at (-0.707*\rtwo,1.707*\rtwo);
    \coordinate (node2b) at (0.707*\rtwo,1.707*\rtwo);
    \sphere[](0,\rtwo)(\rtwo)(\rtwo);
    \node at (0,0) {$\bullet$};
    \torus[](node2a)(40)(1);
    \closedGhost[](node2b)(-30)(1.2)(1.6)(2);
    \node [left] at (0,3) {$\Sigma_2$};
\end{scope}

\coordinate (node3) at ($(sigma0)+(\rzero,0)$);
\openGhost[](node3)(1.5)(2)(2)();
\begin{scope}[shift=(node3)]
	\addGenus[](1.5,1.1)(0)(0.5)
	\node at (1.5,-0.5) {$\Sigma_3$};
\end{scope}

\end{tikzpicture}
\caption{A curve with three ghost branches, modeled on $(1,3,(1,2))$.}
\end{figure}

\begin{definition} \label{partNum}
Fix $g \in \mathbb{N}$ and $h \in \Z_+$. A \emph{partition} of $(g,h)$ is an ordered choice of $\{g_1,\ldots,g_r\}$ and $\{(g_{r+1},h_{r+1}),\ldots,(g_{r+q},h_{r+q})\}$ for some $r,q \geq 0$ so that
\begin{align*}
g & = 0+g_1+\ldots+g_k
\\
h & = 1+(h_{r+1}-1)+\ldots+(h_k-1).
\end{align*}
We require
\begin{enumerate}[(i)]
\item $g_i \geq 1$ for $i \leq r$ and
\item $g_i \geq 0$, $h_i \geq 1$, and $2g_i+h_i-1 \geq 1$ for $i>r$.
\end{enumerate}
\end{definition}

\begin{remark}
Two partitions are equivalent if they are the same up to re-ordering. However, we will ignore this equivalence until Section~\ref{calcS}. We will order partitions so that all the closed ghosts appear before the open ghosts for the sake of notational clarity.
\end{remark}

\begin{definition} \label{partition}
Fix a partition $\lambda=(g_1,\ldots,g_r,(g_{r+1},h_{r+1}),\ldots,(g_{r+q},h_{r+q}))$ of some topological type $(g,h)$. A \emph{domain modeled on $\lambda$} is a nodal domain $\Sigma \in \Mbar_{(g,h),0,\vec{0}}$ so that
\begin{enumerate}[(i)]
\item $\Sigma_0$ is disk with marked points $\{z_1,\ldots,z_{r+q}\}$ so that $z_i \in \partial\Sigma_0$ if and only if $i>r$,
\item for $1 \leq i \leq r$, $\Sigma_i \in \Mbar_{g_i,1}$ is a closed curve with marked point $y_i$,
\item for $r+1 \leq i \leq r+q$, $\Sigma_i \in \Mbar_{(g_i,h_i),0,(1,0,\ldots,0)}$ is an open curve with marked point $y_i \in \partial\Sigma_i$, and
\item $\Sigma_i$ is attached to $\Sigma_0$ by identifying marked points:
\[
\Sigma=\left.\left( \bigsqcup\limits_{i=0}^{r+q} \Sigma_i \right) \middle/ (z_i \sim y_i) \right.
\]
\end{enumerate}
We refer to $\Sigma_0$ as the \emph{main component}; the remaining (possibly nodal) curves $\Sigma_1,\ldots,\Sigma_{r+q}$ are \emph{branches}.

A \emph{holomorphic curve modeled on $\lambda$} is a holomorphic map $u:(\Sigma,\partial\Sigma)\arr(M,L)$ so that
\begin{enumerate}[(i)]
\item $\Sigma$ is a domain modeled on $\lambda$,
\item $u_0=u|_{\Sigma_0}$ satisfies Hypothesis~\ref{hypMain}, and
\item the branches $u_i=u|_{\Sigma_i}$ are constant for $i \geq 1$.
\end{enumerate}
\end{definition}

\begin{remark}
If $\Sigma$ is modeled on a partition $\lambda=(g_1,\ldots,g_r,(g_{r+1},h_{r+1}),\ldots,(g_{r+q},h_{r+q}))$ of $(g,h)$, then a smoothing of $\Sigma$ has genus $g$ with $h$ boundary components. The conditions we impose on $g_i$ for $i \leq r$ and on $2g_i+h_i-1$ when $i>r$ exclude unstable ghost branches.
\end{remark}

\begin{definition}
Let $\Lambda$ be the set of all partitions of $(g,h)$. For each $\lambda \in \Lambda$, the \emph{$\lambda$-cell} $\Nbar_\lambda$ is the moduli space of holomorphic curves modeled on $\lambda$. 

The moduli space of curves of type $(g,h)$ is
\[
\Nbar = \bigcup\limits_\Lambda \Nbar_\lambda.
\]
\end{definition}

Each cell of $\Nbar$ is moduli space of curves in the typical sense: it has a fixed dimension, with one top stratum and various lower-dimensional strata corresponding to degenerations of the domain. However, these cells may have $1$-dimensional boundary, and different cells may have different dimensions.

The advantage of decomposing $\Nbar$ in this manner is that individual cells are straightforward. If $\lambda=(g_1,\ldots,g_r,(g_{r+1},h_{r+1}),\ldots,(g_{r+q},h_{r+q}))$, then
\[
\Nbar_\lambda = \prod\limits_{i=1}^{r} (\Sigma_0 \times \Mbar_{g_i,1}) \times \prod\limits_{i=r+1}^{r+q} (\partial\Sigma_0\times\Mbar_{(g_i,h_i),0,(1,0,\ldots,0)}).
\]
Therefore
\[
\dim_\R(\Nbar_\lambda) = \sumto{i=1}{r} (2+2(3g_i-2)) + \sumto{i=r+1}{r+q} (1+3(2g_i+h_i-1)-2) = 3\tilde{g}-(2r+q),
\]
where $\tilde{g}=2g+h-1$ is the genus of the complex double of a curve modeled on $\lambda$.

Cells corresponding to different partitions intersect when the following phenomena occur:
\begin{enumerate}[(i)]
\item ghost branches collide, or
\item an interior ghost branch approaches $\partial\Sigma_0$.
\end{enumerate}
Note that the second type of collision corresponds to a collision of conjugate ghosts in the complex double of a curve. Thus it is possible to understand all of these collisions using standard techniques for closed curves.

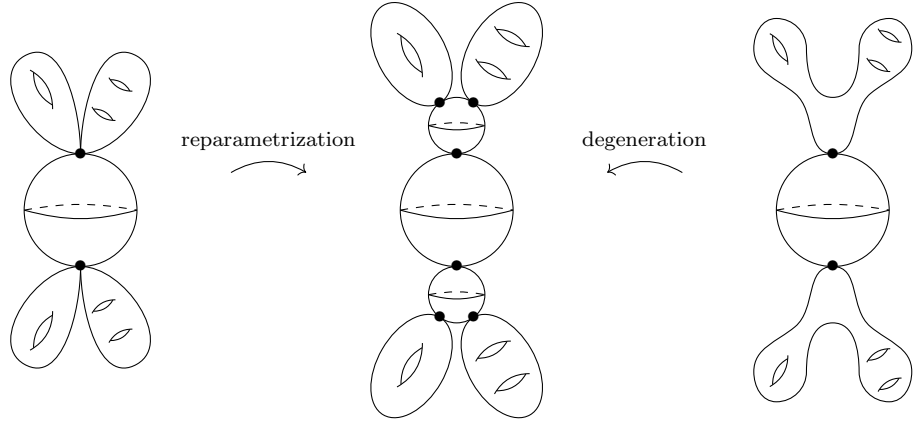
\begin{figure}[ht]
\centering
\begin{tikzpicture}[scale=0.5]

\def\r{1.5} 
\def\w{5} 


\coordinate (coll) at (0,0);

\sphere[](coll)(\r)(\r);

\coordinate (nodeC) at ($(coll)+(0,\r)$);
\node at (nodeC) {$\bullet$};

\begin{scope}[shift=(nodeC), rotate around={30:(nodeC)}]
	\draw (0,0) to [out=60,in=0] (0,3) to [out=180,in=150] (0,0);
    \addGenus[](0,2)(90)(0.5);
\end{scope}
\begin{scope}[shift=(nodeC), rotate around={-30:(nodeC)}]
    \draw (0,0) to [out=120,in=180] (0,3) to [out=0,in=30] (0,0);
    \addGenus[](0,1.2)(0)(0.3);
    \addGenus[](0,2.1)(0)(0.3);
\end{scope}

\coordinate (negNodeC) at ($(coll)+(0,-\r)$);
\node at (negNodeC) {$\bullet$};

\begin{scope}[shift=(negNodeC), yscale=-1, rotate around={30:(negNodeC)}]
	\draw (0,0) to [out=60,in=0] (0,3) to [out=180,in=150] (0,0);
    \addGenus[](0,2)(90)(0.5);
\end{scope}
\begin{scope}[shift=(negNodeC), yscale=-1, rotate around={-30:(negNodeC)}]
    \draw (0,0) to [out=120,in=180] (0,3) to [out=0,in=30] (0,0);
    \addGenus[](0,1.2)(0)(0.3);
    \addGenus[](0,2.1)(0)(0.3);
\end{scope}


\coordinate (bubble) at (2*\w,0);

\sphere[](bubble)(\r)(\r);

\coordinate (nodeB) at ($(bubble)+(0,\r)$);
\coordinate (nodeB1) at ($(nodeB)+(-0.3*\r,0.9*\r)$);
\coordinate (nodeB2) at ($(nodeB)+(0.3*\r,0.9*\r)$);

\sphere[]($(nodeB)+(0,0.5*\r)$)(0.5*\r)(0.5*\r);
\node at (nodeB) {$\bullet$};

\torus[](nodeB1)(30)(1.5);

\closedGhost[](nodeB2)(-30)(0.9)(1.5)(2);

\coordinate (negNodeB) at ($(bubble)+(0,-\r)$);
\coordinate (negNodeB1) at ($(negNodeB)+(-0.3*\r,-0.9*\r)$);
\coordinate (negNodeB2) at ($(negNodeB)+(0.3*\r,-0.9*\r)$);

\sphere[]($(negNodeB)+(0,-0.5*\r)$)(0.5*\r)(0.5*\r);
\node at (negNodeB) {$\bullet$};

\begin{scope}[yscale=-1]
    \torus[](negNodeB1)(30)(1.5);

    \closedGhost[](negNodeB2)(-30)(0.9)(1.5)(2);
\end{scope}


\coordinate (smooth) at (4*\w,0);

\sphere[](smooth)(\r)(\r);

\coordinate (nodeS) at ($(smooth)+(0,\r)$);
\node at (nodeS) {$\bullet$};

\begin{scope}[shift={(nodeS)}]
\draw (0,0) to [out=0,in=-140] (1,1.6) to [out=40,in=-90] (2.1,2.8) to [out=90,in=0] (1.5,3.6) to [out=180,in=90] (0.7,2.8) to [out=-90,in=0] (0,1.5) to [out=180,in=-90] (-0.7,2.8) to [out=90,in=0] (-1.5,3.6) to [out=180,in=90] (-2.1,2.8) to [out=-90,in=140] (-1,1.6) to [out=-40,in=180] (0,0);
\addGenus[](-1.4,2.8)(120)(0.4);
\addGenus[](1.2,2.34)(-30)(0.3);
\addGenus[](1.5,3.12)(-30)(0.3);
\end{scope}

\coordinate (negNodeS) at ($(smooth)+(0,-\r)$);
\node at (negNodeS) {$\bullet$};

\begin{scope}[shift={(negNodeS)}, yscale=-1]
\draw (0,0) to [out=0,in=-140] (1,1.6) to [out=40,in=-90] (2.1,2.8) to [out=90,in=0] (1.5,3.6) to [out=180,in=90] (0.7,2.8) to [out=-90,in=0] (0,1.5) to [out=180,in=-90] (-0.7,2.8) to [out=90,in=0] (-1.5,3.6) to [out=180,in=90] (-2.1,2.8) to [out=-90,in=140] (-1,1.6) to [out=-40,in=180] (0,0);
\addGenus[](-1.4,2.8)(120)(0.4);
\addGenus[](1.2,2.34)(-30)(0.3);
\addGenus[](1.5,3.12)(-30)(0.3);
\end{scope}


\def\h{1}
\coordinate (arr1) at ($0.5*(coll)+0.5*(bubble)+(0,\h)$);
\coordinate (arr2) at ($0.5*(bubble)+0.5*(smooth)+(0,\h)$);
\draw [->, bend left] ($(arr1)+(-1,0)$) to ($(arr1)+(1,0)$);
\draw [->, bend right] ($(arr2)+(1,0)$) to ($(arr2)+(-1,0)$);
\node [above=5pt] at (arr1) {\footnotesize reparametrization};
\node [above=5pt] at (arr2) {\footnotesize degeneration};

\end{tikzpicture}
\caption{(I) The collision of genus $1$ and $2$ ghosts as a degeneration of a genus $3$ ghost.}
\label{collClosedPic}
\end{figure}

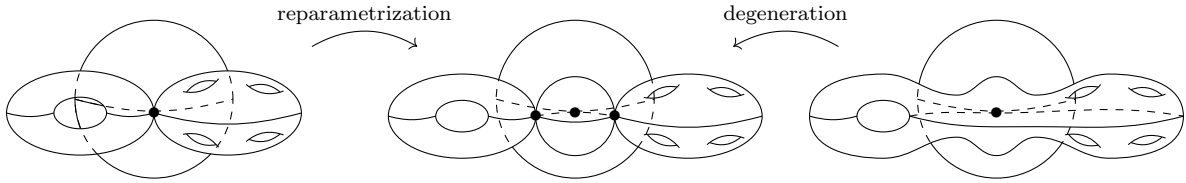
\begin{figure}[ht]
\centering
\begin{tikzpicture}[scale=0.7]

\def\r{1.5}
\def\w{4} 


\coordinate (coll) at (0,0);

\coordinate (nodeC) at ($(coll)+(0,-0.175*\r)$);

\draw (coll) circle (\r);

\begin{scope} 
	\clip ($(nodeC)+(-1.4,0)$) ellipse (1.4 and 0.8);
	\draw [white, fill=white] (coll) circle (\r+1);
	\draw [dashed] (coll) circle (\r);
\end{scope}

\begin{scope} 
	\clip ($(nodeC)+(-1.4,0)$) ellipse (0.5 and 0.3);
	\draw [bend right = 15] ($(coll)+(-\r,0)$) to ($(coll)+(\r,0)$);
	\draw (coll) circle (\r);
\end{scope}

\begin{scope} 
	\clip ($(nodeC)+(1.4,0)$) ellipse (1.4 and 0.8);
	\draw [white, fill=white] (coll) circle (\r+1);
	\draw [dashed] (coll) circle (\r);
\end{scope}

\draw [dashed, bend right = 15] ($(coll)+(-\r,0)$) to ($(coll)+(\r,0)$);
\node at (nodeC) {$\bullet$};

\begin{scope}[shift={(nodeC)}]
	\draw (-1.4,0) ellipse (1.4 and 0.8);
	\draw (-1.4,0) ellipse (0.5 and 0.3);
	\draw [bend right = 15] (-2.8,0) to (-1.9,0);
	\draw [bend right = 15] (-0.9,0) to (0,0);
    
	\draw (1.4,0) ellipse (1.4 and 0.8);
	\draw [bend right = 15] (0,0) to (2.8,0);
	\addGenus[](0.9,0.5)(20)(0.3);
	\addGenus[](2.1,0.45)(-10)(0.3);
	\addGenus[yscale=-1](0.9,-0.5)(20)(0.3);
	\addGenus[yscale=-1](2.1,-0.45)(-10)(0.3);
\end{scope}


\coordinate (bubble) at (2*\w,0);
\coordinate (nodeB) at ($(bubble)+(0,-0.175*\r)$);
\coordinate (sphere) at ($(nodeB)+(0,-0.04*\r)$);
\coordinate (nodeB1) at ($(sphere)+(-0.5*\r,0)$);
\coordinate (nodeB2) at ($(sphere)+(0.5*\r,0)$);

\draw (bubble) circle (\r);

\begin{scope} 
	\clip ($(nodeB1)+(-1.4,0)$) ellipse (1.4 and 0.8);
	\draw [white, fill=white] (bubble) circle (\r+1);
	\draw [dashed] (bubble) circle (\r);
\end{scope}

\begin{scope} 
	\clip ($(nodeB2)+(1.4,0)$) ellipse (1.4 and 0.8);
	\draw [white, fill=white] (bubble) circle (\r+1);
	\draw [dashed] (bubble) circle (\r);
\end{scope}

\draw [dashed, bend right = 15] ($(bubble)+(-\r,0)$) to ($(bubble)+(\r,0)$);

\node at (nodeB) {$\bullet$};
\node at (nodeB1) {$\bullet$};
\node at (nodeB2) {$\bullet$};

\begin{scope}[shift=(sphere)] 
	\draw (0,0) circle (0.5*\r);
	\draw [bend right = 15] (-0.5*\r,0) to (0.5*\r,0);
	\draw [dashed, bend left = 10] (-0.5*\r,0) to (0.5*\r,0);
\end{scope}

\begin{scope}[shift={(nodeB1)}] 
	\draw (-1.4,0) ellipse (1.4 and 0.8);
	\draw (-1.4,0) ellipse (0.5 and 0.3);
	\draw [bend right = 15] (-2.8,0) to (-1.9,0);
	\draw [bend right = 15] (-0.9,0) to (0,0);
\end{scope}

\begin{scope}[shift={(nodeB2)}] 
	\draw (1.4,0) ellipse (1.4 and 0.8);
	\draw [bend right = 15] (0,0) to (2.8,0);
	\addGenus[](0.9,.45)(20)(0.3);
	\addGenus[yscale=-1](0.9,-0.45)(20)(0.3);
	\addGenus[](2.1,0.45)(-10)(0.3);
	\addGenus[yscale=-1](2.1,-0.45)(-10)(0.3);
\end{scope}


\coordinate (smooth) at (4*\w,0);

\coordinate (nodeS) at ($(smooth)+(0,-0.175*\r)$);
\coordinate (sphereS) at ($(nodeS)+(0,-0.04*\r)$);

\draw (smooth) circle (\r);

\begin{scope}[shift={(sphereS)}]
	\path[preaction={draw,fill=white}][clip] (0,0.5*\r) to [out=180, in=0] (-0.5*\r,0.4) to [out=180, in=0] (-0.5*\r-1.4,0.8) to [out=180, in=90] (-0.5*\r-2.8,0) to [out=-90, in=180] (-0.5*\r-1.4,-0.8) to [out=0, in=180] (-0.5*\r,-0.4) to [out=0, in=180] (0,-0.5*\r) to [out=0, in=180] (0.5*\r,-0.4) to [out=0, in=180] (0.5*\r+1.4,-0.8) to [out=0, in=-90] (0.5*\r+2.8,0) to [out=90, in=0] (0.5*\r+1.4,0.8) to [out=180, in=0] (0.5*\r,0.4) to [out=180, in=0] (0,0.5*\r);

	\draw (-0.5*\r-0.9,0) to [out=-15,in=180] (0.5*\r,-0.2) to [out=0,in=195] (0.5*\r+2.8,0);
	\draw [dashed] (-0.5*\r-0.9,0) to [out=10,in=180] (nodeS) to [out=0,in=170] (0.5*\r+2.8,0);

    \draw[dashed] (smooth) circle (\r);
\end{scope}

\draw [dashed, bend right = 15] ($(smooth)+(-\r,0)$) to ($(smooth)+(\r,0)$);

\node at (nodeS) {$\bullet$};

\begin{scope}[shift={($(sphereS)+(-0.5*\r,0)$)}]
	\draw (-1.4,0) ellipse (0.5 and 0.3);
	\draw [bend right = 15] (-2.8,0) to (-1.9,0);
\end{scope}

\begin{scope}[shift={($(sphereS)+(0.5*\r,0)$)}]
	\addGenus[](0.9,0.45)(20)(0.3);
	\addGenus[yscale=-1](0.9,-0.45)(20)(0.3);
	\addGenus[](2.1,0.45)(-10)(0.3);
	\addGenus[yscale=-1](2.1,-0.45)(-10)(0.3);
\end{scope}


\def\h{1}
\coordinate (arr1) at ($0.5*(coll)+0.5*(bubble)+(0,\h)$);
\coordinate (arr2) at ($0.5*(bubble)+0.5*(smooth)+(0,\h)$);
\draw [->, bend left] ($(arr1)+(-1,0)$) to ($(arr1)+(1,0)$);
\draw [->, bend right] ($(arr2)+(1,0)$) to ($(arr2)+(-1,0)$);
\node [above=5pt] at (arr1) {\footnotesize reparametrization};
\node [above=5pt] at (arr2) {\footnotesize degeneration};

\end{tikzpicture}
\caption{(II) The collision of ghosts of type $(0,2)$ and $(2,1)$ as a degeneration of a ghost of type $(2,2)$.}
\label{collOpenPic}
\end{figure}

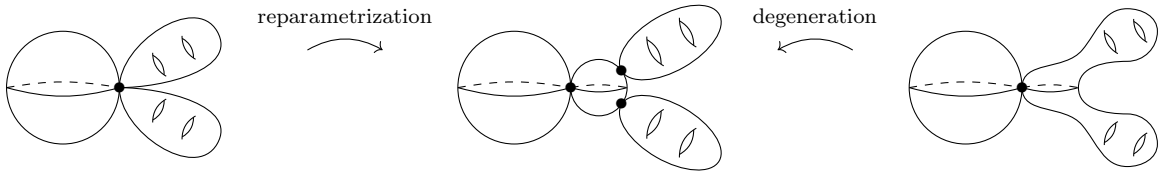
\begin{figure}[ht]
\centering
\begin{tikzpicture}[scale=0.5]

\def\r{1.5}
\def\w{6}


\coordinate (coll) at (0,0);

\sphere[](coll)(\r)(\r);

\coordinate (nodeC) at ($(coll)+(\r,0)$);
\node at (nodeC) {$\bullet$};

\begin{scope}[shift=(nodeC), rotate around={-60:(nodeC)}]
    \draw (0,0) to [out=60,in=0] (0,3) to [out=180,in=150] (0,0);
    \addGenus[](0,1.2)(0)(0.3);
    \addGenus[](0,2.1)(0)(0.3);
\end{scope}

\begin{scope}[yscale=-1, shift=(nodeC), rotate around={-60:(nodeC)}]
    \draw (0,0) to [out=60,in=0] (0,3) to [out=180,in=150] (0,0);
    \addGenus[](0,1.2)(0)(0.3);
    \addGenus[](0,2.1)(0)(0.3);
\end{scope}


\coordinate (bubble) at (2*\w,0);

\sphere[](bubble)(\r)(\r);

\coordinate (nodeB) at ($(bubble)+(\r,0)$);
\node at (nodeB) {$\bullet$};
\sphere[]($(nodeB)+(0.5*\r,0)$)(0.5*\r)(0.5*\r);

\coordinate (nodeB1) at ($(nodeB)+(0.9*\r,0.3*\r)$);
\closedGhost[](nodeB1)(-60)(0.75)(1.5)(2);

\coordinate (negNodeB1) at ($(nodeB)+(0.9*\r,-0.3*\r)$);
\begin{scope}[shift=(negNodeB1), yscale=-1]
    \closedGhost[](0,0)(-60)(0.75)(1.5)(2);
\end{scope}


\coordinate (smooth) at (4*\w,0);

\sphere[](smooth)(\r)(\r);

\coordinate (nodeS) at ($(smooth)+(\r,0)$);
\node at (nodeS) {$\bullet$};

\begin{scope}[shift={(nodeS)}, rotate around={-90:(nodeS)}]
	\draw [bend right = 15] (0,0) to (0,1.5);
	\draw [dashed, bend left = 10] (0,0) to (0,1.5);
	\draw (0,0) to [out=0,in=-140] (1,1.6) to [out=40,in=-90] (2.1,2.8) to [out=90,in=0] (1.5,3.6) to [out=180,in=90] (0.7,2.8) to [out=-90,in=0] (0,1.5) to [out=180,in=-90] (-0.7,2.8) to [out=90,in=0] (-1.5,3.6) to [out=180,in=90] (-2.1,2.8) to [out=-90,in=140] (-1,1.6) to [out=-40,in=180] (0,0);
	\addGenus[](-1.2,2.34)(30)(0.3);
	\addGenus[](-1.5,3.12)(30)(0.3);
	\addGenus[](1.2,2.34)(-30)(0.3);
	\addGenus[](1.5,3.12)(-30)(0.3);
\end{scope}


\def\h{1}
\coordinate (arr1) at ($0.5*(nodeC)+0.5*(nodeB)+(0,\h)$);
\coordinate (arr2) at ($0.5*(nodeB)+0.5*(nodeS)+(0.5,\h)$);
\draw [->, bend left] ($(arr1)+(-1,0)$) to ($(arr1)+(1,0)$);
\draw [->, bend right] ($(arr2)+(1,0)$) to ($(arr2)+(-1,0)$);
\node [above=5pt] at (arr1) {\footnotesize reparametrization};
\node [above=5pt] at (arr2) {\footnotesize degeneration};

\end{tikzpicture}
\caption{(III) The collision of a genus $2$ ghost with the boundary as a degeneration of a ghost of type $(2,1)$.}
\label{collBdryPic}
\end{figure}

Our next goal is to characterize these cell intersections more explicitly. Figures~\ref{collClosedPic}, \ref{collOpenPic}, and \ref{collBdryPic} depict the complex doubles of the three basic intersection types.

If a curve in $\Nbar_\lambda$ has two ghosts attached at $z_i$ and $z_{i'}$, then $\Nbar_\lambda$ intersects another cell when $z_i=z_{i'}$. The collision produces an extra bubble (a sphere in the closed case and a disk in the open case); the resulting ghost branch is a degeneration of a single ghost component attached at $z_i=z_{i'}$ .

If a curve in $\Nbar_\lambda$ has a closed ghost attached at $z_i$, then $\Nbar_\lambda$ intersects another cell when $z_i$ approaches $\partial\Sigma_0$. The collision produces an extra disk bubble; the resulting ghost branch is a degeneration of an open ghost attached at $z_i$.

\begin{enumerate}[(I)]
\item If two closed ghosts $[\Sigma_i,y_i] \in \Mbar_{g_i,1}$ and $[\Sigma_{i'},y_{i'}] \in \Mbar_{g_{i'},1}$ collide, we see a sphere bubble attached to $\Sigma_0$, $\Sigma_i$, and $\Sigma_{i'}$. This is a degeneration of a closed ghost with genus $g_i+g_{i'}$. \label{collisionClosed}
\item If two open ghosts $[\Sigma_i,y_i] \in \Mbar_{(g_i,h_i),0,(1,0,\ldots,0)}$ and $[\Sigma_{i'},y_{i'}] \in \Mbar_{(g_{i'},h_{i'}),0,(1,0,\ldots,0)}$ collide, we see a disk bubble attached to $\partial\Sigma_0$, $\partial\Sigma_i$, and $\partial\Sigma_{i'}$. This is a degeneration of an open ghost of type $(g_i+g_{i'},h_i+h_{i'}-1)$. \label{collisionOpen}
\item If a closed ghost $[\Sigma_i,y_i] \in \Mbar_{g_i,1}$ approaches $\partial\Sigma_0$, we see a disk bubble attached to $\partial\Sigma_0$ and $\Sigma_i$. This is a degeneration of an open ghost of type $(g_i,1)$. \label{collisionBoundary}
\end{enumerate}

Cell intersections may be more complicated if there are many ghost branches: for example, more than two ghosts may collide. However, all intersections may be understood from these rules. For instance, if a closed ghost of genus $g_i$ collides with an open ghost of type $(g_{i'},h_{i'})$, we combine rules~(\ref{collisionOpen}) and (\ref{collisionBoundary}) to see that an open ghost of type $(g_i+g_{i'},h_{i'})$ appears.

\begin{remark} \label{basicInt}
These basic intersections are the only possible intersections of exactly two cells. All cell intersections arise from these three types in the following sense. Given a cell intersection $\Nbar_{\lambda_1} \cap \ldots \cap \Nbar_{\lambda_k} \neq \emptyset$, we can build a graph graph with vertices $\lambda_1,\ldots,\lambda_k$ and edges between $\lambda_i$ and $\lambda_{i'}$ precisely when the intersection $\Nbar_{\lambda_i} \cap \Nbar_{\lambda_{i'}}$ corresponds to one of the three basic types. Then this graph is connected, so any pairwise intersection $\Nbar_{\lambda_i} \cap \Nbar_{\lambda_{i'}}$ can be viewed as a sequence of basic intersection types.
\end{remark}

The following lemma makes precise the way collisions occur in moduli spaces of domains.

\begin{lemma} \label{collisionModuli}
There exist canonical inclusions
\begin{align*}
\rho_{g_i,g_{i'}}:\Mbar_{g_i,1} \times \Mbar_{g_{i'},1} & \hookrightarrow \Mbar_{g_i+g_{i'},1}
\\
\rho_{(g_i,h_i),(g_{i'},h_{i'})}:\Mbar_{(g_i,h_i),0,(1,0,\ldots,0)} \times \Mbar_{(g_{i'},h_{i'}),0,(1,0,\ldots,0)} & \hookrightarrow \Mbar_{(g_i+g_{i'},h_i+h_{i'}-1),0,(1,0,\ldots,0)}
\\
\rho_{g_i}:\Mbar_{g_i,1} & \hookrightarrow \Mbar_{(g_{i},1),0,(1,0,\ldots,0)}.
\end{align*}
Under these inclusions, there are canonical isomorphisms
\begin{align*}
\E_{g_i}^* \oplus \E_{g_{i'}}^* & \cong (\E_{g_i+g_{i'}}^*)|_{Im(\rho_{g_i,g_{i'}})}
\\
\E_{(g_i,h_i)}^* \oplus \E_{(g_{i'},h_{i'})}^* & \cong (\E_{(g_i+g_{i'},h_i+h_{i'}-1)}^*)|_{Im(\rho_{(g_i,h_i),(g_{i'},h_{i'})})}
\\
\E_{g_i}^* & \cong (\E_{(g_{i},1)}^*)|_{Im(\rho_{g_i})}.
\end{align*}
\begin{proof}
For the map $\rho_{g_i,g_{i'}}$, we identify a pair $((\Sigma,z),(\Sigma',z'))$ of closed surfaces with the nodal surface obtained by attaching $\Sigma$ and $\Sigma'$ to a sphere $(S^2,y)$ at $z$ and $z'$. This sphere then has two nodes and one marked point, and since $\Mbar_{0,3}=\{\text{pt}\}$ these nodal curves are in bijection with pairs $((\Sigma,z),(\Sigma',z'))$.
\begin{figure}[ht]
\centering
\begin{tikzpicture}

\def\r{1}
\def\h{1.2}
\def\w{0.8}

\sphere[](0,0)(\r)(\r);
\node at (0,-\r) {$\bullet$};
\node [below] at (0,-\r) {$y$};

\coordinate (node1) at (-0.8*\r,0.6*\r);
\closedGhost[](node1)(50)(\w)(\h)(2);
\begin{scope}[shift={(node1)}, rotate around={50:(node1)}]
	\node [above left] at (0,2*\h) {$\Sigma$};
\end{scope}

\coordinate (node2) at (0.8*\r,0.6*\r);
\closedGhost[](node2)(-50)(\w)(\h)(3);
\begin{scope}[shift={(node2)}, rotate around={-50:(node2)}]
	\node [above right] at (0,2*\h) {$\Sigma'$};
\end{scope}

\end{tikzpicture}
\caption{Embedding $\Mbar_{g_i,1} \times \Mbar_{g_{i'},1}$ in $\Mbar_{g_i+g_{i'},1}$.}
\end{figure}
It is evident that there is an injective map $\coker(\delbar_{\Sigma}) \arr \coker(\delbar_{\rho(\Sigma,\Sigma')}) \oplus \coker(\delbar_{\Sigma'})$. To show that this map is an isomorphism, we only need to observe that these two spaces have the same dimension.

Next we define $\rho_{(g_i,h_i),(g_{i'},h_{i'})}$. Let $\tilde{g}_i=2g_i+h_i-1$ and $\tilde{g}_{i'}=2g_{i'}+h_{i'}-1$. For a pair of open surfaces $((\Sigma,z),(\Sigma',z'))$, let $((\Sigma^{(\C)},z),((\Sigma')^{(\C)},z'))$ be the pair of complex doubles. This pair lives in the real locus of $\Mbar_{\tilde{g}_i,1} \oplus \Mbar_{\tilde{g}_{i'},1}$ (see Section 3.3.3 of\cite{katzLiu}). We can then attach $\Sigma^{(\C)}$ and $(\Sigma')^{(\C)}$ to a sphere using $\rho_{\tilde{g}_i,\tilde{g}_{i'}}$, taking care to pick the special points on the sphere along the real locus. This new nodal curve is the symmetric double of a curve in $\Mbar_{(g_i+g_{i'},h_i+h_{i'}-1),0,(1,0,\ldots,0)}$, and such nodal curves are in bijection with pairs $((\Sigma,z),(\Sigma',z'))$ because $\Mbar_{0,3}=\{\text{pt}\}$.
\begin{figure}[ht]
\centering
\begin{tikzpicture}

\def\r{1}

\node at (-\r,0) {$\bullet$};
\node at (\r,0) {$\bullet$};
\node at (0,-0.175*\r) {$\bullet$};
\node [below] at (0,-0.175*\r) {$y$};

\sphere[](0,0)(\r)(\r);

\begin{scope}[shift={(-2*\r,0)}]
	\sphere[](0,0)(\r)(1.4*\r);
	
	\addGenus[](-.45*\r,0.75*\r)(70)(0.3);
	\addGenus[](0.45*\r,0.6*\r)(-30)(0.3);
	\addGenus[](-.45*\r,-0.75*\r)(110)(0.3);
	\addGenus[](0.45*\r,-0.6*\r)(-150)(0.3);
	
	\node [left] at (-\r,0) {$\Sigma^{(\C)}$};
\end{scope}

\begin{scope}[shift={(2*\r,0)}]
	\draw (0,0) ellipse ({\r} and {1.4*\r});
	\draw [bend right=15] (-\r,0) to (-0.5*\r,0);
	\draw [dashed, bend left=10] (-\r,0) to (-0.5*\r,0);
	\draw [bend right=15] (0.5*\r,0) to (\r,0);
	\draw [dashed, bend left=10] (0.5*\r,0) to (\r,0);
	\draw (0.5*\r,0) arc (0:360:{0.5*\r} and {0.25*\r});
	
	\addGenus[](0,0.8*\r)(0)(0.4);
	\addGenus[](0,-0.8*\r)(180)(0.4);
	
	\node [right] at (\r,0) {$(\Sigma')^{(\C)}$};
\end{scope}

\end{tikzpicture}
\caption{Embedding $\Mbar_{(g_i,h_i),0,(1,0,\ldots,0)} \times \Mbar_{(g_{i'},h_{i'}),0,(1,0,\ldots,0)}$ in $\Mbar_{(g_i+g_{i'},h_i+h_{i'}-1),0,(1,0,\ldots,0)}$.}
\end{figure}

It is evident that there is an injective map $\coker(\delbar_{\Sigma}) \arr \coker(\delbar_{\rho(\Sigma,\Sigma')}) \oplus \coker(\delbar_{\Sigma'})$. To show that this map is an isomorphism, we only need to observe that these two spaces have the same dimension.

Finally, we define $\rho_{g_i}$. Fix a closed surface $(\Sigma,z)$. We define a symmetric nodal curve in $\Mbar_{2g_i,1}$ by attaching $(\Sigma,z)$ and $(\ov{\Sigma},\ov{z})$ to a sphere $(S^2,y)$ at $z$ and $\ov{z}$. This sphere then has two nodes and one marked point, and since $\Mbar_{0,3}=\{\text{pt}\}$ these nodal curves are in bijection with closed surfaces $(\Sigma,z)$. If we pick $y$ in the real locus of $S^2$, this nodal curve is the complex double of a curve in $\Mbar_{(g_{i},1),0,(1,0,\ldots,0)}$.
\begin{figure}[ht]
\centering
\begin{tikzpicture}

\def\r{1}
\def\h{1.2}
\def\w{0.8}

\sphere[](0,0)(\r)(\r);
\node at (-\r,0) {$\bullet$};
\node [left] at (-\r,0) {$y$};

\coordinate (node1) at (0.8*\r,0.6*\r);
\closedGhost[](node1)(-40)(\w)(\h)(2);
\begin{scope}[shift={(node1)}, rotate around={-40:(node1)}]
	\node [above right] at (0,2*\h) {$\Sigma$};
\end{scope}

\coordinate (node2) at (0.8*\r,-0.6*\r);
\closedGhost[](node2)(-130)(\w)(\h)(2);
\begin{scope}[shift={(node2)}, rotate around={-130:(node2)}]
	\node [below right] at (0,2*\h) {$\ov{\Sigma}$};
\end{scope}

\end{tikzpicture}
\caption{Embedding $\Mbar_{g_i,1}$ in $\Mbar_{(g_{i},1),0,(1,0,\ldots,0)}$.}
\end{figure}

It is evident that there is an injective map $\coker(\delbar_{\Sigma}) \arr \coker(\delbar_{\rho(\Sigma)})$. To show that this map is an isomorphism, we only need to observe that these two spaces have the same dimension.
\end{proof}
\end{lemma}

The following examples demonstrate some of the ways in which cells can intersect.

\begin{example}
If $(g,h)=(1,1)$, then there are precisely two partitions of $(g,h)$: 
\begin{align*}
\lambda_1 & = (g_1=1)
\\
\lambda_2 & = ((g_1=1,h_1=1)).
\end{align*}
Partition $\lambda_1$ corresponds to those curves with a closed ghost torus, and partition $\lambda_2$ corresponds to those curves with an open ghost which has genus one and one boundary component. The intersection $\Nbar_{\lambda_1} \cap \Nbar_{\lambda_2}$ is the set of curves in $\Nbar_{\lambda_1}$ where the closed ghost torus hits $\partial\Sigma_0$, as in Figure~\ref{pic11sigma}. See Subsection~\ref{base11ss}.
\end{example}

\begin{example} \label{int4moduli}
Consider the topological type $(3,1)$. We consider the intersection of four partitions:
\[
\begin{array}{rclcrcl}
\lambda_1 & = & (g_1=1, g_2=2) & \qquad\qquad & \dim_\R(\Nbar_{\lambda_{1}}) & = & 14
\\
\lambda_2 & = & ((g_1=1,h_1=1),g_2=2) && \dim_\R(\Nbar_{\lambda_{2}}) & = & 15
\\
\lambda_3 & = & (g_1=1,(g_2=2,h_2=1)) && \dim_\R(\Nbar_{\lambda_{3}}) & = & 15
\\
\lambda_4 & = & ((g_1=1,h_1=1),(g_2=2,h_2=1)) && \dim_\R(\Nbar_{\lambda_{3}}) & = & 16.
\end{array}
\]
We label the nodes of curves in these cells as in Figure~\ref{picFour}.

\begin{figure}[ht]
\centering
\begin{tikzpicture}

\def\r{1}
\def\h{1}
\def\w{0.6}
\def\d{4}       
\def\v{-0.6}    

\coordinate (c1) at (0,\d);
\coordinate (c2) at (1.5*\d,\d);
\coordinate (c3) at (0,0);
\coordinate (c4) at (1.5*\d,0);


\disk[](c1)(\r)(\r);

\coordinate (LNode1) at ($(c1)+(-0.6*\r,0.8*\r)$);
\node [below right] at (LNode1) {$a$};
\torus[](LNode1)(30)(\h);

\coordinate (RNode1) at ($(c1)+(0.6*\r,0.8*\r)$);
\node [below left] at (RNode1) {$b$};
\closedGhost[](RNode1)(-30)(\w)(\h)(2);

\node at ($(c1)+(0,\v)$) {$\lambda_1=(1,2)$};


\disk[](c2)(\r)(\r);

\coordinate (LNode2) at ($(c2)+(-\r,0)$);
\node [below right] at (LNode2) {$c$};
\openGhost[](LNode2)(\w)(\h)(1)(left);
\begin{scope}[shift={(LNode2)}]
	\addGenus[](-1*\w,0.5*\h)(0)(0.25);
\end{scope}

\coordinate (RNode2) at ($(c2)+(0.6*\r,0.8*\r)$);
\node [below] at (RNode2) {$b'$};
\closedGhost[](RNode2)(-30)(\w)(\h)(2);

\node at ($(c2)+(0,\v)$) {$\lambda_2=((1,1),2)$};


\disk[](c3)(\r)(\r);

\coordinate (LNode3) at ($(c3)+(-0.6*\r,0.8*\r)$);
\node [below right] at (LNode3) {$a'$};
\torus[](LNode3)(30)(\h);

\coordinate (RNode3) at ($(c3)+(\r,0)$);
\node [below] at (RNode3) {$d$};
\openGhost[](RNode3)(\w)(1.5*\h)(1)();
\begin{scope}[shift={(RNode3)}]
	\addGenus[](\w,0.4*\h)(0)(0.2);
	\addGenus[](\w,0.8*\h)(0)(0.2);
\end{scope}

\node at ($(c3)+(0,\v)$) {$\lambda_3=(1,(2,1))$};


\disk[](c4)(\r)(\r);

\coordinate (LNode4) at ($(c4)+(-\r,0)$);
\node [below right] at (LNode4) {$c'$};
\openGhost[](LNode4)(\w)(\h)(1)(left);
\begin{scope}[shift={(LNode4)}]
	\addGenus[](-1*\w,0.5*\h)(0)(0.25);
\end{scope}

\coordinate (RNode4) at ($(c4)+(\r,0)$);
\node [below] at (RNode4) {$d'$};
\openGhost[](RNode4)(\w)(1.5*\h)(1)();
\begin{scope}[shift={(RNode4)}]
	\addGenus[](\w,0.4*\h)(0)(0.2);
	\addGenus[](\w,0.8*\h)(0)(0.2);
\end{scope}

\node at ($(c4)+(0,\v)$) {$\lambda_4=((1,1),(2,1))$};

\end{tikzpicture}
\caption{Curves in $\Nbar_{\lambda_1}$, $\Nbar_{\lambda_2}$, $\Nbar_{\lambda_3}$, and $\Nbar_{\lambda_4}$.}
\label{picFour}
\end{figure}

In this example we examine the intersection $\Nbar_{\lambda_1} \cap \Nbar_{\lambda_2} \cap \Nbar_{\lambda_3} \cap \Nbar_{\lambda_4}$ in the moduli space. See Example~\ref{int4glue} for gluing parameters over this intersection.

The intersection $\Nbar_{\lambda_1} \cap \Nbar_{\lambda_2}$ consists of curves in $\Nbar_{\lambda_1}$ where the genus $1$ ghost approaches the boundary. The intersection $\Nbar_{\lambda_1} \cap \Nbar_{\lambda_3}$ consists of curves in $\Nbar_{\lambda_1}$ where the genus $2$ ghost approaches the boundary. The intersections $\Nbar_{\lambda_2} \cap \Nbar_{\lambda_4}$ and $\Nbar_{\lambda_3} \cap \Nbar_{\lambda_4}$ consist of curves in $\Nbar_{\lambda_2}$ and $\Nbar_{\lambda_3}$, respectively, where the interior ghost approaches the boundary. The intersections $\Nbar_{\lambda_1} \cap \Nbar_{\lambda_4}$, $\Nbar_{\lambda_2} \cap \Nbar_{\lambda_3}$, and $\Nbar_{\lambda_1} \cap \Nbar_{\lambda_2} \cap \Nbar_{\lambda_3} \cap \Nbar_{\lambda_4}$ all consist of curves illustrated in Figure~\ref{pic4int}. We compute the dimensions:
\[
\begin{array}{rclcrcl}
\dim_\R(\Nbar_{\lambda_1} \cap \Nbar_{\lambda_2}) & = & 13 & \qquad\qquad & \dim_\R(\Nbar_{\lambda_1} \cap \Nbar_{\lambda_3}) & = & 13
\\
\dim_\R(\Nbar_{\lambda_2} \cap \Nbar_{\lambda_4}) & = & 14 && \dim_\R(\Nbar_{\lambda_3} \cap \Nbar_{\lambda_4}) & = & 14
\\
\dim_\R(\Nbar_{\lambda_1} \cap \Nbar_{\lambda_2} \cap \Nbar_{\lambda_3} \cap \Nbar_{\lambda_4}) & = & 12.
\end{array}
\]

\begin{figure}[ht]
\centering
\begin{tikzpicture}

\def\r{1}
\def\h{1}
\def\w{0.6}
\def\d{5}

\coordinate (main) at (0,0);
\coordinate (lDisk) at (-2*\r,0);
\coordinate (rDisk) at (2*\r,0);

\disk[](main)(\r)(\r);
\disk[](lDisk)(\r)(\r);
\disk[](rDisk)(\r)(\r);

\node at (-\r,0) {$\bullet$};
\node [below] at (-\r,0) {$c=c'$};
\node at (\r,0) {$\bullet$};
\node [below] at (\r,0) {$d=d'$};

\coordinate (LNode1) at ($(lDisk)+(-0.6*\r,0.8*\r)$);
\torus[](LNode1)(30)(\h);
\begin{scope}[shift=(LNode1), rotate around={-60:(LNode1)}]
    \node [left=1pt] at (0,-0.4) {$a=a'$};
\end{scope}

\coordinate (RNode1) at ($(rDisk)+(0.6*\r,0.8*\r)$);
\closedGhost[](RNode1)(-30)(\w)(\h)(2);
\begin{scope}[shift=(RNode1), rotate around={-30:(RNode1)}]
    \node [right=1pt] at (0.4,0) {$b=b'$};
\end{scope}

\end{tikzpicture}
\caption{A curve in $\Nbar_{\lambda_1} \cap \Nbar_{\lambda_2} \cap \Nbar_{\lambda_3} \cap \Nbar_{\lambda_4}$.}
\label{pic4int}
\end{figure}
Given a curve in $\Nbar_{\lambda_1} \cap \Nbar_{\lambda_2} \cap \Nbar_{\lambda_3} \cap \Nbar_{\lambda_4}$, smoothing one node or a pair of nodes will produce a curve in one of the larger moduli spaces, as described in Table~\ref{table4int}
.\begin{table}[ht]
\centering
\begin{tabular}{|c|c|}
\hline
Smoothed nodes & Moduli space
\\
\hline
$a$ & $\Nbar_{\lambda_2} \cap \Nbar_{\lambda_4}$
\\
\hline
$b$ & $\Nbar_{\lambda_3} \cap \Nbar_{\lambda_4}$
\\
\hline
$c$ & $\Nbar_{\lambda_1} \cap \Nbar_{\lambda_3}$
\\
\hline
$d$ & $\Nbar_{\lambda_1} \cap \Nbar_{\lambda_2}$
\\
\hline
$c$ and $d$ & $\Nbar_{\lambda_1}$
\\
\hline
$a$ and $d$ & $\Nbar_{\lambda_2}$
\\
\hline
$b$ and $c$ & $\Nbar_{\lambda_3}$
\\
\hline
$a$ and $b$ & $\Nbar_{\lambda_4}$
\\
\hline
\end{tabular}
\caption{Smoothing nodes in $\Nbar_{\lambda_1} \cap \Nbar_{\lambda_2} \cap \Nbar_{\lambda_3} \cap \Nbar_{\lambda_4}$.}
\label{table4int}
\end{table}
\end{example}

\begin{remark}
If $h>1$, then every cell is of the form $\Nbar_\lambda=S^1 \times \mathcal{X}$, which will force the contribution of $(g,h)$-type covers to be zero (see Section~\ref{calcS}). Even if the main component had some other topological type, we would obtain a similar result whenever $h-h_0 \equiv 1 \pmod{2}$. This condition is equivalent to the doubled curve having odd excess genus; it is well-known in the closed case that such curves do not contribute.

\begin{figure}[ht]
\centering
\begin{tikzpicture}

\disk[](0,0)(2)(2);
\openGhost[](2,0)(1.5)(2)(2)();
\addGenus[](3,1)(40)(0.25);
\addGenus[](4.25,0.75)(-20)(0.25);

\begin{scope}[xshift=8cm]
    \sphere[](0,0)(2)(2);
    \openGhost[](2,0)(1.5)(2)(2)();
    \draw (5,0) arc (0:-180:1.5 and 2);
    \draw (4,0) arc (0:-180:0.5 and 0.5);
    \addGenus[](3,1)(40)(0.25);
    \addGenus[](3,-1)(140)(0.25);
    \addGenus[](4.25,0.75)(-20)(0.25);
    \addGenus[](4.25,-0.75)(200)(0.25);
\end{scope}

\end{tikzpicture}
\caption{A curve modeled on $((2,2))$ and its genus $5$ double.}
\end{figure}
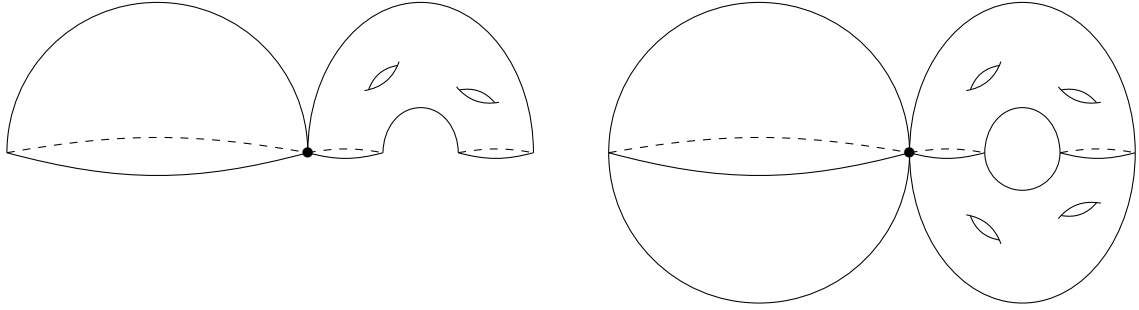
\end{remark}

\section{Obstruction Bundle} \label{obS}

In this section we compute the obstruction bundle over each moduli space of holomorphic curves (cf. Section~5 of \cite{dw}), the fiber of which is the cokernel of the linearization. 
Fix a topological type $(g,h)$ and let $\Nbar$ be the moduli space computed in Section~\ref{baseS}. In Lemmas~\ref{kerClosed} and \ref{kerOpen} we compute the kernels and cokernels of the linearizations over individual ghost components. In Proposition~\ref{fiber} we determine the fiber over the top stratum of each cell. Finally, in Proposition~\ref{ob} we examine the intersections of cells.

\begin{lemma} \label{kerClosed}
Suppose that $(C,\partial C)$ is a disk holomorphically embedded in $(M,L)$ and fix $p_i \in C\setminus\partial C$. Fix a closed domain $\Sigma_i$ and assume that $u_i:\Sigma_i \arr M$ is constant with value $p_i$. 
If we decompose with respect to the tangent space splitting $T_{p_i}M \cong T_{p_i}C \oplus N_{p_i}C$, then the linearization $\hat{D}_i:\Gamma(\Sigma_i;u_i^*TM)\arr\Omega^{0,1}(\Sigma_i,u_i^*TM)$ satisfies
\[
\begin{array}{rclcrcl}
\ker(\hat{D}_i^N) & \cong & N_{p_i}C & \qquad\qquad & \coker(\hat{D}_i^N) & \cong & N_{p_i}C \otimes_\C H^{0,1}(\Sigma_i)
\\
\ker(\hat{D}_i^T) & \cong & T_{p_i}C && \coker(\hat{D}_i^T) & \cong & T_{p_i}C \otimes_\C H^{0,1}(\Sigma_i)
\\
\ker(\hat{D}_i) & \cong & T_{p_i}M && \coker(\hat{D}_i) & \cong & T_{p_i}M \otimes_\C H^{0,1}(\Sigma_i).
\end{array}
\]
\begin{proof}
The arguments for $\hat{D}_i^N$ and $\hat{D}_i^T$ are identical. Let $V \arr \Sigma_i$ be a trivial bundle, equal to either $u_i^*TC$ or $u_i^*NC$ (with fiber $V_{y_i}$), and let $\hat{D}_i^V$ be the part of $\hat{D}$ which corresponds to $V$. Then the domain and codomain of $\hat{D}_i^V$ are
\begin{align*}
\Gamma(\Sigma_i; V) & \cong \Map(\Sigma_i;V_{y_i})
\\
\Omega^{0,1}(\Sigma_i;V) & \cong \Omega^{0,1}(\Sigma_i;\C) \otimes_\C V_{y_i}.
\end{align*}
Under this identification, the linearization is
\begin{align*}
\hat{D}_i(\xi_i,k_i) & = \df{1}{2} \left( \nabla\xi_i + J(u_i) \circ \nabla\xi_i \circ j_i + (\nabla_{\xi_i} J) \circ du_i \circ j_i \right) + \df{1}{2}J \circ du_i \circ k_i
\\
& = \delbar\xi_i.
\end{align*}
Thus $\ker(\hat{D}_i^V)$ is the set of holomorphic functions $\Sigma_i \arr V_{y_i}$, all of which are constant by Lemma~\ref{holFnRS}, and $\coker(\hat{D}_i^V)$ is just $V_{y_i} \otimes_\C \coker(\delbar_{\Sigma_i})$.
\end{proof}
\end{lemma}

\begin{lemma} \label{kerOpen}
Suppose that $(C,\partial C)$ is a disk holomorphically embedded in $(M,L)$ and fix $p_i~\in~\partial~C$. Fix an open domain $\Sigma_i$ and assume that $u_i:(\Sigma_i,\partial\Sigma_i) \arr (M,L)$ is constant with value $p_i \in L$. 
If we decompose with respect to compatible tangent space splittings $T_{p_i}M~\cong~T_{p_i}C~\oplus~N_{p_i}C$ and $T_{p_i}L \cong T_{p_i}\partial C \oplus N_{p_i}\partial C$, then the linearization $\hat{D}_i~:~\Gamma(\Sigma_i,\partial\Sigma_i;u_i^*TM,u_i^*TL)~\arr~\Omega^{0,1}(\Sigma_i,u_i^*TM)$ satisfies
\[
\begin{array}{rclcrcl}
\ker(\hat{D}_i^N) & \cong & N_{p_i}\partial C & \qquad\qquad & \coker(\hat{D}_i^N) & \cong & N_{p_i}\partial C \otimes_\R H^{0,1}(\Sigma_i)
\\
\ker(\hat{D}_i^T) & \cong & T_{p_i}\partial C && \coker(\hat{D}_i^T) & \cong & T_{p_i}\partial C \otimes_\R H^{0,1}(\Sigma_i)
\\
\ker(\hat{D}_i) & \cong & T_{p_i}L && \coker(\hat{D}_i) & \cong & T_{p_i}L \otimes_\R H^{0,1}(\Sigma_i).
\end{array}
\]
\begin{proof}
The arguments for $\hat{D}_i^N$ and $\hat{D}_i^T$ are identical. Let $(V,V^{(\R)}) \arr (\Sigma_i,\partial\Sigma_i)$ be a trivial bundle pair, equal to either $(u_i^*TC,u_i^*T\partial C)$ or $(u_i^*NC,u_i^*N\partial C)$, with fiber $V_{y_i}$ over $\Sigma_i$ and totally real fiber $V_{y_i}^{(\R)}$ over $\partial\Sigma_i$. Let $\hat{D}_i^V$ be the part of $\hat{D}$ which corresponds to $V$. Then the domain and codomain of $\hat{D}_i^V$ are
\begin{align*}
\Gamma(\Sigma_i,\partial\Sigma_i; V,V^{(\R)}) & \cong \Map(\Sigma_i,\partial\Sigma_i;V_{y_i},V_{y_i}^{(\R)})
\\
\Omega^{0,1}(\Sigma_i;V) & \cong \Omega^{0,1}(\Sigma_i;\C) \otimes_\C V_{y_i}.
\end{align*}
Under this identification, the linearization is
\begin{align*}
\hat{D}_i(\xi_i,k_i) & = \df{1}{2} \left( \nabla\xi_i + J(u_i) \circ \nabla\xi_i \circ j_i + (\nabla_{\xi_i} J) \circ du_i \circ j_i \right) + \df{1}{2}J \circ du_i \circ k_i
\\
& = \delbar\xi_i.
\end{align*}
Thus $\ker(\hat{D}_i^V)$ is the set of holomorphic functions $(\Sigma_i,\partial\Sigma_i) \arr (V_{y_i},V_{y_i}^{(\R)})$, all of which are constant by Lemma~\ref{holFnRS}, and $\coker(\hat{D}_i^V)$ is just $V_{y_i}^{(\R)} \otimes_\R \coker(\delbar_{\Sigma_i})$.
\end{proof}
\end{lemma}

\begin{proposition} \label{fiber}
Fix a partition $\lambda=(g_1,\ldots,g_r,(g_{r+1},h_{r+1}),\ldots,(g_{r+q},h_{r+q}))$ of some topological type $(g,h)$. For $[u,\Sigma] \in \Nbar_\lambda$ with $p_i=u(z_i)$ the image of the $i^{\text{th}}$ node, the linearization satisfies
\begin{align*}
\ker(D) & = 0
\\
\coker(D) & = \left( \bigoplus\limits_{i=1}^r T_{p_i}M \otimes_\C H^{0,1}(\Sigma_i) \right) \oplus \left( \bigoplus\limits_{i=r+1}^{r+q} T_{p_i}L \otimes_\R H^{0,1}(\Sigma_i) \right).
\end{align*}
\begin{proof}
We have computed most of the data in Lemmas~\ref{ker0}, \ref{kerClosed}, and \ref{kerOpen}. All that remains is to understand how the operators on each component glue together to form $D$. This follows from an argument using long exact sequences; see Proposition~\ref{nodalLin}.
\end{proof}
\end{proposition}

\begin{proposition} \label{ob}
Let $\N$ be the moduli space of curves of type $(g,h)$ centered around $(u_0,\Sigma_0)$. Let $\Lambda$ be the set of partitions of $(g,h)$. For $\lambda \in \Lambda$, the obstruction bundle over $\Nbar_\lambda$ is
\[
Ob_\lambda = \left( \bigoplus\limits_{i=1}^r u_0^*TM \boxtimes_\C \E_{g_i}^* \right) \oplus \left( \bigoplus\limits_{i=r+1}^{r+q} u_0^*TL \boxtimes_\R \E_{(g_i,h_i)}^* \right)
\]
(where $\E_{g_i}$ and $\E_{g_i,h_i}$ are the Hodge bundles for genus $g_i$ and type $(g_i,h_i)$ curves respectively).

The obstruction bundle over $\Nbar$ is constructed by identifying fibers along intersections using Lemma~\ref{collisionModuli}. 
For collisions of closed ghosts $1 \leq i<i' \leq r$, we identify
\[
\left( u_0^*TM \boxtimes_\C \E_{g_i}^* \right) \oplus \left( u_0^*TM \boxtimes_\C \E_{g_{i'}}^* \right) \cong u_0^*TM \boxtimes_\C \E_{g_i+g_{i'}}^*.
\]
For collisions of open ghosts $r+1 \leq i<i' \leq r+q$, we identify 
\[
\left( u_0^*TL \boxtimes_\R \E_{(g_i,h_i)}^* \right) \oplus \left( u_0^*TL \boxtimes_\R \E_{(g_{i'},h_{i'})}^* \right) \cong u_0^*TL \boxtimes_\R \E_{(g_i+g_{i'},h_i+h_{i'}-1)}^*.
\]
For interior ghosts which approach $\partial\Sigma_0$, we identify
\[
u_0^*TM|_{\partial\Sigma_0} \boxtimes_\C \E_{g_i}^* \cong u_0^*TL \boxtimes_\R \E_{(g_i,1)}^*.
\]
\begin{proof}
By Proposition~\ref{fiber}, the kernel of the linearization at any map $u \in \N$ is zero. It follows that there exist bundles $Ob_\lambda \arr \Nbar_\lambda$ such that the fiber over any given map is precisely the cokernel of the linearization. These fibers were also computed in Proposition~\ref{fiber}. It is clear how the fibers fit together over any given cell, so all that remains is to understand cell intersections.

Fix a partition $\lambda=(g_1,\ldots,g_r,(g_{r+1},h_{r+1}),\ldots,(g_{r+q},h_{r+q}))$. 
There are three basic intersection types:
\begin{enumerate}[(I)]
\item a collision of two closed ghosts of genus $g_i$ and $g_{i'}$ produces a closed ghost with genus $g_i+g_{i'}$,
\item a collision of two open ghosts of toplogical type $(g_i,h_i)$ and $(g_{i'},h_{i'})$ produces an open ghost of type $(g_i~+~g_{i'}, h_i~+~h_{i'}~-~1)$, or
\item a closed ghost of genus $g_i$ approaches $\partial\Sigma_0$ to produce an open ghost of type $(g_i,1)$.
\end{enumerate}
In order to understand how to identify fibers along any cell intersections, we only need to understand identifications corresponding to these three types of cell intersections (see Section~\ref{baseS}).

We begin with the collision of two closed ghosts. Suppose that $(u,\Sigma) \in \Nbar_\lambda$ is a map where two closed ghosts collide. Assume without loss of generality that the colliding ghosts are $\Sigma_{r-1}$ and $\Sigma_r$, so $z_{r-1}=z_r$. Let
\[
\lambda' = (g_1,\ldots,g_{r-2},g_{r-1}+g_r,(g_{r+1},h_{r+1}),\ldots,(g_{r+q},h_{r+q})).
\]
Then there is a map $(v,\Sigma') \in \Nbar_{\lambda'}$ which represents the same curve. This means 
\[
\Sigma' = \Sigma_0 \cup \ldots \cup \Sigma_{r-2} \cup \Sigma_{r-1}' \cup \Sigma_{r+1} \cup \ldots \cup \Sigma_{r+q},
\]
where $\Sigma_i$ is attached to $\Sigma_0$ at $z_i$ and $\Sigma_{r-1}'$ is attached to $\Sigma_0$ at $z_{r-1}=z_r$. Moreover, $u$ and $v$ are identical along all the components they have in common, $(\Sigma_{r-1},\Sigma_r)$ is identified with $\Sigma_{r-1}'$ under the inclusion given by Lemma~\ref{collisionModuli}, and $u(\Sigma_{r-1})=u(\Sigma_r)=v(\Sigma_{r-1}')=p_{r-1} \in M$.

It is clear that in order to identify the fibers $\coker(D_u)$ and $\coker(D_v)$, it is sufficient to check
\[
\coker(D_{v,r-1}) \cong \coker(D_{u,r-1}) \oplus \coker(D_{u,r}).
\]
Using Lemma~\ref{kerClosed}, we identify 
\begin{align*}
\coker(D_{u,r-1}) & \cong T_{p_{r-1}}M \otimes_\C H^{0,1}(\Sigma_{r-1})
\\
\coker(D_{u,r}) & \cong T_{p_{r-1}}M \otimes_\C H^{0,1}(\Sigma_r)
\\
\coker(D_{v,r-1}) & \cong T_{p_{r-1}}M \otimes_\C H^{0,1}(\Sigma_{r-1}').
\end{align*}
All that remains is to apply Lemma~\ref{collisionModuli}.

Now we proceed to the collision of two open ghosts. Suppose that $(u,\Sigma) \in \Nbar_\lambda$ is a map where two open ghosts collide. Assume without loss of generality that the colliding ghosts are $\Sigma_{r+q-1}$ and $\Sigma_{r+q}$, so $z_{r+q-1}=z_{r+q}$. Let
\[
\lambda' = (g_1,\ldots,g_r,(g_{r+1},h_{r+1}),\ldots,(g_{r+q-2},h_{r+q-2}),(g_{r+q-1}+g_{r+q},h_{r+q-1}+h_{r+q}-1)).
\]
Then there is a map $(v,\Sigma') \in \Nbar_{\lambda'}$ which represents the same curve. This means 
\[
\Sigma' = \Sigma_0 \cup \ldots \cup \Sigma_{r+q-2} \cup \Sigma_{r+q-1}',
\]
where $\Sigma_i$ is attached to $\Sigma_0$ at $z_i$ and $\Sigma_{r+q-1}'$ is attached to $\Sigma_0$ at $z_{r+q-1}=z_{r+q}$. Moreover, $u$ and $v$ are identical along all the components they have in common, $(\Sigma_{r+q-1},\Sigma_{r+q})$ is identified with $\Sigma_{r+q-1}'$ under the inclusion given by Lemma~\ref{collisionModuli}, and $u(\Sigma_{r+q-1})=u(\Sigma_{r+q})=v(\Sigma_{r+q}')=p_{r+q-1} \in M$.

It is clear that in order to identify the fibers $\coker(D_u)$ and $\coker(D_v)$, it is sufficient to show
\[
\coker(D_{v,r+q-1}) \cong \coker(D_{u,rq-+1}) \oplus \coker(D_{u,r+q}).
\]
Using Lemma~\ref{kerClosed}, we identify 
\begin{align*}
\coker(D_{u,r+q-1}) & \cong T_{p_{r+q-1}}L \otimes_\R H^{0,1}(\Sigma_{r+q-1})
\\
\coker(D_{u,r+q}) & \cong T_{p_{r+q-1}}L \otimes_\R H^{0,1}(\Sigma_{r+q})
\\
\coker(D_{v,r+q-1}) & \cong T_{p_{r+q-1}}L \otimes_\R H^{0,1}(\Sigma_{r+q-1}')
.
\end{align*}
All that remains is to apply Lemma~\ref{collisionModuli}.

Finally, we examine what happens when a closed ghost approaches $\partial\Sigma_0$. Take $(u,\Sigma) \in \Nbar_\lambda$ and assume without loss of generality that $z_r \in \partial\Sigma_0$. Then there is a map $(v,\Sigma') \in \Nbar_{\lambda'}$  which represents the same curve, where
\[
\lambda'=(g_1,\ldots,g_{r-1},(g_r,1),(g_{r+1},h_{r+1}),\ldots,(g_{r+q},h_{r+q})).
\]
In particular, $(u,\Sigma)$ and $(v,\Sigma')$ are identical except for the $r^{\text{th}}$ ghost branch. The $r^{\text{th}}$ branch $\Sigma_r \in \Mbar_{g_r,1}$ of $\Sigma$ is identified with the $r^{\text{th}}$ branch $\Sigma_r' \in \Mbar_{(g_r,1),0,1}$ of $\Sigma'$ under the inclusion given by Lemma~\ref{collisionModuli}. Moreover, these branches are attached at the same point $z_r \in \Sigma_0$, and they are both sent to some point $p_r=u_0(z_r) \in L$.

It is evident that the only factor of the obstruction bundle which may differ over $u$ versus $v$ is the cokernel of the $r^{\text{th}}$ linearization:
\begin{align*}
D_{u,r}: & \Gamma(\Sigma_r;\Sigma_r \times T_{p_r}M) \arr \Omega^{0,1}(\Sigma_r,\Sigma_r \times T_{p_r}M)
\\
D_{v,r}: & \Gamma(\Sigma_r',\partial\Sigma_r';\Sigma_r' \times T_{p_r}M,\partial\Sigma_r' \times T_{p_r}L) \arr \Omega^{0,1}(\Sigma_r',\Sigma_r' \times T_{p_r}M).
\end{align*}
Since the bundles are trivial, we can identify
\begin{align*}
\Gamma(\Sigma_r;\Sigma_r \times T_{p_r}M) & \cong C^\oo(\Sigma_r, \C^3)
\\
\Omega^{0,1}(\Sigma_r,\Sigma_r \times T_{p_r}M) & \cong \Omega^{0,1}(\Sigma_r) \otimes_\C T_{p_r}M
\\
\Gamma(\Sigma_r',\partial\Sigma_r';\Sigma_r' \times T_{p_r}M,\partial\Sigma_r' \times T_{p_r}L) & \cong C^\oo(\Sigma_r',\partial\Sigma_r'; \C^3,\R^3)
\\
\Omega^{0,1}(\Sigma_r',\Sigma_r' \times T_{p_r}M) & \cong \Omega^{0,1}(\Sigma_r') \otimes_\R T_{p_r}M.
\end{align*}
Then the linearizations are identified with the $(0,0)$-Dolbeault operators for $\Sigma_r$ and $\Sigma_r'$, respectively. Since
\[
T_{p_r}M \cong \C \otimes_\R T_{p_r}L,
\]
all that remains is to apply Lemma~\ref{collisionModuli}:
\begin{align*}
T_{p_r}M \otimes_\C H^{0,1}(\Sigma_r) & \cong T_{p_r}L \otimes_\R \C \otimes_\C H^{0,1}(\Sigma_r)
\\
& \cong T_{p_r}L \otimes_\R H^{0,1}(\Sigma_r').
\end{align*}
\end{proof}
\end{proposition}

\section{Gluing Parameters} \label{glueS}

Fix a topological type $(g,h)$ of curve and let $\Nbar$ be the moduli space described in Section~\ref{baseS}. 
Our goal now is to determine the relationship between invariants of the bundle $Ob \arr \Nbar$ we built in Section~\ref{obS} and the contribution of curves in $\N$ to Gromov-Witten invariants. We perturb the equation $\delbar u = 0$ and count the resulting solutions. More precisely, fix some section $\nu$ of the bundle of $(0,1)$-forms over $\N$. Our goal is to count those $P \in \Nbar$ which perturb to a $t\nu$-holomorphic map for all small $t$.

We would like to view the solution space as the zero locus of a generic section of a bundle. Unfortunately, the rank of the obstruction bundle is too large: in the cell associated to a particular partition, we have $\rk_\R(Ob)=3\tilde{g}$ but $\dim_\R(\N)=3\tilde{g}-[(2r)+q]$ (where $\tilde{g}=2g+h-1$ is the genus of the double, $r$ is the number of closed ghost branches, and $q$ is the number of open ghost branches). Thus for each ghost branch we must introduce an extra line bundle (complex for interior nodes and real for boundary nodes) in order to resolve this difference.

This bundle will consist of gluing parameters, which give ways to smooth out nodes to yield new (non-holomorphic) curves. In this section we construct the bundle; in Section~\ref{lotS} we will examine its relationship to the contribution we wish to compute.

\begin{definition}
Fix a partition $\lambda=(g_1,\ldots,g_r,(g_{r+1},h_{r+1}),\ldots,(g_{r+q},h_{r+q}))$ of some topological type $(g,h)$.

For $1 \leq i \leq r$, the \emph{bundle of (complex) gluing parameters for the $i^{\text{th}}$ node} is $\L_i = T\Sigma_0 \boxtimes_\C \mathcal{T}_{g_i,1}$, where $\mathcal{T}_{g_i,1}$ is the relative tangent bundle over $\Mbar_{g_i,1}$.

For $r+1\leq i \leq r+q$, set $\sigma_i=((g_i,h_i),0,(1,0,\ldots,0))$. The \emph{bundle of (real) gluing parameters for the $i^{\text{th}}$ node} is $\L_i = T\partial\Sigma_0 \boxtimes_\R \mathcal{T}_{\sigma_i}$, where $\mathcal{T}_{\sigma_i}$ is the relative tangent bundle over $\Mbar_{\sigma_i}$. A gluing parameter $\tau_0 \otimes_\R \tau_1 \in T_{z_i}\partial\Sigma_0 \otimes_\R T_{y_i}\partial\Sigma_i$ is \emph{positive} if $-j(\tau_0) \in T_z\Sigma_0$ and $j(\tau_1) \in T_y\Sigma_1$ are both inward pointing or both outward pointing.

The \emph{bundle of gluing parameters} over $\Nbar_\lambda$ is
\[
\L_\lambda = \bigoplus\limits_{i=1}^{r+q} \L_i.
\]

The \emph{bundle of gluing parameters} $\pi_\L:\L\arr\Nbar$ is obtained by attaching the bundles over the cells along intersections. The fiber $\L_i$ is identified with the normal direction along the moduli space of curves whose $i^{\text{th}}$ node has been smoothed, or perpendicular to the zero section if such a smoothed curve does not lie in $\Nbar$ (see Remark~\ref{howToGlue}).
\end{definition}

\begin{remark} \label{howToGlue}
We see that $\L_i$ is a complex line bundle whose fiber over $(u,\Sigma)$ is $T_{z_i}\Sigma_0 \otimes_\C T_{y_i}\Sigma_i$ when $i \leq r$ and a real line bundle whose fiber over $(u,\Sigma)$ is $T_{z_i}\partial\Sigma_0 \otimes_\R T_{y_i}\partial\Sigma_i$ when $i>r$.

Here we explain how to attach the bundles $\L_{\lambda}$ so that the total space of $\L \arr \Nbar$ has constant dimension equal to the rank of the obstruction bundle. Over any individual cell this dimension criterion is clearly satisfied; all that remains is to examine cell intersections.

Consider a set $A=\{\lambda_1,\ldots,\lambda_k\}$ of partitions and an intersection $\Nbar_A= \bigcap_{\lambda \in A} \Nbar_{\lambda}$ of cells. We need to attach the bundles $\L_\lambda$ over this intersection. See Examples~\ref{glueEx11}, \ref{closedCollGlue}, and \ref{int4glue} for concrete interpretations.

First, we identify fibers of gluing parameters whenever nodes are identified. That is, all ghosts except those involved in the collision will be identified in a straightforward manner. After dividing $\bigoplus_{\lambda \in A} \L_{\lambda}$ by this equivalence, we are left with exactly one copy of $\C$ for each interior node and one copy of $\R$ for each boundary node in a generic curve in $\Nbar_A$ (note that this may include nodes outside of $\Sigma_0$, which are ignored throughout the rest of this paper).

Lemma~\ref{glue} gives instructions for identifying smoothing parameters with maps. A complex gluing parameter can be used to smooth an interior node, and a real gluing parameter can be used to smooth a boundary node. If smoothing the $i^{\text{th}}$ node yields a curve in $\bigcap_{\lambda \in A'} \Nbar_{\lambda}$ for $A' \subset A$, we attach the fiber of $\L_i$ along the normal bundle to $\Nbar_A$ in $\Nbar_{A'}$. Otherwise, we leave it as a fiber of $\L$. Since every direction perpendicular to $\Nbar_A$ in $\Nbar$ is obtained by smoothing at least one node, we see that a neighborhood of $\Nbar_A$ in $\L$ has dimension $\rk_\R(Ob)$.
\end{remark}

\begin{example} \label{glueEx11}
We return to the case presented in Section~\ref{11s} to understand gluing parameters over intersections arising from collision with $\partial\Sigma_0$. The intersection $\Nbar_{1,1} \cap \Nbar_\sigma$ has real dimension three. Fix $(u,\Sigma) \in \Nbar_{1,1} \cap \Nbar_\sigma$ with ghost $\rho_{1}([\Sigma_{1,1},y_{1,1}])=[\Sigma_\sigma,y_\sigma]$ attached at $z \in \Sigma_0$. Then the direct sum of fibers is
\[
(\L_{1,1})_u \oplus (\L_{\sigma})_u = \C_u \oplus \R_u,
\]
where
\begin{align*}
\C_u & = T_z\Sigma_0 \otimes_\C T_{y_{1,1}}\Sigma_{1,1}
\\
\R_u & = T_z\partial\Sigma_0 \otimes_\R T_{y_\sigma}\partial\Sigma_\sigma.
\end{align*}
We can split
\[
T\Nbar|_{\Nbar_{1,1} \cap \Nbar_\sigma} \cong T(\Nbar_{1,1} \cap \Nbar_\sigma) \oplus V_{1,1} \oplus V_\sigma,
\]
where $V_{1,1}$ is the normal bundle to $\Nbar_{1,1} \cap \Nbar_\sigma$ in $\Nbar_{1,1}$ and $V_\Sigma$ is the normal bundle to $\Nbar_{1,1} \cap \Nbar_\sigma$ in $\Nbar_\sigma$. Observe that $\dim_\R(V_{1,1})=1$ and $\dim_\C(V_\sigma)=1$; we wish to identify these bundles with $\R_u$ and $\C_u$, respectively.

Lemma~\ref{glue11} gives instructions for identifying smoothing parameters with maps. A complex gluing parameter $\tau \in \C_u$ can be used to smooth the interior node of $\Sigma$, and a real gluing parameter $\tau \in \R_u$ can be used to smooth the boundary node of $\Sigma$ (see Figure~\ref{pic11sigma}). Thus we can identify 
\begin{align*}
T_u\Nbar_{1,1} & \cong T_u(\Nbar_{1,1} \cap \Nbar_\sigma) \oplus \R_u
\\
T_u\Nbar_\sigma & \cong T_u(\Nbar_{1,1} \cap \Nbar_\sigma) \oplus \C_u.
\end{align*}
This process allows us to build the bundle $\L$ of gluing parameters over all of $\Nbar$. Although $\Nbar$ has two pieces of different dimensions, the total space of $\L$ has constant real dimension $6$ (which is also the real rank of the obstruction bundle).
\end{example}

\begin{example} \label{closedCollGlue}
In order to understand what happens to gluing parameters when ghosts collide, consider the case where $(g,h)=(4,1)$, $\lambda_1=(g_1=1,g_2=2,(g_3=1,h_3=1))$, and $\lambda_2=(g_1=3,(g_2=1,h_2=1))$. Label the nodes as in Figure~\ref{pic31two}. 
\begin{figure}[ht]
\centering
\begin{tikzpicture}

\def\r{1.5}
\def\h{1.2}
\def\w{0.8}


\coordinate (c1) at (0,0);
\disk[](c1)(\r)(\r);

\coordinate (Lnode1) at (-0.6*\r,0.8*\r);
\node [below right] at (Lnode1) {$a$};
\torus[](Lnode1)(30)(\h);

\coordinate (Rnode1) at (0.6*\r,0.8*\r);
\node [below left] at (Rnode1) {$b$};
\closedGhost[](Rnode1)(-30)(\w)(\h)(2);

\coordinate (Bnode1) at (\r,0);
\node [below] at (Bnode1) {$c$};
\openGhost[](Bnode1)(0.5*\r)(0.75*\r)(1)();
\begin{scope}[shift={(Bnode1)}]
	\addGenus[](0.5*\r,0.3*\r)(0)(0.25);
\end{scope}


\coordinate (c2) at (6,0);
\disk[](c2)(\r)(\r);

\coordinate (node2) at ($(c2)+(0,\r)$);
\node [below] at (node2) {$d$};
\closedGhost[](node2)(0)(\w)(\h)(3);

\coordinate (Bnode2) at ($(c2)+(\r,0)$);
\node [below] at (Bnode2) {$e$};
\openGhost[](Bnode2)(0.5*\r)(0.75*\r)(1)();
\begin{scope}[shift={(Bnode2)}]
	\addGenus[](0.5*\r,0.3*\r)(0)(0.25);
\end{scope}

\end{tikzpicture}
\caption{Two curves of type $(4,1)$.}
\label{pic31two}
\end{figure}
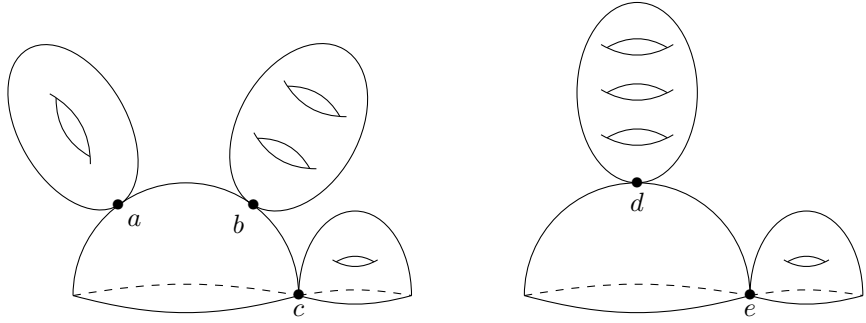

We examine the intersection $\Nbar_{\lambda_1} \cap \Nbar_{\lambda_2}$, which is precisely the set of curves $(u,\Sigma) \in \Nbar_{\lambda_1}$ where the two ghost curves $(\Sigma_1,y_1)$ and $(\Sigma_2,y_2)$ collide at $z_1=z_2$ (see Figure~\ref{pic31coll}).

\begin{figure}[ht]
\centering
\begin{tikzpicture}

\def\r{1.5}
\def\rb{1}
\def\h{1.2}
\def\w{0.8}

\coordinate (c1) at (0,0);
\disk[](c1)(\r)(\r);

\coordinate (c2) at (0,\r+\rb);
\coordinate (node2) at ($(c2)+(0,-\rb)$);
\node at (node2) {$\bullet$};
\node [below] at (node2) {$d$};
\sphere[](c2)(\rb)(\rb);

\coordinate (Lnode1) at ($(c2)+(-0.6*\rb,0.8*\rb)$);
\node [below right] at (Lnode1) {$a$};
\torus[](Lnode1)(30)(\h);

\coordinate (Rnode1) at ($(c2)+(0.6*\rb,0.8*\rb)$);
\node [below left] at (Rnode1) {$b$};
\closedGhost[](Rnode1)(-30)(\w)(\h)(2);

\coordinate (Bnode1) at (\r,0);
\node [below] at (Bnode1) {$c=e$};
\openGhost[](Bnode1)(0.5*\r)(0.75*\r)(1)();
\begin{scope}[shift={(Bnode1)}]
	\addGenus[](0.5*\r,0.3*\r)(0)(0.25);
\end{scope}

\end{tikzpicture}
\caption{A curve in $\Nbar_{(1,2,(1,1))} \cap \Nbar_{(3,(1,1))}$.}
\label{pic31coll}
\end{figure}

We begin with the direct sum of gluing parameters from $\Nbar_{\lambda_1}$ and $\Nbar_{\lambda_2}$:
\[
(\L_{\lambda_1})_u \oplus (\L_{\lambda_2})_u = (\C_a \oplus \C_b \oplus \R_c) \oplus (\C_d \oplus \R_e).
\]
The ghosts which are not involved in the collision are matched in a straightforward manner, so we first identify $\R_c$ with $\R_e$. Smoothing this node in a curve in $\Nbar_{\lambda_1} \cap \Nbar_{\lambda_2}$ does not yield a curve in $\Nbar$, so this fiber will be perpendicular to the zero section in $\L$.

Smoothing node $d$ gives a curve in $\Nbar_{\lambda_1}$, and since $\Nbar_{\lambda_1} \cap \Nbar_{\lambda_2}$ has real codimension $2$ in $\Nbar_{\lambda_1}$ it makes sense to identify $\C_d$ with the normal bundle $V_1$ (using Lemma~\ref{glue}). Similarly, curves such as $u$ exist in a real codimension $4$ subset of $\Nbar_{\lambda_2}$, so it makes sense to identify the gluing parameters $\C_a \oplus \C_b$ with the bundle $V_2$ normal to $\Nbar_{\lambda_1} \cap \Nbar_{\lambda_2}$ in $\Nbar_{\lambda_2}$. Now
\[
T_u\L \cong T_u(\Nbar_{\lambda_1} \cap \Nbar_{\lambda_2}) \oplus (V_1)_u \oplus (V_2)_u \oplus \R_c.
\]
We compute
\[
\begin{array}{rclcrcl}
\dim_\R(\Nbar_{\lambda_1}) & = & 19 & \qquad\qquad & \dim_\R(\Nbar_{\lambda_2}) & = & 21
\\
\dim_\R(\Nbar_{\lambda_1} \cap \Nbar_{\lambda_2}) & = & 17 && \codim_\R(\Nbar_{\lambda_1} \cap \Nbar_{\lambda_2},\Nbar_{\lambda_1}) & = & 2
\\
\codim_\R(\Nbar_{\lambda_1} \cap \Nbar_{\lambda_2},\Nbar_{\lambda_2}) & = & 4 && \dim_\R(\R_c) & = & 1.
\end{array}
\]
Therefore a neighborhood of $\Nbar_{\lambda_1} \cap \Nbar_{\lambda_2}$ in $\L$ has dimension
\begin{align*}
17+2+4+1 = 24 = \rk_\R(Ob).
\end{align*}
\end{example}

\begin{example} \label{int4glue}
We examine gluing parameters over the intersection $\Nbar_{\lambda_1} \cap \Nbar_{\lambda_2} \cap \Nbar_{\lambda_3} \cap \Nbar_{\lambda_4}$ described in Example~\ref{int4moduli}. Label the nodes on curves in $\Nbar_{\lambda_1}$, $\Nbar_{\lambda_2}$, $\Nbar_{\lambda_3}$, and $\Nbar_{\lambda_4}$ as in Figure~\ref{picFour}. The fibers of $\L_{\lambda_1}$, $\L_{\lambda_2}$, $\L_{\lambda_3}$, and $\L_{\lambda_4}$ are $\C_a \oplus \C_b$, $\C_{b'} \oplus \R_c$, $\C_{a'} \oplus \R_d$, and $\R_{c'} \oplus \R_{d'}$, respectively.

We first examine pairwise intersections. Over $\Nbar_{\lambda_1} \cap \Nbar_{\lambda_2}$, we have fibers $\C_a \oplus \C_b$ and $\C_{b'} \oplus \R_c$. Since nodes $b$ and $b'$ are not involved in the collision, we identify their fibers, leaving one gluing parameter per node: $\C_a \oplus \C_b \oplus \R_c$. Smoothing $b=b'$ does not yield a curve in $\Nbar$, so $\C_b=\C_{b'}$ will be perpendicular to the zero section in $\L$. Smoothing $a$ yields a curve in $\Nbar_{\lambda_2}$, so $\C_a$ is identified with the normal bundle to $\Nbar_{\lambda_1} \cap \Nbar_{\lambda_2}$ in $\Nbar_{\lambda_2}$. Smoothing $c$ yields a curve in $\Nbar_{\lambda_1}$, so $\R_c$ is identified with the normal bundle to $\Nbar_{\lambda_1} \cap \Nbar_{\lambda_2}$ in $\Nbar_{\lambda_1}$.

Similarly, over $\Nbar_{\lambda_1} \cap \Nbar_{\lambda_3}$ the fiber $\C_a=\C_{a'}$ is perpendicular to the zero section in $\L$, $\C_b$ is normal to $\Nbar_{\lambda_1} \cap \Nbar_{\lambda_3}$ in $\Nbar_{\lambda_3}$, and $\R_d$ is normal to $\Nbar_{\lambda_1} \cap \Nbar_{\lambda_3}$ in $\Nbar_{\lambda_1}$. Over $\Nbar_{\lambda_2} \cap \Nbar_{\lambda_4}$ the fiber $\R_c=\R_{c'}$ is perpendicular to the zero section in $\L$, $\C_{b'}$ is normal to $\Nbar_{\lambda_2} \cap \Nbar_{\lambda_4}$ in $\Nbar_{\lambda_4}$, and $\R_{d'}$ is normal to $\Nbar_{\lambda_2} \cap \Nbar_{\lambda_4}$ in $\Nbar_{\lambda_2}$. Over $\Nbar_{\lambda_3} \cap \Nbar_{\lambda_4}$ the fiber $\R_d=\R_{d'}$ is perpendicular to the zero section in $\L$, $\C_{a'}$ is normal to $\Nbar_{\lambda_3} \cap \Nbar_{\lambda_4}$ in $\Nbar_{\lambda_4}$, and $\R_{c'}$ is normal to $\Nbar_{\lambda_3} \cap \Nbar_{\lambda_4}$ in $\Nbar_{\lambda_3}$.

We now proceed to the intersection of all four cells. 
We begin with the direct sum of gluing parameters the four cells:
\[
(\C_a \oplus \C_b) \oplus (\C_{b'} \oplus \R_c) \oplus (\C_{a'} \oplus \R_d) \oplus (\R_{c'} \oplus \R_{d'}).
\]
We first identify fibers for nodes which are matched: $\C_a=\C_{a'}$, $\C_b=\C_{b'}$, $\R_c=\R_{c'}$, and $\R_d=\R_{d'}$. These four fiber directions are identified with the normal bundles to $\Nbar_{\lambda_1} \cap \Nbar_{\lambda_2}  \cap \Nbar_{\lambda_3} \cap \Nbar_{\lambda_4}$ in $\Nbar_{\lambda_2} \cap \Nbar_{\lambda_4}$, $\Nbar_{\lambda_3} \cap \Nbar_{\lambda_4}$, $\Nbar_{\lambda_1} \cap \Nbar_{\lambda_3}$, and $\Nbar_{\lambda_1} \cap \Nbar_{\lambda_2}$, respectively (see Table~\ref{table4int}). These are the only possible directions of movement in $\Nbar$, so a neighborhood of $\Nbar_{\lambda_1} \cap \Nbar_{\lambda_2}  \cap \Nbar_{\lambda_3} \cap \Nbar_{\lambda_4}$ in $\L$ has dimension
\[
\dim_\R(\Nbar_{\lambda_1} \cap \Nbar_{\lambda_2}  \cap \Nbar_{\lambda_3} \cap \Nbar_{\lambda_4})+2+2+1+1=18=\rk_\R(Ob).
\]
\end{example}

\begin{remark}
All cell intersections can be seen as straightforward ghost collisions in the complex double of the curve. For example, if an interior ghost collides with a boundary ghost, in the complex double we see this interior ghost collide simultaneously with its complex conjugate and a ghost attached along the real locus. Therefore, the principles explained in Example~\ref{int4glue} can be extended to all non-basic intersections.
\end{remark}

While smoothing a given curve, we choose some small constant $R_0>0$ which satisfies all the hypotheses for gluing in \cite{dw}. We also add the hypothesis that $R_0$ is small enough to guarantee that the ball of radius $4R_0$ around any interior node $z \in \Sigma_0\setminus\partial\Sigma_0$ does not intersect $\partial\Sigma_0$.

\begin{lemma} \label{glue}
Fix a partition $\lambda=(g_1,\ldots,g_r,(g_{r+1},h_{r+1}),\ldots,(g_{r+q},h_{r+q}))$ of some topological type $(g,h)$. For $(u,\Sigma)$ in the top stratum of $\Nbar_\lambda$, fix an element $\tau=(\tau_1,\ldots,\tau_{r+q})$ of the fiber $\L_{(u,\Sigma)}$ and assume
\begin{enumerate}[(i)]
\item $|\tau|<R_0$, and
\item $\tau_i$ is positive for $r+1 \leq i \leq r+q$.
\end{enumerate}
Then $\tau$ yields a Riemann surface $(\Sigma_\tau,j_\tau)$ and a smooth map $\tilde{u}_\tau:(\Sigma_\tau,\partial\Sigma_\tau) \arr (M,L)$ such that $\norm{\delbar(\tilde{u}_\tau)}_{L^p}$ is small in the sense of Proposition~5.8 of \cite{dw}.
\begin{proof}
For $1 \leq i \leq r$, we can use $\tau_i$ to smooth the node $z_i \sim y_i$ as in \cite{dw}. Because the domain and map are only altered near the node, the analysis in Sections~4.2 and 5.2 of \cite{dw} still applies.

However, we must take care when the node sits in $\partial\Sigma_0$. We may apply the results of \cite{dw} only to the double of the curve in the case of a boundary node.

Fix $(u,\Sigma) \in \Nbar_\lambda$ and let $(\Sigma^{(\C)},c,\pi)$ be the complex double of $\Sigma$ (see Definition~\ref{cplxDouble}). Choose a metric on $\Sigma^{(\C)}$ so that the fixed locus of the involution $c$ is totally geodesic. Then we can write
\[
\Sigma^{(\C)} = \left. \left( \bigsqcup\limits_{i=0}^{r+q} \Sigma_i^{(\C)} \right) \middle/ \left( z_i \sim y_i \right) \right. .
\]

We must analyze smoothings of $\Sigma^{(\C)}$ to understand smoothings of $\Sigma$. We focus here on the $i^{\text{th}}$ node for some $i>r$. Gluing parameters for the node $z_i \sim y_i$ in $\Sigma^{(\C)}$ are (small) elements of $T_{z_i}\Sigma_0^{(\C)} \otimes_\C T_{y_i}\Sigma_i^{(\C)}$. We smooth the $i^{\text{th}}$ node of $\Sigma^{(\C)}$ by removing small neighborhoods of $z_i$ and $y_i$ from $\Sigma_0^{\C}$ and $\Sigma_i^{\C}$, respectively, and then identifying small collars $A_{z_i}$ and $A_{y_i}$ around these removed neighborhoods via a map $\iota_{\tau_i}$ (see Figure~\ref{collars}).
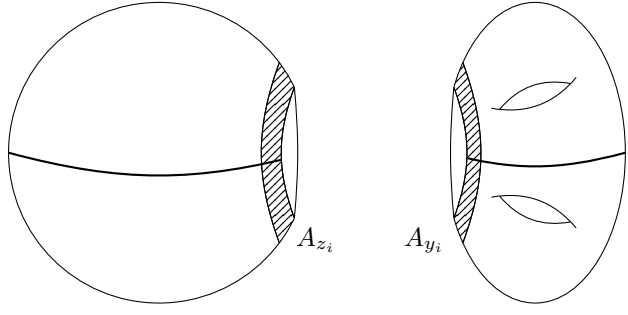
\begin{figure}[ht]
\centering
\begin{tikzpicture}

\def\r{2}
\def\w{0.6*\r}

\coordinate (sigma0) at (-1,0);

\begin{scope}[shift=(sigma0)]
	\draw (0,0) circle (\r);
	\draw [thick, bend right = 15] (-\r,0) to (\r,0);
\end{scope}

\begin{scope}[shift=(sigma0)]
	\draw [white, fill=white] (1.8,-0.872) to [bend left=20] (1.8,0.872) -- (2.1,1) -- (2.1,-1) -- cycle;
	
	\draw (1.8,-0.872) to [bend right=5] (1.8,0.872);
	\draw (1.8,-0.872) to [bend left=20] (1.8,0.872);
	\draw (1.6,-1.2) to [bend left=20] (1.6,1.2);
\end{scope}
\begin{scope}[shift=(sigma0)]
	\clip (sigma0) circle (\r);
	\draw [pattern = north east lines] (1.8,-0.872) to [bend left=20] (1.8,0.872) to [bend right=50] (1.6,1.2) to [bend right=20] (1.6,-1.2) to [bend right=50] cycle;
\end{scope}
\node [below right] at ($(sigma0)+(1.7,-0.872)$) {$A_{z_i}$};

\coordinate (sigma1) at (4,0);
\begin{scope}[shift=(sigma1)]
	\draw (0,0) ellipse ({\w} and {\r});
	\draw [thick, bend right = 15] (-\w,0) to (\w,0);
	\addGenus[](0,0.375*\r)(20)(0.5);
	\addGenus[](0,-0.375*\r)(160)(0.5);
\end{scope}

\begin{scope}[shift={(-\w+0.2,0)}]
    \draw [white, fill=white] (3.92,-0.872) to [bend right=20] (3.92,0.872) -- (3.7,1) -- (3.7,-1) -- cycle;
    
    \draw (3.92,-0.872) to [bend left=5] (3.92,0.872);
    \draw (3.92,-0.872) to [bend right=20] (3.92,0.872);
    \draw (4.04,-1.2) to [bend right=20] (4.04,1.2);
\end{scope}
\begin{scope}[shift={(-\w+0.2,0)}]
    \node [below left] at (3.92,-0.872) {$A_{y_i}$};
    \clip (sigma1) ellipse ({\w} and {\r});
    \draw [pattern = north east lines] (3.92,-0.872) to [bend right=20] (3.92,0.872) to [bend left=50] (4.04,1.2) to [bend left=20] (4.04,-1.2) to [bend left=50] cycle;
    \node [below left] at (3.92,-0.872) {$A_{y_i}$};
\end{scope}

\end{tikzpicture}
\caption{Collars near nodes.}
\label{collars}
\end{figure}

If $\tau_i=\tau_{z_i} \otimes_\C \tau_{y_i} \in T_{z_i}\Sigma_0 \otimes_\C T_{y_i}\Sigma_i$, then
\[
v \otimes_\C (\exp_{y_i}^{-1} \circ \iota_{\tau_i} \circ \exp_{z_i}(v)) = \tau_{z_i} \otimes_\C \tau_{y_i}
\]
for all $v \in T_{z_i}\Sigma_0$. In particular, $\iota_{\tau_i}(\exp_{z_i}(t\tau_{z_i}))=\exp_{y_i}(\tfrac{1}{t}\tau_{y_i})$ and $\iota_{\tau_i}(\exp_{z_i}(-j(t\tau_{z_i})))=\exp_{y_i}(j(\tfrac{1}{t}\tau_{y_i}))$ (see Figure~\ref{iota}).
\begin{figure}[ht]
\centering
\begin{tikzpicture}

\def\d{1}
\def\w{1}
\def\r{3}
\def\h{2}

\draw (-\d,-1) to [bend right=10] (-\d,1);
\draw (-\d,-1) to [bend left=15] (-\d,1);
\draw [dashed] (-\d-\w,-1) to [bend right=10] (-\d-\w,1);
\draw (-\d-\w,-1) to [bend left=15] (-\d-\w,1);
\draw (-\d-\w,-1) -- (-\d,-1);
\draw (-\d-\w,1) -- (-\d,1);
\node at (-\d-0.75*\w,0.4) {$\bullet$};
\node [left] at (-\d-\w,0.4) {$\exp_{z_i}(-j(s\tau_{z_i}))$};
\node at (-\d-0.25*\w,-0.8) {$\bullet$};
\node [below] at (-\d-0.25*\w,-0.9) {$\exp_{z_i}(t\tau_{z_i})$};
\node at (-\d-0.5*\w,-2) {$A_{z_i}$};

\coordinate (center0) at (-\r-3,-\h);
\coordinate (jt0) at ($(center0)+(-0.2,0.5)$);
\coordinate (t0) at ($(center0)+(-0.4,-0.4)$);
\draw [->] (center0) to (t0);
\draw [->] (center0) to (jt0);
\node [below] at (t0) {$\tau_{z_i}$};
\node [above] at (jt0) {$-j(\tau_{z_i})$};
\draw ($(center0)+(-1,-2)$) -- ($(center0)+(-1,1)$) -- ($(center0)+(1,2)$) -- ($(center0)+(1,-1)$) -- cycle;
\node at ($(center0)+(0,-2)$) {$T_{z_i}\Sigma_0$};

\coordinate (arrow0) at ($0.5*(center0)+0.5*(-\d,0)$);
\draw [->] ($(arrow0)+(-0.75,-0.25)$) to [bend left] ($(arrow0)+(0.75,0.25)$);
\node [above left] at ($(arrow0)+(0,0.1)$) {$\exp_{z_i}$};

\draw [dashed] (\d,-1) to [bend right=10] (\d,1);
\draw (\d,-1) to [bend left=15] (\d,1);
\draw (\d+\w,-1) to [bend right=10] (\d+\w,1);
\draw (\d+\w,-1) to [bend left=15] (\d+\w,1);
\draw (\d+\w,-1) -- (\d,-1);
\draw (\d+\w,1) -- (\d,1);
\node at (\d+0.25*\w,0.4) {$\bullet$};
\node [right] at (\d+\w,0.4) {$\exp_{y_i}(j(\tfrac{1}{s}\tau_{y_i}))$};
\node at (\d+0.75*\w,-0.8) {$\bullet$};
\node [below] at (\d+0.75*\w,-0.9) {$\exp_{y_i}(\tfrac{1}{t}\tau_{y_i})$};
\node at (\d+0.5*\w,-2) {$A_{y_i}$};

\coordinate (center1) at (\r+3,-\h);
\coordinate (jt1) at ($(center1)+(0.2,0.5)$);
\coordinate (t1) at ($(center1)+(0.4,-0.4)$);
\draw [->] (center1) to (t1);
\draw [->] (center1) to (jt1);
\node [below] at (t1) {$\tau_{y_i}$};
\node [above] at (jt1) {$j(\tau_{y_i})$};
\draw ($(center1)+(1,-2)$) -- ($(center1)+(1,1)$) -- ($(center1)+(-1,2)$) -- ($(center1)+(-1,-1)$) -- cycle;
\node at ($(center1)+(0,-2)$) {$T_{y_i}\Sigma_i$};

\coordinate (arrow1) at ($0.5*(center1)+0.5*(\d,0)$);
\draw [->] ($(arrow1)+(0.75,-0.25)$) to [bend right] ($(arrow1)+(-0.75,0.25)$);
\node [above right] at ($(arrow1)+(0,0.1)$) {$\exp_{y_i}$};

\draw [<->] (-0.5*\d,0) -- (0.5*\d,0);
\node [above] at (0,0) {$\iota_{\tau_i}$};

\end{tikzpicture}
\caption{Identifying collars via $\tau_i$.}
\label{iota}
\end{figure}
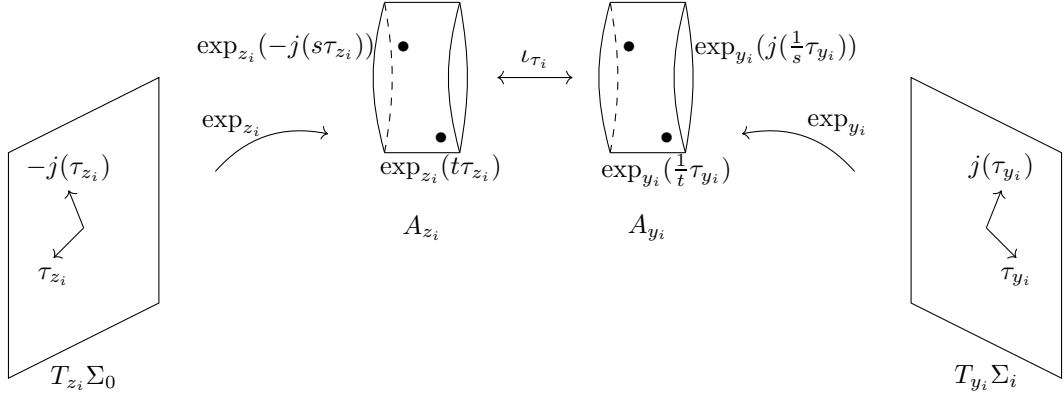

Now we must determine whether this smoothing of $\Sigma^{(\C)}$ yields a smoothing of $\Sigma$ near the $i^{\text{th}}$ node. The doubled curve $\Sigma^{(\C)}$ is equipped with an anti-holomorphic involution $c$ and a double cover $\pi:\Sigma^{(\C)}\arr\Sigma$. In order for the smoothing of $\Sigma^{(\C)}$ to yield a smoothing of $\Sigma$, the smoothing must respect these structures. That is, when we identify the $i^{\text{th}}$ collars in $\Sigma_0^{(\C)}$ and $\Sigma_i^{(\C)}$, the halves of the collars which lie in $\Sigma_0$ and $\Sigma_i$ must be identified. This occurs precisely when the gluing parameter $\tau_i$ lies in the positive half of the real locus of $T_{z_i}\Sigma_0^{(\C)} \otimes_\C T_{y_i}\Sigma_i^{(\C)}$.

Indeed, the smoothing $\Sigma_\tau^{\C}$ yields a smoothing of the $i^{\text{th}}$ node of $\Sigma$ precisely when
\[
\iota_{\tau_i}(A_{z_i} \cap \Sigma_0)=A_{y_i} \cap \Sigma_i.
\]
If $\tau_i=\tau_{z_i} \otimes_\C \tau_{y_i}$, we can assume without loss of generality that $\tau_{y_i}$ is tangent to the fixed locus $\text{Fix}(c)$ and that $j(\tau_{y_i})$ points inward along $\Sigma_i$. It follows that $\exp_{y_i}(\tfrac{1}{t}\tau_{y_i})$ must lie in the fixed locus and that $\exp_{y_i}(j(\tfrac{1}{t}\tau_{y_i}))$ must lie in $\Sigma_i$ (for appropriate values $t \in \R^+$). Since these two points are identified under $\iota_{\tau_i}$ with $\exp_{z_i}(t\tau_{z_i})$ and $\exp_{z_i}(-j(t\tau_{z_i}))$, respectively, we can smooth the $i^{\text{th}}$ node of $\Sigma$ via $\tau_i$ precisely when $\tau_{z_i}$ is also tangent to $\text{Fix}(c)$ and $-j(\tau_{z_i})$ is inward pointing along $\Sigma_0$. When we embed $\Sigma \arr \Sigma^{(\C)}$, we see that $\text{Fix}(c)$ is the image of $\partial\Sigma$, so a smoothing of $\Sigma_\tau^{(\C)}$ yields a smoothing of the $i^{\text{th}}$ node of $\Sigma$ if and only if the gluing parameter lies in the positive part of $T_{z_i}\partial\Sigma_0 \otimes_\R T_{y_i}\partial\Sigma_i$.

When $\tau$ satisfies the hypotheses of this lemma, let $(\Sigma_\tau,\partial\Sigma_\tau)$ be the smoothed domain. We define a smoothed map $u_\tau:(\Sigma_\tau,\partial\Sigma_\tau)\arr(M,L)$ precisely as in \cite{dw}. This map still sends $\partial\Sigma$ to $L$ because $L$ is totally geodesic. It is evident from the construction that the estimates computed in \cite{dw} still apply.
\end{proof}
\end{lemma}

\section{Leading Order Term} \label{lotS}

We have introduced gluing parameters in order to reconcile the difference between the rank of the obstruction bundle and the dimension of the base (see Remark~\ref{dimFinZeros}). In order to relate these parameters to the obstruction bundle, we pull the obstruction bundle back over the map $\pi_{\L}:\L\arr \Nbar$ and build a section of this new bundle.

\begin{lemma} \label{lot}
There is a section $\alpha:\L \arr \pi_{\L}^*Ob$ such that a curve $P \in \Nbar$ perturbs to a $t\nu$-holomorphic map if and only if there exists a gluing parameter $\tau$ (with $\tau_i$ positive for $i>r$) such that $a(P;\tau)=t\ov{\nu}_P$. Given a partition $\lambda=(g_1,\ldots,g_r,(g_{r+1},h_{r+1}),\ldots,(g_{r+q},h_{r+q}))$ of some topological type $(g,h)$, the restriction of $\alpha$ to the interior gluing parameters $\L_1 \oplus \ldots \oplus \L_r$ is the leading order term of the obstruction map constructed in Section~5.7 of \cite{dw}.
\begin{proof}
Fix a partition $\lambda=(g_1,\ldots,g_r,(g_{r+1},h_{r+1}),\ldots,(g_{r+q},h_{r+q}))$ and a map $(\Sigma,u) \in \Nbar_\lambda$. We define the map $\alpha$ on each factor by the relation \footnote{The formula given in Lemma~5.46 of \cite{dw} includes a map $u_{\sigma,\tau_1,0;J;\kappa}$. The parameter $\sigma$ represents a variation in the domain and $\kappa$ an element of the kernel of the linearization--- together, these correspond to our choice of map in $\Nbar_\lambda$. The gluing parameter $(\tau_1,0)$ in \cite{dw} is a gluing parameter for nodes internal to ghost branches (our $\tau$ corresponds to their $\tau_2$). We ignore the $J$ parameter because we do not vary the almost complex structure.}
\[
\langle \alpha_{\lambda,i}(\Sigma,u;\tau), \zeta_i \otimes_\C v_i \rangle_{L^2} = \langle (\ov{\zeta_i} \otimes_\C d_{y_i}u)(\tau_{i,1}), v_i \rangle
\]
for $i \leq r$ and 
\[
\langle \alpha_{\lambda,i}(\Sigma,u;\tau), \zeta_i \otimes_\R v_i \rangle_{L^2} = \langle (\ov{\zeta_i} \otimes_\R d_{y_i}(u|_{\partial\Sigma}))(\tau_{i,1}), v_i \rangle
\]
for $i>r$. Set
\[
\alpha_\lambda = \bigoplus\limits_{i=1}^{r+q} \alpha_{\lambda,i}.
\]
This matches the leading term of the Kuranishi map constructed in Section~5.6 of \cite{dw}. With minor modifications, the analysis of this map in Sections~5.6 and 5.7 of \cite{dw} applies to curves with boundary in the context of Ruan-Tian perturbations.
\end{proof}
\end{lemma}

\section{Calculation} \label{calcS}

To compute the contribution $C(g,h)$ of degree one covers of the main component $(\Sigma_0,u_0)$ to type $(g,h)$ Gromov-Witten invariants, we need to count those maps which can be perturbed in the sense of \cite{rt}. More specifically, we fix a generic section $\nu$ of the bundle $\mathcal{E}$ of $(0,1)$-forms over $\Nbar$ and let $\ov{\nu}$ be its projection to $Ob$. The contribution of $\Nbar$ to $C(g,h)$ is the number of maps $P=(\Sigma,u) \in \Nbar$ for which there exists a perturbation $P_\tau$ which satisfies $\delbar_J(P_\tau)=t\ov{\nu}(P_\tau)$, for all small $t>0$.

Let $\alpha$ be the (linear) section of $\pi_L^*Ob$ from Lemma~\ref{lot}, and let $Ob^F$ be the complement of its image $\im(\pi_L \circ \alpha)$ in Ob. 
\begin{equation}
\begin{tikzcd}
\Omega^{0,1}(\Nbar) \arrow[r, "\pi_{Ob}"] & Ob \arrow[d] & \pi_{\L}^*Ob \arrow[l, "\pi_{\L}"] \arrow[d]
\\
& \Nbar \arrow[ul, bend left, dashed, "\nu"] \arrow[u, bend left, dashed, "\ov{\nu}"] & \L \arrow[l, "\pi_{\L}"] \arrow[u, bend left, dashed, "\alpha"]
\end{tikzcd} \label{diagram}
\end{equation}
By Lemma~\ref{lot}, a map $P\in\Nbar$ perturbs precisely when it satisfies the following two conditions:
\begin{enumerate}[(i)]
\item $\ov{\nu}(P)$ lies in $\im(\pi_L \circ \alpha)$, and
\item the projection of $\ov{\nu}(P)$ to the factor $u_0^*TL \otimes_\R \E_{(g_i,h_i)}^*$ corresponding to any boundary ghost is positive.
\end{enumerate}
We will show that the positivity condition is irrelevant because we can build a generic non-vanishing section over factors with boundary ghosts. The condition $\ov{\nu}(P) \in \im(\pi_L \circ \alpha)$ is equivalent to $\proj_{Ob^F}(\ov{\nu})=0$. Since $\proj_{Ob^F}(\ov{\nu})$ is a generic section of $Ob^F$, the contribution is the count of the zero locus of a generic section of $Ob^F$.

\begin{remark} \label{dimFinZeros}
One function of gluing parameters is to fix the dimension problem. Let $\tilde{g}=2g+h-1$, so the obstruction bundle has real rank $3\tilde{g}$ over every cell. If a partition $\lambda$ has $r$ interior ghosts and $q$ boundary ghosts, then $\dim_\R(\Nbar_\lambda)=3\tilde{g}-(2r+q)$. However, since $2r+q$ is precisely the rank of the bundle of gluing parameters $\L_\lambda$, the bundle $Ob_\lambda^F$ has rank equal to the dimension of the base. Therefore, a generic section of this bundle has a finite number of zeros.

We stipulate that our section be non-vanishing along factors corresponding to boundary ghosts, which leaves only factors of the form $u_0^*TM \otimes_\C \E_{g_i}^*$. Since we will specify boundary conditions for $u_0^*TM \arr \Sigma_0$ and $\E_{g_i}^* \arr \Mbar_{g_i,1}$ has no codimension one boundary, we will eliminate any ambiguity arising from the choice of section.
\end{remark}

Proposition~\ref{calc} contains the bulk of the calculation; it applies when we count ordered constant branches. Corollary~\ref{perm} addresses the issue of permuting ghost branches.

\begin{remark}
The calculation in Proposition~\ref{calc} is essentially the same as that done in \cite{pand} and \cite{niuZinger}, but expressed differently. Because our curves have boundary, we write everything in terms of sections of the obstruction bundle rather than Chern classes. 
\end{remark}

\begin{proposition} \label{calc}
The contribution of degree one covers of $(\Sigma_0,u_0)$ with ordered ghost branches is
\[
\widetilde{C}(g,1) = \sum\limits_{g_1+\ldots+g_r=g} \prod\limits_{i=1}^r \left( \df{1}{2}\mu(N_0,N_0^{(\R)}) \cdot \alpha_{g_i} \right),
\]
where
\[
\alpha_{g_i} = \displaystyle\int_{\Mbar_{g_i,1}} c_{g_i}(\E_{g_i})c_{g_i-1}(\E_{g_i}) \left( \sumto{l=0}{g_i-1} (-1)^l c_l(\E_{g_i})\psi_{g,1}^{g_i-1-l} \right)
\]
for $\psi_{g,1}$ the first Chern class of the cotangent line over $\Mbar_{g,1}$. 
The contribution for $h>1$ is zero.
\begin{proof}
What is written below applies if we first pass to a smooth cover of each orbifold. Since each space of domains has a finite smooth cover, we can ignore the orbifold structure entirely.

We will construct a generic section $\rho$ of $Ob^F$ which is non-vanishing on each cell with boundary ghosts. In essence, this is possible because gluing a ghost along the boundary of the main component always yields a factor of $S^1$, which forces part of the obstruction bundle to be trivial. After eliminating the cells with boundary ghosts, the remaining contribution is more straightforward because all interior ghosts live in moduli spaces of closed curves (and in particular, the positivity criterion for gluing parameters does not apply).

The exclusion of boundary ghosts means that curves of type $(g,h)$ for $h>1$ cannot contribute. Moreover, the contribution arising from interior ghosts is very closely related to the corresponding contribution for degree one covers of a sphere. That is, we can separate the contribution coming from the factors $u_0^*TM \arr \Sigma_0$ and $\E_{g_i}^* \arr \Mbar_{g_i,1}$ of the obstruction bundle. The analysis of this second factor is slightly complicated, but it is essentially identical to that presented in \cite{pand}. On the other hand, the bundle $u_0^*TM \arr \Sigma_0$ is relatively easy to handle; we merely see a Maslov index instead of the Chern number we would get in the closed case.

We now proceed  the details of the proof. 
We computed the obstruction bundle $Ob$ in Proposition~\ref{ob} and bundle $\L$ of gluing parameters in Lemma~\ref{glue}. Let $\alpha$ be the leading order term from Lemma~\ref{lot}. Its image in $Ob_\lambda$ is 
\begin{equation}
Im(\pi_{\L_\lambda} \circ \alpha_\lambda) = \left( \bigoplus\limits_{i=1}^r T\Sigma_0 \boxtimes_\C F_{g_i}^\perp \right) \oplus \left( \bigoplus\limits_{i=r+1}^{r+q} T\partial\Sigma_0 \boxtimes_\R F_{g_i,h_i}^\perp \right), \label{imAlpha}
\end{equation}
where $F_{g_i}^\perp$ is the complex rank $1$ complement of the bundle generated by $\zeta_{y_i}=0$ in $\E_{g_i}^*$ and $F_{g_i,h_i}^\perp$ is the real rank $1$ complement of the bundle generated by $\zeta_{y_i}=0$ in $\E_{g_i,h_i}^*$. Then $Ob_\lambda^F$ is the complement of~(\ref{imAlpha}).

Whenever a partition $\lambda$ has ghosts along the boundary, the corresponding cell $\Nbar_\lambda$ has factors of the form $S^1 \times \Mbar_{\sigma_{\lambda,i}}$ for $\sigma_{\lambda,i}=((g_i,h_i),0,(1,0,\ldots,0))$. In general a moduli space of open curves, such as $\Mbar_{\sigma_{\lambda,i}}$, may have codimension one boundary (meaning that the zero count may vary from one section to the next). However, we will use $S^1$ to construct a non-vanishing section for such a factor, which eliminates any ambiguity in the zero count of a section.

We consider each factor of $Ob$ separately. First we decompose the tangent bundles $TM|_{\Sigma_0}$ and $TL|_{\partial\Sigma_0}$. We can split into directions tangent and normal to the curve. Because the normal bundle is a complex rank two bundle over a surface, we can split off a trivial line bundle. Therefore, we can write
\begin{align*}
TM|_{\Sigma_0} & = V_1 \oplus V_2 \oplus V_3
\\
TL|_{\partial\Sigma_0} & = V_1^{(\R)} \oplus V_2^{(\R)} \oplus V_3^{(\R)},
\end{align*}
where $V_1=T\Sigma_0$, $V_j^{(\R)}=V_j \cap TL$, and $(V_3,V_3^{(\R)})$ is a trivial bundle pair.

Pick generic sections $v_1$ of $V_1$, $v_2$ and $\tilde{v}_2$ of $V_2$, and $v_3$ of $V_3$ so that
\begin{enumerate}[(i)]
\item $v_1$, $v_2$, and $\tilde{v}_2$ have pairwise disjoint zero loci,
\item $v_3$ is non-vanishing,
\item the restriction of each section to $\partial\Sigma_0$ lands in the totally real sub-bundle and is non-vanishing as a section of that bundle.
\end{enumerate}
The first condition is generic for dimension reasons. It is possible to insist that the projections onto the real sub-bundles be non-vanishing because every (orientable) bundle over $\partial\Sigma_0 \cong S^1$ is trivial. Observe that the zero loci of $v_2$ and $\tilde{v}_2$ are Poincar\'{e} dual to $\tfrac12\mu(N_0,N_0^{(\R)})$.

In the case $q>0$, take any index $i>r$ corresponding to an open ghost. Then the $i^{\text{th}}$ factor of $Ob$ is
\[
u_0^*TL \boxtimes_\R \E_{g_i,h_i}^*.
\]
By Proposition~\ref{eSquared} it is possible to choose generic sections $\eta_{i,1},\eta_{i,2},\eta_{i,3}$ of $\E_{g_i,h_i}^*$ so that $\eta_{i,2}$ and $\eta_{i,3}$ have disjoint zero loci \footnote{Technically, choosing sections of duals of Hodge bundles is somewhat complicated because they must match along cell intersections. One way to resolve this issue would be recognize each domain modeled on $\lambda$ as an element of $\Mbar_{(g,h),0,\vec{0}}$ using arguments similar to those in Lemma~\ref{collisionModuli}.}. The section
\[
\rho_{\lambda,i} = (v_1|_{\partial\Sigma_0} \boxtimes_\R \proj_{F_{g_i,h_i}}(\eta_{i,1})) \oplus (v_2|_{\partial\Sigma_0} \boxtimes_\R \eta_{i,2}) \oplus (v_3|_{\partial\Sigma_0} \boxtimes_\R \eta_{i,3})
\]
of the $i^{\text{th}}$ factor of $Ob^F$ is non-vanishing. Picking any generic sections of the other factors of $Ob^F$ will yield a non-vanishing section of $Ob^F$ over $\Nbar_\lambda$. It follows that the contribution from $\Nbar_\lambda$ is zero whenever there are open ghosts.

We have now shown that the total contribution from $\Nbar$ is equal to the count of perturbable maps in $\Nbar$ with only closed ghosts. If $h>1$ then every partition of $(g,h)$ has at least one open ghost, so the contribution is zero except when $h=1$.

We are left with partitions of the form $\lambda=(g_1,\ldots,g_r)$; maps modeled on these partitions have only closed ghosts. But now the complicated part of the computation reduces to the case of closed invariants. We can choose generic sections $\eta_{i,1},\eta_{i,2},\eta_{i,3}$ of $E_{g_i}$ so that if $\zeta_{i,j}=\proj_{F_{g_j}}(\eta_{i,j})$ and $\zeta_{i,j}^\perp=\proj_{F_{g_j}^\perp}(\eta_{i,j})$, then
\begin{enumerate}[(i)]
\item $Z(\eta_{i,2}) \cap Z(\eta_{i,3}) = \emptyset$,
\item $Z(\zeta_{i,1}) \cap Z(\zeta_{i,2}^\perp) \cap Z(\eta_{i,3}) = \emptyset$, and
\item $Z(\zeta_{i,1}) \cap Z(\zeta_{i,2}) \cap Z(\eta_{i,3})$ is Poincar\'{e} dual to $\alpha_{g_i}$.
\end{enumerate}
The first two conditions are made possible by Proposition~\ref{eSquared}. The last condition is a result of the techniques of \cite{pand} and \cite{niuZinger} (this is where it is crucial to observe that all ghosts are closed, so existing techniques apply).

We assemble this data to build a generic section $\rho$ of $Ob^F$ whose restriction to the $\lambda$-cell is the following:
\[
\rho_{\lambda,i} = (v_1 \boxtimes_\R \zeta_{i,1}) \oplus (v_2 \boxtimes_\R \zeta_{i,2}) \oplus (\tilde{v}_2 \boxtimes_\R \zeta_{i,2}^\perp) \oplus (v_3 \boxtimes_\R \eta_{i,3}).
\]
Then the zero locus of $\rho$ is Poincar\'{e} dual to 
\[
\prod\limits_{i=1}^r \left( \df{1}{2}\mu(N_0,N_0^{(\R)}) \cdot \alpha_{g_i} \right).
\]
Finally, if we choose these sections carefully so that they agree along cell intersections, we can compute the entire contribution of $\Nbar$ by adding across partitions, as in \cite{pand} and \cite{niuZinger}.
\end{proof}
\end{proposition}

\begin{corollary} \label{perm}
Fix $g \geq 0$ and let $\Lambda_{c}$ be the set of partitions of $(g,1)$ with only closed ghosts. Then the contribution of degree one covers of a disk to Gromov-Witten invariants of type $(g,1)$ is
\[
C(g,1) = \sum\limits_{\lambda \in \Lambda_c} \df{1}{|\Aut(\lambda)|} \prod\limits_{i=1}^r \left( \df{1}{2}\mu(N_0,N_0^{(\R)}) \cdot \alpha_{g_i} \right).
\]
The generating function is
\[
\sumto{g=0}{\oo} C(g,1)t^{2g-1} = \left( \df{\sin(t/2)}{t/2} \right)^{-1}.
\]
\end{corollary}

\begin{corollary} \label{zero}
If $h>1$, then the contribution $C(g,h)$ is zero.
\end{corollary}

\begin{remark} \label{doubleCt}
Suppose that $L$ is the fixed locus of an anti-symplectic involution on $M$.  Then the Schwarz reflection principle allows us to double maps, bundles, and sections, as in Section~3.3.3 of \cite{katzLiu}. These doubled maps have no boundary, which allows us to make sense of the contribution using formal properties of characteristic classes. If we decompose $TM|_{\Sigma_0}$ as in the proof of Proposition~\ref{calc}, then the contribution of the $i^{\text{th}}$ factor is
\[
\df{1}{2}\mu(N_0,N_0^{(\R)}) \cdot \alpha_{g_i} = \#Z(v_2)\alpha_{g_i} = \df{1}{2}c_1(V_2^{(\C)})\alpha_{g_i} = \df{1}{2}c_1(N_0^{(\C)})\alpha_{g_i}
\]
(cf. \cite{niuZinger}).
\end{remark}

\section*{List of Symbols}

\renewcommand{\arraystretch}{1.5}
\def\colWidth{1.2in}

\begin{longtable}{l>{\raggedright}p{3.7in}l} \label{dict}
{\bf Symbol} & {\bf Meaning} & {\bf Location}
\\
$\alpha$ & leading order term of obstruction map & \parbox[t]{\colWidth}{Lemma~\ref{lot11} \\ Lemma~\ref{lot}}
\\
$\mathcal{B}$ & space of domains and smooth maps & Section~\ref{prelim}
\\
$\boxtimes$ & exterior tensor product of bundles & Definition~\ref{extTensor}
\\
$C$, $C_0$ & image of $u$, $u_0$ & 
\\
$C(g,h)$ & contribution of a $(0,1)$ curve to $(g,h)$ invariants & \parbox[t]{\colWidth}{Section~\ref{intro} \\ Proposition~\ref{calc11} \\ Proposition~\ref{calc}}
\\
$D$, $D^T$, $D^N$ & linearization of $\delbar_J$ and its tangent and normal components & Section~\ref{prelim}
\\
$\hat{D}$ & linearization of $\delbar_J$ with fixed domain &
\\
$\delbar_J$ &  & Section~\ref{prelim}
\\
$\delbar_S$ & $(0,0)$-Dolbeault operator for $S$ & Definition~\ref{hodge}
\\
$\mathcal{E}$ & bundle of $(0,1)$-forms & Section~\ref{prelim}
\\
$\E_g$, $\E_{(g,h)}$ & Hodge bundle & Definition~\ref{hodge}
\\
$\Gamma(\Sigma,\partial\Sigma;E,F)$ & vector fields on $\Sigma$ with values in $E$ which take values in $F$ along $\partial\Sigma$ & Section~\ref{prelim}
\\
$H^{0,1}(\Sigma)$ & $\Omega^{0,1}(\Sigma)/\im(\delbar_\Sigma)$ & Definition~\ref{hodge}
\\
$\L$ & bundle of gluing parameters & \parbox[t]{1.1in}{Subsection~\ref{glue11ss} \\ Section~\ref{glueS}}
\\
$\lambda$ & partition of a topological type & Definition~\ref{partNum}
\\
$\Mbar_{g,k}$ & moduli space of closed genus $g$ curves with $k$ marked points & Definition~\ref{modDomain}
\\
$\Mbar_{(g,h),n,\vec{m}}$ & moduli space of open genus $g$ curves with $h$ boundary components and $(n,\vec{m})$ marked points & Definition~\ref{modDomain}
\\
$\mu$ & Maslov class & 
\\
$\Nbar$ & moduli space of holomorphic maps with main component $(\Sigma_0,u_0)$ & \parbox[t]{1.1in}{Subsection~\ref{base11ss} \\ Section~\ref{baseS}}
\\
$\Nbar_\lambda$ & moduli space of holomorphic maps modeled on $\lambda$ & Definition~\ref{partition}
\\
$N_i$ & normal bundle to $T_i$ in $M$ & 
\\
$\nu$ & generic section of $\mathcal{E}$ & 
\\
$\ov{\nu}$ & projection of $\nu$ to $Ob$ & 
\\
$Ob$ & obstruction bundle & \parbox[t]{1.1in}{Subsection~\ref{ob11ss} \\ Section~\ref{obS}}
\\
$Ob^F$ & complement of $\im(\alpha)$ in $Ob$ & \parbox[t]{1.1in}{Subsection~\ref{calc11ss} \\ Section~\ref{calcS}}
\\
$p$, $p_i$ & image of node in $M$ &
\\
$q$ & number of open ghosts & 
\\
$r$ & number of closed ghosts & 
\\
$R_0$ & radius for gluing & \parbox[t]{1.1in}{Subsection~\ref{glue11ss} \\ Section~\ref{glueS}}
\\
$\Sigma^{(\C)}$ & complex double of $\Sigma$ & Definition~\ref{cplxDouble}
\\
$T_i$ & tangent space to $C$ over $\Sigma_i$ & 
\\
$\tau$ & gluing parameter & \parbox[t]{1.1in}{Subsection~\ref{glue11ss} \\ Section~\ref{glueS}}
\\
$\mathcal{T}$ & relative tangent bundle & Definition~\ref{relTan}
\\
$V^{(\R)}$ & the real locus of a vector bundle $V$ & 
\\
$y$, $y_i$ & point in a ghost branch which is attached to main component & 
\\
$z$, $z_i$ & point in main component at which a ghost branch is attached & 
\end{longtable}

\listoffigures

\printbibliography

@book{audinDamian,
      AUTHOR = {Audin, Mich\`ele and Damian, Mihai},
     TITLE = {Morse theory and {F}loer homology},
    SERIES = {Universitext},
      NOTE = {Translated from the 2010 French original by Reinie Ern\'{e}},
 PUBLISHER = {Springer, London; EDP Sciences, Les Ulis},
      YEAR = {2014},
     PAGES = {xiv-596},
      ISBN = {978-1-4471-5495-2},
   MRCLASS = {53-02 (53D40 58E05)},
  MRNUMBER = {3155456},
MRREVIEWER = {Sonja Hohloch},
       DOI = {10.1007/978-1-4471-5496-9},
       URL = {https://doi-org.stanford.idm.oclc.org/10.1007/978-1-4471-5496-9},
}

@misc{dw,
      title={Counting embedded curves in symplectic 6-manifolds}, 
      author={Aleksander Doan and Thomas Walpuski},
      year={2019},
      eprint={1910.12338},
      archivePrefix={arXiv},
      primaryClass={math.SG}
}

@article{ortn,
   title={The orientability problem in open Gromov–Witten theory},
   volume={17},
   ISSN={1465-3060},
   url={http://dx.doi.org/10.2140/gt.2013.17.2485},
   DOI={10.2140/gt.2013.17.2485},
   number={4},
   journal={Geometry \& Topology},
   publisher={Mathematical Sciences Publishers},
   author={Georgieva, Penka},
   year={2013},
   month={8},
   pages={2485–2512}
}

@misc{splitting,
      title={Splitting formulas for the local real Gromov-Witten invariants}, 
      author={Penka Georgieva and Eleny-Nicoleta Ionel},
      year={2021},
      eprint={2005.05928},
      archivePrefix={arXiv},
      primaryClass={math.SG},
      note={Preprint}
}

@book{harris,
	AUTHOR = {Harris, Joe and Morrison, Ian},
     TITLE = {Moduli of curves},
    SERIES = {Graduate Texts in Mathematics},
    VOLUME = {187},
 PUBLISHER = {Springer-Verlag, New York},
      YEAR = {1998},
     PAGES = {xiv-366},
      ISBN = {0-387-98438-0},
   MRCLASS = {14H10 (14-02 14D20 14D22)},
  MRNUMBER = {1631825},
MRREVIEWER = {R. F. Lax},
}

@article{katzLiu,
   title={Enumerative geometry of stable maps with Lagrangian boundary conditions and multiple covers of the disc},
   ISSN={1464-8997},
   url={http://dx.doi.org/10.2140/gtm.2006.8.1},
   DOI={10.2140/gtm.2006.8.1},
   journal={The interaction of finite-type and Gromov--Witten
invariants (BIRS 2003)},
   publisher={Mathematical Sciences Publishers},
   author={Katz, Sheldon and Liu, Chiu-Chu Melissa},
   year={2006},
   month={4}
}

@phdthesis{liu,
AUTHOR = {Liu, Chiu-Chu Melissa},
     TITLE = {Moduli of {J}-holomorphic curves with {L}agrangian boundary
              conditions},
      NOTE = {Thesis (Ph.D.)--Harvard University},
      SCHOOL = {Harvard University},
 PUBLISHER = {ProQuest LLC, Ann Arbor, MI},
      YEAR = {2002},
     PAGES = {151},
      ISBN = {978-0493-65757-8},
   MRCLASS = {Thesis},
  MRNUMBER = {2703393},
       URL =
              {http://gateway.proquest.com.stanford.idm.oclc.org/openurl?url_ver=Z39.88-2004&rft_val_fmt=info:ofi/fmt:kev:mtx:dissertation&res_dat=xri:pqdiss&rft_dat=xri:pqdiss:3051225},
}

@book{msBig,
      AUTHOR = {McDuff, Dusa and Salamon, Dietmar},
     TITLE = {{$J$}-holomorphic curves and symplectic topology},
    SERIES = {American Mathematical Society Colloquium Publications},
    VOLUME = {52},
   EDITION = {Second},
 PUBLISHER = {American Mathematical Society, Providence, RI},
      YEAR = {2012},
     PAGES = {xiv-726},
      ISBN = {978-0-8218-8746-2},
   MRCLASS = {53D45 (32Q65 53D35)},
  MRNUMBER = {2954391},
MRREVIEWER = {Mark Alan Branson},
}

@incollection {mumford,
    AUTHOR = {Mumford, David},
     TITLE = {Towards an enumerative geometry of the moduli space of curves},
 BOOKTITLE = {Arithmetic and geometry, {V}ol. {II}},
    SERIES = {Progr. Math.},
    VOLUME = {36},
     PAGES = {271--328},
 PUBLISHER = {Birkh\"{a}user Boston, Boston, MA},
      YEAR = {1983},
   MRCLASS = {14H10 (14C15)},
  MRNUMBER = {717614},
MRREVIEWER = {Werner Kleinert},
}

@article{niuZinger,
      AUTHOR = {Niu, Jingchen and Zinger, Aleksey},
     TITLE = {Lower bounds for enumerative counts of positive-genus real
              curves},
   JOURNAL = {Adv. Math.},
  FJOURNAL = {Advances in Mathematics},
    VOLUME = {339},
      YEAR = {2018},
     PAGES = {191--247},
      ISSN = {0001-8708},
   MRCLASS = {14N35 (53D45)},
  MRNUMBER = {3866896},
MRREVIEWER = {Hsian-Hua Tseng},
       DOI = {10.1016/j.aim.2018.09.024},
       URL = {https://doi-org.stanford.idm.oclc.org/10.1016/j.aim.2018.09.024},
}

@article{pand,
   title={Hodge Integrals and Degenerate Contributions},
   volume={208},
   ISSN={1432-0916},
   url={http://dx.doi.org/10.1007/s002200050766},
   DOI={10.1007/s002200050766},
   number={2},
   journal={Communications in Mathematical Physics},
   publisher={Springer Science and Business Media LLC},
   author={Pandharipande, R.},
   year={1999},
   month={12},
   pages={489–506}
}

@article{rt,
author = {Yongbin Ruan and Gang Tian},
title = {{A mathematical theory of quantum cohomology}},
volume = {42},
journal = {Journal of Differential Geometry},
number = {2},
publisher = {Lehigh University},
pages = {259 -- 367},
year = {1995},
doi = {10.4310/jdg/1214457234},
URL = {https://doi.org/10.4310/jdg/1214457234}
}

@book{wendl,
      AUTHOR = {Wendl, Chris},
     TITLE = {Holomorphic curves in low dimensions},
    SERIES = {Lecture Notes in Mathematics},
    VOLUME = {2216},
      NOTE = {From symplectic ruled surfaces to planar contact manifolds},
 PUBLISHER = {Springer, Cham},
      YEAR = {2018},
     PAGES = {xiii-292},
      ISBN = {978-3-319-91369-8},
   MRCLASS = {57R17 (32Q65)},
  MRNUMBER = {3821526},
MRREVIEWER = {Richard Keith Hind},
       DOI = {10.1007/978-3-319-91371-1},
       URL = {https://doi-org.stanford.idm.oclc.org/10.1007/978-3-319-91371-1},
}

\end{document}